\documentclass{article}
\usepackage{amsfonts, amsbsy, amssymb}
\usepackage{graphicx}

\begin{document}


\newcommand{\eps}{\varepsilon}
\newcommand{\seps}{\sqrt{\eps}}
\def\omeps{{\omega \over \seps}}
\def\omepsi{{\omega \over \eps}}

\newcommand{\x}{\boldsymbol x}
\newcommand{\s}{\boldsymbol s}
\newcommand{\io}{\boldsymbol \iota}
\newcommand{\p}{\boldsymbol p}
\newcommand{\q}{\boldsymbol q}
\newcommand{\g}{\boldsymbol g}
\newcommand{\ar}{\boldsymbol a}
\newcommand{\br}{\boldsymbol b}
\newcommand{\ur}{\boldsymbol u}
\newcommand{\ka}{\boldsymbol k}
\newcommand{\w}{\boldsymbol \omega}
\newcommand{\ps}{\boldsymbol \psi}
\newcommand{\xii}{\boldsymbol \xi}

\newcommand{\fb}{\mathfrak b}
\newcommand{\fB}{\mathfrak B}
\newcommand{\fD}{\mathfrak D}
\newcommand{\fN}{\mathfrak N}
\newcommand{\fW}{\mathfrak W}
\newcommand{\fS}{\mathfrak S}
\newcommand{\fI}{\mathfrak I}
\newcommand{\fri}{\mathfrak i}
\newcommand{\fp}{\mathfrak{p}}

\newcommand{\de}{{\scriptstyle \Delta}}
\newcommand{\te}{{\scriptstyle T}}
\newcommand{\ie}{{\scriptstyle I}}

\newcommand{\ld}{,\ldots,}
\newcommand{\spr}[2]{\langle#1,#2\rangle}
\newcommand{\ti}[1]{{}^{\mathfrak t}{#1}^{\mathtt{-1}}}
\def\={\stackrel{\rm def}{=}}

\def\B{{\mathbb B}}
\def\R{{\mathbb R}}
\def\Rn{{\mathbb R}^n}
\def\Rnp{{\mathbb R}^{n+1}}
\def\T{{\mathbb T}}
\def\N{{\mathbb N}}
\def\M{{\mathbb M}}
\def\C{{\mathbb C}}
\def\Z{{\mathbb Z}}
\def\Q{{\mathbb Q}}

\def\Dlw{D_{\lambda,\omega}}
\newcommand{\Diff}{\lambda D_s+\spr{\omega }{D_\varphi} }

\newcommand{\It}{{\cal I}\times\T^n}
\newcommand{\cyl}{\cal C}

\newcommand{\di}{{\fW}^n_{\tau,\gamma}}
\newcommand{\defp}{{\fD}_{\fp,\de}({\cal C})}
\newcommand{\defpit}{{\fD}_{\fp,\de}(\It)}
\newcommand{\defppi}{{\fD}_{\fp,\de}(\Pi\times\T^n)}

\newcommand{\bst}{{\fB}_\sigma(\T^n)}
\newcommand{\bstprime}{{\fB}_{\sigma'}(\T^n)}
\newcommand{\bstk}{{\fB}_{\kappa,\sigma}(T^*\T^n)}

\newcommand{\bsitb}{{\fB}_{\beta,\fp}({\cal I}\times\T^n)}
\newcommand{\bsitbj}{{\fB}_{\beta,\fp}^j({\cal I}\times\T^n)}
\newcommand{\bsit}{{\fB}_{\fp}({\cal I}\times\T^n)}
\newcommand{\bsitj}{{\fB}_{\fp}^j({\cal I}\times\T^n)}
\newcommand{\bsitprime}{{\fB}_{\fp'}({\cal I}\times\T^n)}

\newcommand{\bskit}{{\fB}_{\kappa,\mathfrak{p}}[T^*(\It)]}
\newcommand{\bskitprime}{{\fB}_{\kappa',\mathfrak{p}'}[T*(\It)]}

\newcommand{\bs}{{\fB}_{\mathfrak{p}}({\cal C})}
\newcommand{\bsone}{{\fB}^1_{\mathfrak{p}}({\cal C})}
\newcommand{\bstwo}{{\fB}^2_{\mathfrak{p}}({\cal C})}
\newcommand{\bsj}{{\fB}^j_{\mathfrak{p}}({\cal C})}

\newcommand{\bsprime}{{\fB}_{\mathfrak{p}'}({\cal C})}
\newcommand{\bsprimeb}{{\fB}_{\mathfrak{p}'}({\cal C_\beta})}
\newcommand{\bsoneprime}{{\fB}^1_{\mathfrak{p}'}({\cal C})}
\newcommand{\bstwoprime}{{\fB}^2_{\mathfrak{p}'}({\cal C})}
\newcommand{\bsjprime}{{\fB}^j_{\mathfrak{p}'}({\cal C})}

\newcommand{\bsm}{{\fB}^{-}_{\mathfrak{p}}({\cal C})}
\newcommand{\bsmprime}{{\fB}_{\mathfrak{p}'}^{-}({\cal C})}

\newcommand{\bsg}{{\fB}^{(-,1)}_{\mathfrak{p}}({\cal C})}
\newcommand{\bsgprime}{{\fB}^{(-,1)}_{\mathfrak{p}'}({\cal C})}

\newcommand{\bsz}{{\fB}^{(0,1)}_{\fp}({\cal C})}
\newcommand{\bszprime}{{\fB}^{(0,1)}_{\fp'}({\cal C})}

\newcommand{\bsk}{{\fB}_{\kappa,\mathfrak{p}}(T^*{\cal C})}
\newcommand{\bsbk}{{\fB}_{\beta,\kappa,\mathfrak{p}}(T^*{\cal C})}
\newcommand{\bskprime}{{\fB}_{\kappa',\mathfrak{p}'}(T^*{\cal C})}
\newcommand{\bsbkprime}{{\fB}_{\beta,\kappa',\mathfrak{p}'}(T^*{\cal C})}

\newcommand{\bscheck}{ {\fB}_{\fp} (\check{\cal C}) }
\newcommand{\bsinf}{{\fB}_{r,\infty,\rho,\sigma}({\cal C})}
\newcommand{\bskinf}{{\fB}_{r,\infty,\rho,\sigma}({\cal C})}
\newcommand{\bskcheck}{ {\fB}_{\kappa,\fp} (T^*\check{\cal C}) }
\newcommand{\bshat}{ {\fB}_{\fp} (\hat{\cal C}) }
\newcommand{\bshatprime}{ {\fB}_{\fp'} (\hat{\cal C}) }
\newcommand{\bskhat}{ {\fB}_{\kappa,\fp} (T^*\hat{\cal C}) }

\newcommand{\bspi}{\fB_{\fp}(\Pi\times\T^n)}
\newcommand{\bspiprime}{\fB_{\fp'}(\Pi\times\T^n)}

\newcommand{\bsp}{{\fB}_\sigma(\T^n)\times[\bs]^{n+1}}
\newcommand{\bspprime}{{\fB}_{\sigma'}(\T^n)\times[\bsprime]^{n+1}}

\newtheorem{definition}{Definition}
\newtheorem{assumption}{Assumption}
\newtheorem{hyp}{Hypothesis}
\newtheorem{theorem}{Theorem}
\newtheorem{lemma}{Lemma}[section]
\newtheorem{corollary}{Corollary}[theorem]
\newtheorem{proposition}{Proposition}[subsection]

\newcounter{pups}[section]
\setcounter{pups}{0}

\newenvironment{remark}{

\medskip\noindent\addtocounter{pups}{1}{\tt Remark \arabic{section}.\arabic{pups}:}}{

\medskip\noindent}

\title{Hamilton-Jacobi
method for a simple resonance}
\author{Mischa Rudnev\thanks{Partially supported by the NSF grant DMS 0072153, the Nuffield Foundation grant
NAL/00485/A, and the EPSRC grant GR/S13682/01.\newline
Contact address: Department of Mathematics University of Bristol
University Walk, Bristol BS8 1TW, UK;\newline  e-mail: {\tt m.rudnev@bris.ac.uk}}}

\maketitle

\begin{abstract}
It is well known that a generic small perturbation of a Liouville-integrable
Hamiltonian system causes breakup of resonant and near-resonant
invariant tori. A general approach to the simple resonance
case in the convex real-analytic setting is developed, based on a new technique for solving
the Hamilton-Jacobi equation. It is shown that a generic perturbation creates in the core of a
resonance a partially hyperbolic lower-dimensional invariant torus, whose
Lagrangian stable and unstable manifolds, described as global solutions of the Hamilton-Jacobi
equation, split away from this torus at exponentially small angles. Optimal
upper bounds with best constants are obtained for exponentially small
splitting in the general case.
\end{abstract}

\tableofcontents

\section{Introduction}
The notion of Arnold diffusion refers to a generic instability of
Hamiltonian systems with three and higher degrees of freedom
\cite{{A1},{A3}}. A notable exception are Liouville-integrable
systems allowing the construction of global action-angle variables
\cite{A2}. Small perturbations of such systems provide a natural
set-up to study the instability. Recently Mather \cite{Ma3} using
methods of analysis in the large (which to a great extent had been
created by himself) \cite{{Ma1},{Ma2}} announced the proof of the
existence of Arnold diffusion in the three degrees of freedom
(convex, real-analytic) case.

Consider a Hamiltonian system of $n+1,\,n\geq2$ degrees of freedom
in the cotangent bundle $T^*\T^{n+1}\cong \R^{n+1}\times\T^{n+1}$
of a torus $\T^{n+1}\equiv(\R/2\pi\Z)^{n+1}$. Take an open convex
domain $\Omega\subseteq\R^{n+1}$. The phase space ${\cal
M}=\Omega\times\T^{n+1}$ has a natural exact symplectic structure
$\mathfrak w$. Consider a Hamiltonian function
\begin{equation}
H:\,{\cal M}\rightarrow\R, \mbox{ such that }|H-H_0|_{\cal
M}=\eps\ll1,\mbox{ for some }H_0:\,\Omega\rightarrow\R.
\label{one}
\end{equation}
Suppose $H$ is real-analytic, i.e. it can be extended holomorphically into a
neighborhood of ${\cal M}$ in $\C^{2n+2}$, let $|\cdot|_{\cal M}$ above be the
supremum-norm. Also suppose that $H_0$ is strictly convex. Then one can simply
take $\Omega=\{\p\in\R^{n+1}:\,H_0(\p)<E_0\}$ for some $E_0>1$.

If $(\p,\q)$ are (global) canonical coordinates on ${\cal M}$, or the
action-angle variables \cite{A2} with $\mathfrak w=d\q\wedge d\p$, the Hamiltonian
(\ref{one}) has an expression
\renewcommand{\theequation}{\arabic{equation}$'$}
\setcounter{equation}{0}
\begin{equation}
H=H(\p, \q,\eps)=H_0(\p)+\eps H_1(\p, \q,\eps),\,\,\,(\p, \q,\eps)\in
\Omega\times\T^{n+1}\times\R_+.
\label{oneprime}
\end{equation}
The Hessian matrix $D^2 H_0(\p)$ is positive definite for every $\p\in \Omega$;
$0\leq\eps\ll1$ is a small parameter. The perturbation $H_1$ is $2\pi$-periodic
in each angle $ q_j,\,j=0\ld n$. The system (\ref{one}) is
autonomous\footnote{The case when $H_1$ depends on time periodically can be
treated in the usual way \cite{A2} whereupon the convexity assumption about
$H_0$ should be substituted by quasi-convexity \cite{Lo} and the non-degeneracy
assumption in Theorem \ref{mt} - by isoenergetic non-degeneracy \cite{DS}.
Convexity is far the easiest non-degeneracy assumption to deal with; for more
subtle non-degeneracy settings in the KAM theory see e.g. \cite{Se}. \label{ft1}}. If $\eps=0$, it is
Liouville integrable. Its phase space is foliated by invariant tori, whereupon
$\p(t)=$ const., and $\q(t)={\rm const.}+\w (\p)t$, where $\w(\p)=D H_0(\p)$ is
a frequency. Each torus is a Lagrangian manifold.

The central question of {\em local analysis} of system (\ref{one}) is what geometric objects
replace the invariant tori when $\eps\neq0$.
The KAM theorem \cite{{Ko},{Z2},{Ne},{Po1}} asserts that as
$\eps\rightarrow0_+$, an asymptotically  full measure set of these tori is stable.

{\em Resonant} unperturbed tori are foliated by tori of lower
dimension. The property of a torus being resonant or non-resonant
is clearly intrinsic, as well as the notion of the multiplicity
$m$ of a resonance, i.e. the difference in the dimensions of the
original resonant torus and the minimum foliation torus. In the
above coordinate representation the resonances correspond to the
values of the action $\p$, when the components of the frequency
vector $\w(\p)\in\R^{n+1}$ are linearly dependent over the
integers $\Z$, $m$ being the dimension of the kernel of a linear
map $\w[\ka]=\spr{\ka}{\w},$ for $\ka\in\Z^{n+1}$. Resonant tori
as well as the non-resonant ones sufficiently close to the former,
typically get destroyed for any $\eps\neq0$. The set of destroyed
tori is residual on the unperturbed energy surface $H_0^{-1}(E)$
for a regular value of $E$. It is known \cite{Tr} that given a
specific resonance, the  majority (in the sense of the Lebesgue
measure in $\R^{n+1-m}$) of the corresponding resonant tori result
in particular in the appearance of partially hyperbolic, or
whiskered tori of dimension less by $m=1,\ldots,n$. Characteristic
exponents of these tori are typically $O(\seps)$. Singular
perturbation theory for manifolds asymptotic to these tori has a
number of subtleties \cite{{Ch},{CG},{Ga},{GGM},{RW4},{DG}} which
would not be there, were the above characteristic exponents
$O(1),$ see also \cite{{CG},{BT},{DLS}}.

Non-resonant tori, sufficiently close to resonances experience a complicated
topological perestroika. If $n=1$, the result is a cantorus \cite{Ma1}
supporting an invariant action-minimizing measure \cite{Ma2}.
Higher-dimensional relatives of cantori are not so well understood, unless the local
analysis can be in a sense reduced to the $n=1$ case \cite{Xi}.

\medskip
\noindent The purpose of this paper is to develop from scratch the local theory
for a simple resonance, $m=1$. A resonance is identified by an
integer lattice point $\ka_0\in\Z^{n+1}\setminus\{0\}$. An unperturbed torus,
marked by $\p=\p_0\in\Omega$ is resonant with respect to $\ka_0$ iff the
corresponding frequency $\w_0=\w(\p_0)$ lies on the ``resonant hyperplane''
$\{ \w\in\R^{n+1}:\,\spr{\ka_0}{\w}=0\}$. As $H_0$ is smooth and strictly convex, the
``frequency map'' $\p\rightarrow\w(\p)$ is a global diffeomorphism. The values
of the action $\p$ satisfying the above resonance condition lie on a smooth
hypersurface in $\R^{n+1},$ which  intersects each regular level set of $H_0$
transversely (for otherwise $\spr{\ka_0}{D^2H_0(\p_0)\ka_0}=0$) and is a graph
over the hyperplane $\spr{\ka_0}{\p}=0$. Thus metric and topological
properties of sets on the resonant hypersurface can be described in
terms of their images in the resonant hyperplane, via the frequency map.

Given $\ka_0$, one chooses a value $\p_0$ on the intersection of a
regular level set of $H_0$ with the resonance hypersurface, such
that the corresponding frequency $\w_0$ is non-resonant over
$\Z^{n+1}$ modulo one-dimensional sub-lattice generated by
$\ka_0$. If one denotes the corresponding unperturbed simple
resonance $(n+1)$-torus as $\boldsymbol{\cal T}_0$, the latter is
foliated by a one-parameter family of $n$-tori, which can be
parameterized by some $x\in\T$:
\renewcommand{\theequation}{\arabic{equation}}
\setcounter{equation}{1}
\begin{equation}
\boldsymbol{\cal T}_0=\bigcup_{x\in\T}{\cal T}_x.
\label{fol}
\end{equation}
It is assumed that $\w_0$
is ``far enough'' from higher multiplicity resonances. To express the latter
property, Kolmogorov's Diophantine condition \cite{Ko} over the quotient
lattice is used. The set of all such frequencies $\w_0$ has a positive Lebesgue
measure on the resonance hyperplane \cite{DS}.

Study of simple resonances and their role in global dynamics for
the general system (\ref{one}) had begun at least as early as
Poincar\'{e} \cite{P}. Arnold \cite{A1} used a simple resonance
model to suggest a local mechanism for universal instability, or
diffusion, based on the existence of intersections of Lagrangian
manifolds, asymptotic to whiskered tori, alias the splitting of
separatrices phenomenon. Splitting in a more general context was
studied by Chirikov \cite{Ch} emphasizing its role in the general
diffusion scenario and conjecturing a number of generic asymptotic
exponentially small bounds apropos of the splitting and the
diffusion speed. For the latter, the theorem of Nekhoroshev
\cite{{Nek},{Lo},{Po}} gives the upper bound
$\sim\exp\left(\eps^{-{{\rm const.}\over2(n+1)}}\right)$.

More recently models for simple resonances and splitting of
separatrices have been investigated in a great number of works,
see \cite{{CG},{Ga},{El},{GGM},{RW3},{DG},{LMS}} among others. For
a more extensive bibliography list see the treatise \cite{LMS} by
Lochak et al, to which one can add some 30 more titles which have
become available since the year 2000. The latter work \cite{LMS}
among other things develops a normal form theory for local
near-resonance dynamics, see also \cite{{Lo},{Po}}). However, the
underlying multiple step averaging procedure is rather general and
does not allow to study the splitting in all the detail. As an
alternative Lochak et al advocate the Hamilton-Jacobi method,
which they illustrate for a particular Hamiltonian from the Arnold
example \cite{A1} (also published separately as \cite{Sa}) and
draft formulations of a number of theorems, which are proved
herein.

A fundamental question apropos of exponentially small splitting
(to which the Nekhoroshev-like normal form theory fails to provide
an answer) is one of the best constants for the upper estimates
involved. Such constants have been obtained for various cuts of a
specific model, coupling a pendulum-like one degree of freedom
Hamiltonian system with a bunch of rotators
\cite{{CG},{Ga},{DGJS},{GGM},{RW4},{DG}}. Ideally, the upper
bounds would be supported by lower bounds, which constitute a very
delicate issue and are available only for a few particular
examples \cite{{DGJS},{RW1},{GGM},{LMS}}. The issue is not
addressed in this paper.

An important result concerning the splitting problem in the
general simple resonance context is due to Eliasson \cite{El} (see
also \cite{DG}) who proved the estimate $2n+2$ on the minimum
number of homoclinic orbits to a whiskered torus of dimension $n$
at the resonance core, but not the exponentially small splitting
estimate. The main building blocks for the splitting theory near a
simple resonance are presented in \cite{LMS} although many are
without proofs, apparently due to a variety of technical
difficulties. This paper attempts to do it, as it turns out that
most of these difficulties can be bypassed owing to a technique,
rather different from those used in the above listed references
(for the exception of \cite{Sa}) and appears to be more
``natural'' for the problem involved. The technique certainly
applies to the above mentioned model, for which it gives the
(known) best constant ${\pi\over2}$ and also shows that the latter
is the largest value that the best constant in question can assume
in principle.

In essence, our technique is the  Hamilton-Jacobi approach prompted by
Poincar\'{e} \cite{P} cast as a ``hyperbolic KAM theorem". However it is
developed in an entirely different geometric context than the traditional one
founded by Graff \cite{Gr}. The present geometric scenario was founded in
\cite{RW2}; this paper shapes it into ``KAM theory on semi-infinite bi-cylinders
over tori''.

\medskip
\noindent
The paper is organized as follows.  Section 2 starts out with the
preliminaries in order to describe the standard normal form near a chosen simple
resonance (\ref{nfh}) and Lemma \ref{nfl}. A non-degeneracy Assumption \ref{gp} is
made concerning the ``hyperbolic part'' of the truncated (integrable) normal
form, whereupon the splitting problem is set up somewhat heuristically for the
localization (\ref{wham}) of the normal form Hamiltonian near the truncated
normal form separatrix. The set-up emphasizes what is called a ``sputnik''
property (\ref{trans}) thereof, combining the $2\pi$-periodicity of the foliation (\ref{fol})
and reversibility of the truncation
(\ref{tnfh}) of the normal form Hamiltonian (\ref{nfh}). At that point one of the main
results, Theorem \ref{main} of the paper is formulated. The formulation is
still somewhat heuristic, due to the necessity of developing a certain
amount of machinery.

Section 3 develops this machinery, underlying the aforesaid
version of KAM theory. It is based on a simple holomorphic map
introducing ``energy-time'' coordinates (\ref{stime}, \ref{time})
in which the base space is not compact. Section \ref{et} is almost
entirely dedicated to the relevant formalism. The hyperbolic KAM
theorem, Theorem \ref{mt} follows, providing global generating
functions for perturbed separatrices as the solutions of the
Hamilton-Jacobi equation, Corollary \ref{hjc}. The proof of
Theorem \ref{mt} incorporates a two-parameter trick, yielding an
optimal (in the sense of the parameter dependence) smallness
condition (\ref{munot}).

Section 4 presents the theory for the splitting, based on
application of Theorem \ref{mt} and a global sputnik property,
Assumption \ref{SP}. The role of the sputnik is to ensure that the
one-form giving the splitting distance be exact, which allows one
immediately to give a lower bound for a number of homoclinic
orbits \cite{El}. This is inherent in the homoclinic splitting
problems, being expressed by relations (\ref{trans}, \ref{transs},
\ref{transss}). The main theorem of the section, Theorem
\ref{estimate} claims the principal exponentially small estimate
(\ref{expsmall}) which is adapted as (\ref{sir}) to the simple
resonance normal form in Theorem \ref{main}, concluding its proof.
The estimate (\ref{sir}) contains a pair of best constants
$(\rho,\sigma_2)$, well defined for a specific Hamiltonian with a
given analyticity domain, see (\ref{rho}, \ref{sone}).

\section{Normal form near separatrix at a simple resonance}
\renewcommand{\theequation}{\arabic{section}.\arabic{equation}}\setcounter{equation}{0}
This section prepares the Hamiltonian (\ref{one}) for the set-up
of the theory developed in the sequel. It consists in choosing a
simple resonance action value and restricting the Hamiltonian to
its small neighborhood, where a suitable normal form can be
produced. Under generic assumptions, this normal form can be
viewed as a perturbation of an integrable reversible system
containing a separatrix. These steps are standard, see e.g.
\cite{{Tr},{El}}. However the further analysis is essentially
different from that one traditionally encounters in the
literature. The Hamiltonian gets localized near one branch of the
separatrix of the truncated normal form Hamiltonian. Localization
near the other branch can be seen from the former one as a
symmetry, referred to as a {\em sputnik}.

\subsection{Preliminaries}
\subsubsection*{Canonical transformations of the phase space} Throughout the paper,
a number of canonical transformations is introduced. These transformations
belong to an ``affine'' class, corresponding to the phase space
bundle structure.

For the system (\ref{one}) the phase space is a subset of $T^*\T^{n+1}$,
with canonical coordinates $(\p, \q)\in\R^{n+1}\times\T^{n+1}$. An
automorphism $\ar$ of the base space  induces a family of canonical transformations

\begin{equation}
{\Xi}\,=\,{\Xi}(\ar,S):\ \left\{
\begin{array}{llllllll}
\q&=&\ar(\q'),\\ \p&=&\ti{(d\ar)}\p'+dS,
\end{array}
\right. \label{zed}
\end{equation}
parameterized by a closed one-form $dS$ on the base space. As the latter is
$\T^{n+1}$, the one-form $dS$ is described by the generating function
$S(\q)=\spr{\boldsymbol \xi}{\q}+\hat{S}(\q)$, with some $\boldsymbol
\xi\in\R^{n+1}\equiv[dS]$ specifying the $H^1(\T^{n+1},\R)$ cohomology class of $dS$,
and some function $\hat{S}(\q),$ which is a zero-form on $\T^{n+1}$, i.e. is
$2\pi$-periodic in each component $ q_j,\,j=0\ld n$ of $\q$.

The notation $\spr{}{}$ stands for the canonical coupling between
$\R^{n+1}$ (and later $\R^n$) and its dual space and is identified
with the Euclidean scalar product; $\ti{(d\ar)}$ denotes the
transpose inverse of the Jacobi matrix $d\ar$.

The case when the map $\ar$ is linear, so $\ar$ and $d\ar$ can be
identified with a matrix from $SL(n+1,\Z)$ is referred to as a
symplectic rotation. The one-form $dS$ effects a shift of the
origin within each fiber. This shift is fiber-independent if
$S(\q)=\spr{\p_0}{\q},$ given $\p_0\in\R^{n+1}$, which then
becomes the origin in each fiber.

\subsubsection*{Simple resonance Diophantine condition}

A resonance is identified by a minimal integer lattice point
$\ka_0\in\Z^{n+1}\setminus\{0\}$. I.e. there is
no $\ka\in\Z^{n+1}$, such that $\ka_0= j\ka,\,j=2,3,\ldots.$
Then in any lattice basis  $\boldsymbol e_0,\ldots,\boldsymbol e_n$, the components of $\ka_0$
are relatively prime.
The choice of the basis $\{\boldsymbol e_j\}_{j=0,\ldots,n}$ determines a
coordinate chart $(\p,\q)$: take $q_j\in\R \boldsymbol e_j/2\pi\Z \boldsymbol
e_j,\,j=0,\ldots,n$, let $p_j$ be a momentum canonically conjugate to $q_j$.

Consider a one-dimensional lattice $\Z \ka_0$ and a direct sum decomposition
$\Z^{n+1}=\Z \ka_0\oplus\ \Z^{n+1}/\Z \ka_0 $. Choose a lattice basis in the
quotient lattice $\Z^{n+1}/\Z \ka_0$ and let $\ka_1\ld \ka_n\in \Z^{n+1}$
represent it in $\Z^{n+1}$, so
$\{\ka_j\}_{j=0,\ldots,n}$ is a lattice basis in $\Z^{n+1}$. For a moment, let
us call it a {\em direct sum decomposition basis,} generated by $\ka_0$. One
can expand each  $\ka_j$ over the ``old'' basis $\{\boldsymbol
e_i\}_{i=0,\ldots,n}$. As in the first equation in (\ref{zed}) let us write it
as $\ka=\ar(\boldsymbol e)$, where a linear operator $\ar$ is identified with a
matrix from $SL(n+1,\Z)$. The first row of this matrix simply gives the
coordinates of $\ka_0$ in the basis $\{\boldsymbol e_j\}_{j=0,\ldots,n}$; the
rest of the rows depend on a particular choice of the basis in the quotient
lattice $\Z^{n+1}/\Z \ka_0$. Clearly $\ka_0$ defines a direct sum decomposition
basis modulo  $SL(n,\Z)$.

\begin{definition} A vector $\w\in\R^{n+1}$ is Diophantine modulo $\ka_0$
with an exponent $\tau\geq n-1$ and a constant $\gamma>0$ (one writes $\w\in \mathfrak W^{n+1}_{\tau,\gamma}(\ka_0)$)
if  $\spr{\ka_0}{\w}=0$ and there exists
a direct sum
decomposition basis $\{\ka_j\}_{j=0,\ldots,n}$
generated by $\ka_0$, such that for all $\ka\in\Z^{n+1}$,
represented as $\ka=(k_0\ld k_n)\in\Z^{n+1}$ in this basis, one has
\begin{equation}
|\spr{\ka}{\w}|\geq\gamma|\ka|_{\hat 0}^{-\tau},\,\,{\rm where}\ |\ka|_{\hat
0}=\sum_{j=1}^n|k_j|. \label{diomod}
\end{equation}\label{dmd}
\end{definition}
Note that given $\ka\in\Z^{n+1}$ and $\w\in\R^{n+1}$, the ``small divisor"
$\spr{\ka}{\w}$ does not depend on the choice of the lattice basis $\{\boldsymbol
e_j\}_{j=0,\ldots,n}$. However, the quantity $|\ka|_{\hat 0}$ in the right hand
side of (\ref{diomod}) is not $SL(n,\Z)$-invariant and can attain any positive
integer value for a given $\ka$.

Definition \ref{dmd} is rather unwieldy, and in order to describe metric
properties of Diophantine vectors one is forced to fix the lattice basis
$\{\boldsymbol e_j\}_{j=0,\ldots,n}$. As the lattice element $\ka_0$ is
considered fixed throughout the paper, one should naturally render
$\{\boldsymbol e_j\}_{j=0,\ldots,n}$ and $\{\ka_j\}_{j=0,\ldots,n}$ the same
basis. Then $\ka_0=(1,0,\ldots,0)$ and for any $\w$ such that
$\spr{\ka_0}{\w}=0$, clearly $\w=(0,\omega)$ for some $\omega\in\R^n$. Given
the pair $(\tau,\gamma),$ let us further use the notation ${\mathfrak
W}^n_{\tau,\gamma}$ for a set of $\omega\in\R^n$, satisfying a stronger
definition than Definition \ref{dmd}:

\begin{equation}
\di=\left\{\omega\in\R^n:\;\forall\,k\in\Z^{n}\setminus\{0\},\;\;
|\spr{k}{\omega}|\geq\gamma|k|^{-\tau}, \;\;
|k|=\sum_{i=1}^n|k_i|\right\}.\label{dio}
\end{equation}
The above is the standard Kolmogorov Diophantine condition \cite{Ko}. Some
well known metric properties of the set $\di$ (with a sufficiently small
$\gamma$, see e.g. \cite{{Py},{DS}}, also \cite{Do} for complementary results)
are that it is non-empty if $\tau\geq n-1$, and for $\tau>n-1$ has a full
Lebesgue measure as $\gamma\rightarrow0$.

\subsection{Localization near a simple resonance}
For the Hamiltonian (\ref{one}) given a resonance $\ka_0,$ let
$\p_0\in\Omega$ lie on a regular level set $H_0^{-1}(E)$ and
$\w_0=DH_0(\p_0)$ satisfy Definition \ref{dmd}. Then one can
choose the lattice basis $\{\boldsymbol e_j\}_{j=0,\ldots,n}$ with
$\boldsymbol e_0=\ka_0$, such that $\w_0=(0,\omega_0)$, with
$\omega_0\in\di$, with some $\tau$, which is fixed throughout the
paper and $\gamma=\gamma_0$. Moreover, the origin for the action
variables can be set at $\p_0$, i.e. $\p_0=(0,\ldots,0)$. As $H$
is defined modulo a constant, let $E=H_0(0)=0$.

This fixes the choice of the action-angle variables $(\p,\q)$ and denoting
$\p=(y,I)\in \R\times\R^n,\,\q=(x,\varphi)\in\T\times\T^n$, one can write down
the following representation for the Hamiltonian (\ref{one})\footnote{If
$\ka_0\neq(1,0,\ldots,0)$ in the given lattice basis and the Hamiltonian $H$ is
given in the form (\ref{oneprime}), the representation (\ref{localmod}) can
certainly be achieved by means of a canonical transformation in the form
(\ref{zed}) combining a symplectic rotation and a shift of the action origin to
$\p_0$ \cite{{DS},{Tr}}. With a specific $H$ and $\ka_0$ in mind, the
adaptation of the parameters in the main estimate (\ref{sir}) in Theorem
\ref{main} is straightforward, similar to how it will have embraced the
parameter $\theta$ in the sequel.}:

\begin{equation}
H(\p,\q,\eps)\,=\,\spr{\omega_0}{I}+{1\over2}\spr{Q_0\p}{\p}+O_3(\p)
\;+\;\eps H_1(\p,\q,\eps),
\label{localmod}
\end{equation}
where $Q_0$ is a constant
positive definite matrix.

\renewcommand{\theenumi}{N.\arabic{enumi}}

\medskip
\noindent {\bf Notation:}
\begin{enumerate}
\item The set of non-negative or positive, integer or real numbers is denoted
as $\Z_+$ or $\Z_{++}$, $\R_+$ or $\R_{++}$, respectively.

\item Bold lowercase symbols usually denote $(n+1)$-vector quantities.
Uppercase symbols often but not always denote $n$-vector quantities. E.g.
above $\p=(y,I),\,\q=(x,\varphi)$; further $\g(\q)=(g(\q),G(\q))$.

\item The symbol notation $u(z,\cdot)=O_{\alpha}(z;\cdot)$, with
$\alpha\in\R$, implies that  $\displaystyle\lim_{z\rightarrow0}
{|u(z,\cdot)|\over \|z\|^\alpha}$, where $\|\cdot\|$ is the
Euclidean norm, exists and is uniformly bounded from above for the
whole range of the variables $(\cdot)$ (which may be omitted in
the notation, as well as $\alpha=1$) by some constant $C$ which
may depend on $(n,\tau)$ and perhaps other quantities fundamental
for the problem, to be specified. $C$ will be ``as large as
necessary'' and may increase without notice. To suppress $C$ (or
$C^{-1}$) in estimates, the $\lesssim$ sign is often used instead
of $\leq$.

\item For real $\kappa,\sigma>0$ and $j\in\Z_{++}$ ($j=1$ usually being omitted) let
\[
\begin{array}{lll}
\B^j_\kappa &\=&\{z\in\C^j:\,\|z\|\leq\kappa\},\\ \T^j_\sigma & \=
&\{z\in\C^j:\,\Re{z}\in\T^j,\,|\Im{z}|\leq\sigma\}
\end{array}
\]
be complex extensions of a disk and a torus.

It will always be assumed by default that $\kappa>1$, and it will not enter the
estimates. All the analyticity and non-degeneracy parameters are by default
positive, as well as $\varepsilon$. In addition, if $\delta$ for instance
denotes the analyticity loss in the variable $\varphi\in\T^n_\sigma$, it will be assumed by
default that $0<\delta<\sigma$.

\item Scalar functions $2\pi$-periodic in each variable,
holomorphic and uniformly bounded inside $\T^j_\sigma$, whose
restrictions on $\T^j$  are real-analytic form Banach spaces
$\fB_\sigma(\T^j)$,  with topology induced by the supremum norm
$|\cdot|_{\sigma}$. The space of all Taylor series with
coefficients in $\fB_\sigma(\T^j)$, uniformly convergent inside
$\B^j_\kappa$, with the supremum norm $|\cdot|_{\kappa,\sigma}$ is
denoted as $\fB_{\kappa,\sigma}(T^*\T^j).$ Referring to
real-analytic functions on complex domain in the sequel means
referring to their holomorphic extensions.

The same notation stands for the supremum norms of vector functions.
At places the subscripts in the norm notation can be omitted.
If a function, whose norm is evaluated
depends on additional parameters, omitting these dependencies in the
estimates implies their uniformity.

\end{enumerate}
Thus, in the convex real-analytic set-up, there exists a set of parameters
$\{\gamma_0,\kappa_0,\sigma_0,R_0,M_0,\eps_0\}$, with $\kappa_0>1$ and
$\eps_0\ll1$, such that the Hamiltonian (\ref{localmod}) satisfies the
following.

\renewcommand{\theenumi}{\arabic{enumi}}

\medskip \noindent {\bf Model statement:} {\em \begin{enumerate}
\item For all $\eps\in[0,\eps_0]$,  $H(\p,\q,\eps)\in{\fB}_{\kappa_0,\sigma_0}(T^*
\T^{n+1})$ and $|H_0|_{\kappa_0}\leq M_0,\,|H_1|_{\kappa_0,\sigma_0}\leq1.$
\item The frequency $\omega_0\in{\mathfrak W}^n_{\tau,\gamma_0}$.
\item The constant  matrix $Q_0$ is positive definite, $\|Q_0^{-1}\|\leq
R_0^{-1}$.
\end{enumerate} }
\noindent
For a specific Hamiltonian (\ref{oneprime}) the analyticity
considerations may somehow single out the choice of the original coordinates
$(\p,\q)$, see e.g. \cite{DGJS}. As a result, the parameters above as well as
the bounding constants may also depend on $\ka_0$.

\subsection{Normal form near a resonance}
It is well known that one can come up with a  normal form near a resonance \cite{DS}.
Such a normal form for the Hamiltonian (\ref{localmod}) was
used as a motivation for the results of \cite{{El},{RW2},{RW3}} among others.
Note that the normal form
transformation belongs to the class (\ref{zed}) and cannot be iterated in this
form.

As $\kappa_0>1$, one can accept $M_0$ and $1$ respectively as
bounds for several orders of derivatives of $H_0$ and $H_1$ in
$\p$. Assuming $R_0<1$, let us get rid of this parameter as far as
the quadratic part of $H_0$ is concerned. Rewrite (\ref{localmod})
scaling the time and actions by factor $\sqrt{\eps/ R_0}$, i.e.
$\p\rightarrow \sqrt{\eps/R_0}\, \p,\,H\rightarrow
\sqrt{R_0/\eps}\,H$, and then divide the Hamiltonian by
$\sqrt{\eps R_0}$ (tantamount to yet another time scaling). Then
(\ref{localmod}) changes to

\[
H(\p,\q,\eps)=\spr{\omega_1}{I}+{1\over2 }\spr{Q_1\p\,}{\p}
 +\eps ^{-1} O_3(\sqrt{\eps/R_0}\p) +
 H_1(\sqrt{\eps/R_0}\p, \q,\eps),
\]
with the notations
\[
\omega_1= {\omega_0\over \sqrt{\eps R_0}},\;\;\;
Q_1=R_0^{-1}Q_0.
\]
The analyticity domain of the scaled Hamiltonian $H$ in the scaled action
$\p$ is now a complex ball of radius $O(\sqrt{R_0/\eps})$. The matrix $Q_1$ is
such that its spectrum is contained in $[1,C(n)M_0/R_0]$.

Decompose the quantity $H_1(0,\q,\eps)=H_1(0,x,\varphi,\eps) $ into a $\varphi$-mean $U(x,\eps)$
and an oscillatory part:
\[ \begin{array}{llllllll}
H_1(0,x,\varphi,\eps)&=&\int_{\T^n}H_1(0,x,\varphi,\eps)d\varphi&+&\left(H_1(0,x,\varphi,\eps)-
\int_{\T^n}H_1(0,x,\varphi,\eps)d\varphi \right) \\ \hfill \\
&\equiv & U
(x,\eps)&+&\{H_1(0,x,\varphi,\eps)\}.
\end{array}
\]
To get rid of the $\varphi$-oscillatory term
$\{H_1(0,x,\varphi,\eps)\}$ consider a canonical transformation
$\Xi_\nu$, which is tantamount to the shift
$\p\rightarrow\p+dS_\nu(\q,\eps)$. The 1-form $dS_\nu$ is exact
and is given by a $2\pi$-periodic in each component of
$\q=(x,\varphi)$ function $S_\nu$, satisfying a PDE
\[ \displaystyle D_{\omega_1}
S_\nu(x,\varphi,\eps)\,=\,-\{H_1(0,x,\varphi,\eps)\}, \] with the
general notation for $\omega\in\R^n$
\begin{equation}
D_\omega\= \spr{\omega}{D_\varphi}. \label{dom}
\end{equation}
Note that $x,\eps$ enter the above equation  as parameters, in particular $S_\nu$
is defined modulo a function of $x$. The solution of this equation exists in a
somewhat larger space than that for the right hand side. The following  result
is well known \cite{Ru}.

\begin{proposition}
Let $\omega\in\di$. For a function $v\in\bst$ with zero average on $\T^n$, the solution of the equation
$D_\omega u=v$ exists in the space $\bstprime$ for any
$\sigma'<\sigma$. If
$\sigma-\sigma'=\delta,\,\varsigma=\gamma\delta^{\tau}$, then
\[
|u|_{\sigma'}\,\lesssim\,
\varsigma^{-1}|v|_{\sigma},\;\;\;|du|_{\sigma'}\,\lesssim\,
(\varsigma\delta)^{-1}|v|_{\sigma}.
\]
\label{rus}
\end{proposition}
Then given $\eps$ small enough to ensure that the transformation
$\Xi_\nu$ be near identity, i.e. $|dS_\nu|\ll 1$, suppressing
$\eps$ in the notation one gets for $H_\nu\equiv H\circ\Xi_\nu$:
\begin{equation}
H_\nu(\p,\q)\;=\;\spr{\omega_1}{I}+ \,{1\over2}\spr{Q_1\p}{\p} +
U(x)\,\,+\;\;
[ f_1(\q) +
\spr{{\g}_1(\p,\q)}{\p}],
\label{nfh}
\end{equation} where
\[
|f_1|_{\sigma_1}\lesssim
M_0R_0^{-1}|dS_\nu|_{\sigma_1}^2,\;\;\;|\g_1|_{\sigma_1}\lesssim
M_0 \sup(R_0^{-1}|dS_\nu|_{\sigma_1},\sqrt{\eps R_0^{-3}}),
\]
for some $\sigma_1<\sigma_0$, with  $\g_1$ having absorbed the
momentum-super-quadratic term. In view of Proposition \ref{rus}
the estimate for $|dS_\nu|_{\sigma_1}$ is proportional to
$\gamma_1\delta_0^{\tau+1}$, where $\delta_0=\sigma_0-\sigma_1$ is
the analyticity loss and the parameter
$\gamma_1=\gamma_0(\sqrt{\eps R_0})^{-1}$ characterizes the
``fast'' frequency $\omega_1\in{\mathfrak W}^n_{\tau,\gamma_1}$ in
the sense of (\ref{dio}).

A characteristic feature of the simple resonance normal form $H_\nu$ is that
its truncation
\begin{equation}
H_{\nu,{\rm t}}(\p,\q)\;=\;\spr{\omega_1}{I}+ \,{1\over2}\spr{Q_1\p}{\p} +
U(x)
\label{tnfh}
\end{equation}
is integrable. The pair $(f_1,\g_1)$ can thus be treated as a
perturbation, which requires it to be sufficiently small in comparison with in particular $U(x)=O(1)$, i.e.
\begin{equation}
\eps \lesssim M_0^{-2}R_0[\inf(\varsigma_0\delta_0), R_0]^{2},
\label{sc}
\end{equation}
with the notation $\varsigma_0=\gamma_0\delta_0^\tau$. The results
so far are summarized as follows.

\begin{lemma}[Normal form lemma] Let $\p_0\in \Omega$ lie on a regular level set of $H_0$ and
for some $\ka_0\in\Z^{n+1}$, let $\w_0=DH_0(\p_0)$ be Diophantine
modulo $(\ka_0).$ Suppose, the localization (\ref{localmod}) of
Hamiltonian (\ref{one}) near $\p=\p_0$ satisfies the Model
statement with the frequency $\omega_0\in \fW_{\tau,\gamma_0}^n$
and the set of parameters
$\{\gamma_0,\kappa_0,\sigma_0,M_0,\eps_0\}$. For
$\sigma_1<\sigma_0,$ let
$\delta_0=\sigma_0-\sigma_1,\,\varsigma_0=\gamma_0\delta_0^\tau,$
and suppose
\[
\eps_0\lesssim (\varsigma_0\delta_0)^2.
\]
For any $\eps\in(0,\eps_0)$, the Hamiltonian (\ref{one}) can be cast into the normal form (\ref{nfh})
where:
\begin{enumerate} \item
$H_\nu\in{\fB}_{\kappa_1,\sigma_1}(T^*\T^{n+1}),$ with
$\kappa_1=O(R_0/\seps)$; \item The constant matrix $ Q_1$ is
positive definite, with $\|Q_1^{-1}\|\lesssim 1,\;\|Q_1\|\lesssim
M_0/R_0$; \item One has
\begin{equation}
|f_1|_{\sigma_1}\lesssim \eps M_0
(\varsigma_0\delta_0)^{-2},\;\;\;|\g_1|_{\sigma_1,\kappa_1}\lesssim\seps
(M_0/\sqrt{R_0})[\inf(\varsigma_0\delta_0), R_0]^{-1}.
\label{muandnu}
\end{equation}
\end{enumerate} \label{nfl}
\end{lemma}

\subsection{Localization near separatrix}
Let us take a closer look at the truncated normal form (\ref{tnfh}). The second
clause of Lemma \ref{nfl} implies that without loss of generality one can
assume
\[
\spr{ Q_1\p}{\p}=y^2 + 2y\spr{\theta}{I} +  \spr{\Theta
I}{I},
\]
with a constant vector $\theta\in\R^n$,  and a constant positive
definite matrix $\Theta\in\R^{n^2}$, whose smallest eigenvalue is
at least one (a greater than one coefficient multiplying $y^2$
being favorable). In order to proceed one needs the following
assumption.

\begin{assumption}[Perturbation of general position] The function $U=U(x,\eps)$
possesses a unique uniformly non-degenerate absolute maximum on $\T$ for all
$\eps\in[0,\eps_0]$, with a characteristic exponent $\lambda\in\R_{++}$.
\label{gp}
\end{assumption}
Without loss of generality let the maximizer be $x=0$ for each $\eps$ (which as
far as the above assumption is concerned is non-essential and will be further
omitted in the notation) with $U(0)=0$. Assumption \ref{gp} then is tantamount
to the claim
\begin{equation}
\begin{array}{c}
U(0)=0,\; U_x(0)=0, \;U_{xx}(0)=-\lambda^2, \;\lambda^{-1}=O(1); \\ \hfill \\
\forall\,x_c\in\T\setminus\{0\}:\;U_x(x_c)=0, \;U(x_c)<0.
\end{array}
\label{assonu}
\end{equation}
For the truncated normal form Hamiltonian (\ref{tnfh}) the action
$I$ is an integral of motion. For $I=0$ one can single out a
one-dimensional natural integrable system, whose Hamiltonian is
${\displaystyle y^2/2+U(x)}$. This is a reversible Hamiltonian
system in $T^*\T \cong \R\times\T$. Near a zero energy level, its
phase portrait is reminiscent of the classical pendulum, Fig. 1.
There is a saddle $(x,y)=(0,0)$ connected to itself by a pair of
simple non-contractable curves $y=\pm\sqrt{-2\,U(x)}$, forming a
single $\infty$-shaped curve, further referred to as the
separatrix.

Let $U(x)=\lambda^2 U_1(x)$. Conditions (\ref{assonu}) allow one
to define a $4\pi$-periodic {\em separatrix function} $\psi(x)$,
determined in general as well as the constant $\lambda$ by the
pair $(\ka_0,H_1)$ and possibly depending on $\eps$:
\begin{equation}
\psi(x)\, = \, \left\{
\begin{array}{rl}-\sqrt{-2\,U_1(x)}, & x\in[-2\pi,0), \\
\sqrt{-2\, U_1(x)},& x\in[0,2\pi)
\end{array}\right.
\label{tchi}
\end{equation}
(having chosen the branch of the square root where $\sqrt{1}=1$). Thus
$\psi(0)=0,\,\psi_x(0)=1,\,\psi_{xx}(0)=O(1)$ and
\begin{equation}
\psi\circ l_{2\pi}\,=\,-\psi,\;\;\;l_{2\pi}:\,x\rightarrow x+2\pi.\label{ant}
\end{equation}
By (\ref{assonu}) the function $\psi$ has no other zeroes on the
real axis, but even multiples of $\pi$. For instance in the
classical pendulum case $U(x)=\cos{x}-1$, $\psi(x)\,=\,2\sin x/2$.
In the general case one can write $U_1(x)=(\cos{x}-1)V(x)$, with
some real-analytic $2\pi$-periodic function $V$, which has no
zeroes on the real axis and $V(0)=1$. Therefore
$\psi(x)=2\sin(x/2)\psi_1(x)$, where the function $\psi_1$ is
real-analytic, $2\pi$-periodic, has no zeroes on the real axis,
and $\psi_1(0)=1$. Thus the function $\psi$ is real-analytic and
has no zeroes in some neighborhood of the real axis, except even
multiples of $\pi$. In particular, this property will be valid for
$|\Im x|\leq\sigma_2\leq\sigma_1$, for some $\sigma_2$.

In the full phase space $T^*\T^{n+1}$ the separatrix is represented by a
Lagrangian manifold
\begin{equation}
W_{\tt t}=\{(\p,\q)\in\R^{n+1}\times \T^{n+1}:\; \p =dS_{\tt
t}(\q)=(\lambda\psi(x),0)\}, \label{sep0}
\end{equation}
where the function $\psi$ is viewed as a double-valued function on
$\T$, corresponding to an exact double-valued one-form $dS_{\tt
t}$ on $\T^{n+1}$, given by a $\varphi$-independent generating
function $S_{\tt t}(x,\varphi)=\lambda\int \psi(x)dx$, modulo a
constant. The separatrix forms a coinciding unstable-stable
manifold to an invariant torus ${\cal T}_{\tt t}$ at $x=0$, see
Fig. 1.

One can localize (\ref{nfh}) near the manifold $W_{\tt t}$. Let us
make a formal change\footnote{Note that one should not worry here
about the analyticity domains, as $\kappa_1$ in Lemma \ref{nfl} is
large enough.} $y\rightarrow y+\lambda\psi(x)$ and denote
${L}_\psi$ the corresponding canonical transformation, acting as
the identity on the pair $(I,\varphi)$. The transformation
$L^{-1}_\psi$ acts on the base space variable $x\in\T=\R/2\pi\Z$
as a period doubling map $\T\rightarrow\T'$, where
\[
\T'=\R/4\pi\Z.
\]
The transformation $L_\psi$ changes the phase space to
$T^*(\T'\times\T^n)$ and incurs a topological change on the
separatrix $W_{\tt t}$, doubling the point $(x,y)=(0,0)$ on its
projection on the $(x,y)$-plane, see the following Fig. 1. The
manifold $W_{\tt t}$ now corresponds to the zero section of the
bundle $T^*(\T'\times\T^n)$, with the identical zero generating
function.

Let $H_\psi\equiv H_\nu\circ H_\psi$, now a $4\pi$-periodic function of $x$:
\begin{equation}
H_\psi(\p,\q)\;=\;\lambda\psi(x)y + \spr{\omega_1 +\lambda\psi(x)\theta}{I} +
{1\over2}\spr{Q_1\p}{\p} + [ f_\psi(\q) +
\spr{{\g}_\psi(\p,\q)}{\p}].
\label{wham}
\end{equation}
In particular $f_\psi= f_1+\psi \tilde f_\psi$, where the function
$\tilde f_\psi$ is determined by ${\g}_1$ in (\ref{nfh}). At the
first glance after the transformation $L_\psi$ the bound for both
$f_\psi$ and $\g_\psi$ will be that for $\g_1$ in (\ref{muandnu}).
However this is not quite the case, as one may recall that the
generating function $S_\nu$ of the canonical transformation
$\Xi_\nu$ in Lemma \ref{nfl} is defined modulo a function of $x$.
Thus one can combine the two transformations $\Xi_\nu\circ L_\psi$
into one with the generating function $S_\nu+\lambda\int
\psi(x)dx$, which enables one to improve the above estimates as
follows:
\begin{equation}
\begin{array}{lll}
|f_{\psi}|_\sigma &\lesssim& \seps M_0 \sup[\seps(\varsigma_0\delta_0)^{-2},\, R_0^{-3/2}],\\  \hfill \\
|{\g}_\psi|_{\kappa,\sigma}&\lesssim& \seps M_0
[\sqrt{R_0}\,\inf(\varsigma_0\delta_0,R_0)]^{-1}.
\end{array}
\label{pb}
\end{equation}
\begin{remark}\label{r1} Without loss of generality, $\lambda\leq1$.
As far as the power of the parameter $R_0$ is concerned (the
second entry in the $\sup$ and $\inf$ above) the estimates depend
on whether or not the original unperturbed Hamiltonian $H_0$ in
(\ref{one}) contains super-quadratic terms.\footnote{The
perturbation in (\ref{nfh}) is evaluated not only at $y=0$, but
also quite far away from it on the separatrix. The characteristic
size of the separatrix in the original action variables of the
Hamiltonian (\ref{localmod}) is $\de\sim\sqrt{\eps R_0^{-1}}$, and
in order that the super-quadratic term be considered as a
perturbation of the (quadratic in momenta) resonant normal form,
one should have $\eps>>\de^3$, thus $\eps\ll R_0^3$. It is
nevertheless irrelevant apropos of the estimates, regarding the
preservation of the invariant torus at $x=0$, in particular
because $\psi(0)=0$ (one would have to look still closer at the
structure of the acquired term $\tilde f_\psi$ to see that).}
\end{remark}Along with $L_\psi$ let us denote ${L}_{-\psi}$ a
canonical transformation, effecting the shift $y\rightarrow
y-\lambda\psi(x)$, corresponding to the ``lower'' separatrix on
the phase portrait of the truncated normal form Hamiltonian, Fig.
1. Studying the Hamiltonian $H_{-\psi}\equiv H_\nu\circ
L_{-\psi}$, further referred to as the {\em sputnik} of $H_\psi$
in essence adds nothing new, as  ${L}_{-\psi}$ is tantamount to
${L}_{\psi}$ followed by a shift $l_{2\pi}$ of the $x$-variable,
in view of $2\pi$-antiperiodicity of the function $\psi(x)$ and
$2\pi$-periodicity in $x$ of the normal form Hamiltonian $H_\nu$.
On the other, the sputnik $H_{-\psi}$ turns out to be a convenient
way to describe $H_\psi$ on the interval $x\in[2\pi,4\pi)$.
Denoting $L_{2\pi }$ a canonical transformation, corresponding to
the extension $\boldsymbol l_{2\pi}:\,(x,\varphi)\rightarrow
(x+2\pi,\varphi)$ of the shift $l_{2\pi}$ in the base space (be it
$\T^{n+1},\,\T'\times\T^n$ or further $\R\times\T^n$; also let
$l_{-2\pi}\equiv l_{2\pi}^{-1}$) one has
\[
{L}_{-\psi}\circ L_{2\pi}=L_{2\pi}\circ
{L}_{\psi}\;\mbox{ and } \;H_{-\psi}\circ L_{2\pi}=H_\psi,
\]
after applying $H_\nu$ to the first relation.
Using
the same symbols $L_{\pm\psi}$ for the transformations effecting the change
$y\rightarrow y\pm\lambda\psi(x)$ on $T^*(\T'\times \T^n)$, one also has
\[
H_\nu=H_\psi\circ{L}_{-\psi}
=H_{-\psi}\circ{L}_{\psi},
\]
and
a useful identity follows:
\begin{equation}
H_\psi=
H_\psi\circ{L}_{2\pi}\circ{L}^{2}_{\psi}= H_\psi\circ L^{-2}_\psi\circ L_{2\pi}.
\label{trans}
\end{equation}
For the flow of the Hamiltonian $H_{\psi,{\tt t}}\equiv H_{\nu,{\tt t}}\circ
L_\psi$ on $T^*(\T'\times\T^{n})$, the invariant manifold $W_{\tt t}$ contains
a pair of invariant whiskered tori ${\cal T}_{u,{\tt t}}$ and ${\cal T}_{s,{\tt
t}}$  such $x=0$ on the former torus and $x=2\pi$ on the latter one. $W_{\tt
t}$ is the unstable manifold for ${\cal T}_{u,{\tt t}}$ and  the stable
manifold for ${\cal T}_{s,{\tt t}}$. The two tori can be identified via the
transformation $L^{-1}_\psi$.

The transformation $L_\psi$ not only ``doubles'' the base space but also
the separatrix. Indeed, the two branches thereof in the truncated normal form
$H_{\nu,{\tt t}}$ are not only graphs over the base space, but also one
over the other. Thus the tori ${\cal T}_{u,{\tt t}}$ and ${\cal T}_{s,{\tt
t}}$ not only possess an unstable/stable manifold respectively, which is $W_{\tt
t}$, the zero section of $T^*(\T'\times\T^n)$, but also its sputnik
\begin{equation}
W_{\tt t}'=\{(\p,\q)\in\R^{n+1}\times(\T'\times\T^n):\; \p
=-2dS_{\tt t}(\q)=(-2\lambda\psi(x),0)\}, \label{sepstar}
\end{equation}
which is the stable manifold for ${\cal T}_{u,{\tt t}}$ and  the
unstable manifold for ${\cal T}_{s,{\tt t}}$. This fact is in
essence reflected by (\ref{trans}); this is the symmetry in the
Hamiltonian $H_\psi$ which enables one to identify\footnote{By
looking locally at the Hamiltonian vector field generated by
$H_\psi$, an observer won't be able to tell whether it is applied
at a point $(y,I,x,\varphi)$ or
$(y-2\lambda\psi(x),I,x+2\pi,\varphi)$. Thus they won't be able to
tell $W_{\tt t}$ and $W_{\tt t}'$ as well as ${\cal T}_{u,{\tt
t}}$ and ${\cal T}_{s,{\tt t}}$ apart.} the manifolds $W_{\tt t}$
and $W_{\tt t}'$. On the other hand, $W_{\tt t}'$ is clearly a
flow-invariant zero section for the Hamiltonian $H_{-\psi,\tt
t}=H_{\nu,\tt t}\circ L_{-\psi}$. Further a sputnik Hamiltonian
will be marked by a prime, e.g. $H_{-\psi}=H'_\psi$.

\subsubsection*{Splitting problem}
As the perturbation in (\ref{wham}) is not identically zero (more
precisely the zero order term thereof in the Taylor expansion in
$\p$) the manifold $W_{\tt t}$ no longer lies inside the energy
surface of $H_\psi$ (for which it is the zero section). One can
expect the following scenario, Fig. 1.

\medskip
\noindent A perturbation of general position causes the manifold $W_{\tt t}$ to
bifurcate, or split into a pair of distinct Lagrangian manifolds, denoted as
$W_u$ and $W_s$. If the perturbation is small enough, the two manifolds can be
described as graphs of closed one-forms $dS_{u}$ and $dS_{s}$. The forms
$dS_{u,s}$ are well defined over the cylinders ${\cal I}_{u,s}\times\T^n$
respectively, where ${\cal I}_u=[-2\pi+r,2\pi-r]$, for some $r<1$
(characteristic of the quantity $\psi$) and ${\cal I}_s =
l_{-2\pi}({\cal I}_u)$. Each manifold $W_{u,s}$ contains an invariant torus
${\cal T}_{u,s}$, being the unstable manifold for ${\cal T}_u$ and the stable one
for ${\cal T}_s$.

One should be able to identify the tori ${\cal T}_{u,s}$ via the transformation
$L^{-1}_\psi$. This in particular requires that both one-forms $dS_{u,s}$
belong to the same cohomology class:
\[[dS_u]=[dS_s]=\xi\in\R^n\cong H^1({\cal
I}_{u,s}\times\T^n,\R).\] This fact is easy to establish due to the fact that
the sputnik manifold $W'_{\tt t}$ will split just the same, to which there will
correspond a pair of closed one-forms $dS'_{u}$ and $dS'_{s}$. It will be easy
to see that say  $dS_{u}$ and $dS'_{u}$ have the same cohomology class, as they
in particular describe the same torus ${\cal T}_{u}$. On the other hand
(\ref{trans}) claims a congruency between the graphs of the forms $dS'_{u}$
and $dS_{s}$.

To measure the distance between the two manifolds $W_{u}$ and
$W_{s}$ and to study their intersections, one naturally uses an
exact 1-form  $d\fS=dS_u - dS_s$, well defined on the union of two
disjoint cylinders $(\hat{\cal I}_-\cup\hat{\cal I_+})\times\T^n$
where $\hat{\cal I}_-\equiv[-2\pi+r,-r]$ and $\hat{\cal
I}_+=-\hat{\cal I}_-$. Let $\beta\in\Z_2\equiv\{+,-\}$ be the sign
of $x$, then the notation $d\fS_\beta$ stands for the restriction
of $d\fS$ on $\hat{\cal I}_\beta\times\T^n$. The values of
$\beta=+,-$ respectively correspond to the splitting of the
``upper'' and ``lower'' separatrices of the truncated  normal form
Hamiltonian $H_{\nu,{\tt t}},$ see Fig. 1. The splitting of the
sputnik manifold $W'_{\tt t}$, alias the zero section for
$H'=H_{-\psi}$ is described by $-d\fS_{-\beta}$, by (\ref{trans}).

The generating function $\fS_\beta$ is real-analytic on $\hat{\cal
I}_\beta\times\T^n$, where it satisfies a linear homogeneous
Hamilton-Jacobi equation, whose coefficients can be made constant
via a change of variables. Analyzing the result of the latter
change, one has the following theorem.

\begin{theorem}
Let $\p_0\in \Omega$ lie on a regular level set of $H_0$ and for
some $\ka_0\in\Z^{n+1}$, let $\w_0=DH_0(\p_0)$ be in ${\mathfrak
W}^n_{\tau,\gamma_0}(\ka_0),$ corresponding to a simple resonance
torus $\boldsymbol{\cal T}_0$ with the foliation (\ref{fol}) in
terms of $x\in\T$. Suppose the localization (\ref{localmod}) of
Hamiltonian (\ref{one}) near $\p=\p_0$ satisfies the Model
statement with the frequency $\omega_0\in\fW_{\gamma_0,\tau}$ and
parameters $\{\kappa_0,\sigma_0,R_0,M_0,\eps_0\}$. Suppose the
perturbation $H_1$ satisfies Assumption \ref{gp} with the
characteristic exponent $\lambda$ and the separatrix function
$\psi$. For $\sigma_1<\sigma_0,$ let
$\delta_0=\sigma_0-\sigma_1,\,\varsigma_0=\gamma_0\delta_0^\tau;$
suppose
\begin{equation}
\eps_0\leq (CM_0)^{-2}R_0[\lambda^2\inf(\varsigma_0\delta_0,
\,R_0)]^2\equiv\eta^2,
\label{vsmall}
\end{equation}
for some large enough $C=C(n,\tau,\ka_0,\sigma_0,\psi)$.

For any $\eps\in(0,\eps_0)$, continuously if $H_1(\cdot,\eps)$ is continuous,
there exists a pair of analytic Lagrangian manifolds $W_{u,s}$, intersecting at
an invariant $n$-torus ${\cal T}$, on which the flow of (\ref{one}) is
conjugate to a rotation by $\omega_0$. Each manifold is locally a graph over
$\T\times\T^n$. Away from ${\cal T}$, the distance $\fS$
between $W_u$ and $W_s$ is bounded by

\begin{equation}
| {\fS} | \;\leq\;\seps \eta^{-1}
\sum_{k\in\Z^n\setminus\{0\}}
 \exp\left(-|\spr{k}{\rho{R_0\over\lambda} {\omega_0\over \seps}+\sigma_2\theta}| - |k|\sigma_1\right),
\label{sir}
\end{equation}
where ${\displaystyle \rho<{\pi\over2}},\,\sigma_2\leq\sigma_1$,
and the quantities $\rho,\sigma_2,\theta$ are well defined for
$H$. The manifolds $W_u$ and $W_s$ also intersect along at least
$2n+2$ orbits, biasymptotic to ${\cal T}$. \label{main}
\end{theorem}
\begin{remark} In view of the smallness condition (\ref{vsmall}) and the fact that $\theta\lesssim M_0$ in
(\ref{sir}), the contribution of the quantity $\sigma_2\theta$ in
the estimate (\ref{sir}) becomes important for a Diophantine
$\omega_0$ when $|k|\sim \eps^{-{1\over 2(\tau+1)}}$ and will play
an extra role if one attempts to estimate $|\fS|$ from below
\cite{RW1}. Where exactly the quantities $\rho$ and $\sigma_2$
arise is explained further, see in particular (\ref{rho},
\ref{sone}) and Fig. 2.\end{remark}

\noindent
It is possible to simplify (\ref{wham}) further
 by eliminating the constant $\theta\in \R^n$ therein, letting
\begin{equation}
{\Xi}_\theta:\;\left\{
\begin{array}{lll}
x&=&x', \\ \varphi&=&\varphi'+ \theta x',
\end{array}
\right.
\hspace{1mm}  \begin{array}{lll}
y&=&y'-\spr{\theta}{I},
 \\
I&=&I'.
\end{array}
\label{two}
\end{equation}
The base space transformation, corresponding to ${\Xi}_\theta$ will be denoted
as $\ar_\theta$. Unless $\theta\in\Q^n$, the pre-image of the base space
$\T'\times\T^n$ under the transformation $\ar_\theta$ is a bi-infinite cylinder
$\R\times\T^n$. In other words, the transformation ${\Xi}_\theta$ almost surely results in the loss of
$4\pi$-periodicity in the ``hyperbolic coordinate'' $x$.

If ${H}_\theta\equiv H_\psi\circ\Xi_\theta$  then
\begin{equation}
{H}_\theta(\p,\q)\;=\;\lambda\psi(x)y + \spr{\omega_1}{I} +
{1\over2}\spr{Q_2\p}{\p} + [f_\theta(\q) +
\spr{{\g}_\theta(\p,\q)}{\p}],
\label{wwham}
\end{equation}
where the matrix $Q_2$ arises from $Q_1$ as a result of Gaussian elimination of off-diagonal elements
in the first row and the first column. So $Q_2$ is non-degenerate
with the determinant at least one and the eigenvalue of largest absolute value being bounded in terms of $M_0
R_0^{-1}$.

An apparent change of the analyticity domain of Hamiltonian
(\ref{wwham}) as far as the variables $\varphi$ are concerned is
easy to take into account; this will be done in Section 4. For now
let us assume that $H_\theta(\cdot, \varphi)$ is real-analytic for
$\varphi\in\T^n_\sigma$ for some $\sigma$. The only inevitable
analyticity loss so far has been $\delta_0$ in the application of
Lemma \ref{nfl}. The actions $\p$ live in a complex ball around
the origin, whose radius is ``as large as necessary'', provided
(\ref{sc}) is satisfied.

It is easy to see that
\[
{\Xi}_\theta\circ{L}_\psi={L}_\psi\circ{\Xi}_\theta,\;\;L^j_{2\pi}\circ{\Xi}_\theta
={\Xi}_\theta\circ L_{2\pi\theta}^{-j}\circ L_{2\pi}^j,
\]
where for $j\in\Z,$ $L^j_{2\pi\theta}$ is a canonical
transformation, corresponding to the base space diffeomorphism
$\boldsymbol
l^j_{2\pi\theta}:\,(x,\varphi)\rightarrow(x,\varphi+2\pi
j\theta)$. Then (\ref{trans}) gets modified to
\begin{equation}
{H}_\theta=
{H}_\theta\circ{L}^{-j}_{2\pi}\circ{L}^{2\alpha(j)}_\psi\circ
L^j_{2\pi\theta}={H}_\theta\circ{L}^{-2\alpha(j)}_\psi\circ{L}^{-j}_{2\pi}\circ
L^j_{2\pi\theta},
\label{transs}
\end{equation}
where ${\displaystyle \alpha(j)\=\left\{\begin{array}{cl}
  0, & j\mbox{ even}, \\
  1, & j\mbox{ odd}
\end{array}  \right.}$ is the parity of $j$. In other words (\ref{transs}) reads
\[\begin{array}{lll}
H_\theta(y,I,x,\varphi)&=&
H_{\theta}(y+2\alpha(j)\lambda\psi(x),I,x-2\pi j,\varphi+2\pi j\theta)\\
&=&
H_{\theta}(y-2\alpha(j)\lambda\psi(x-2\pi j),I,x-2\pi j,\varphi+2\pi j\theta).
\end{array}
\]
The manifold $W_{\tt t}$ is now represented by the zero section of
the bundle $T^*(\R\times \T^n)$, which contains unstable or stable
tori ${\cal T}_{j,{\tt t}}$ for respectively even or odd values of
$j$. All the unstable [stable] tori can be identified with one
another via the transformation $\Xi_\theta^{-1}$. In Fig. 1 the
tori ${\cal T}_{0,{\tt t}},{\cal T}_{\pm,{\tt t}}$ correspond to
$j=0,\pm1$ respectively.

\begin{figure}[ht]
\includegraphics{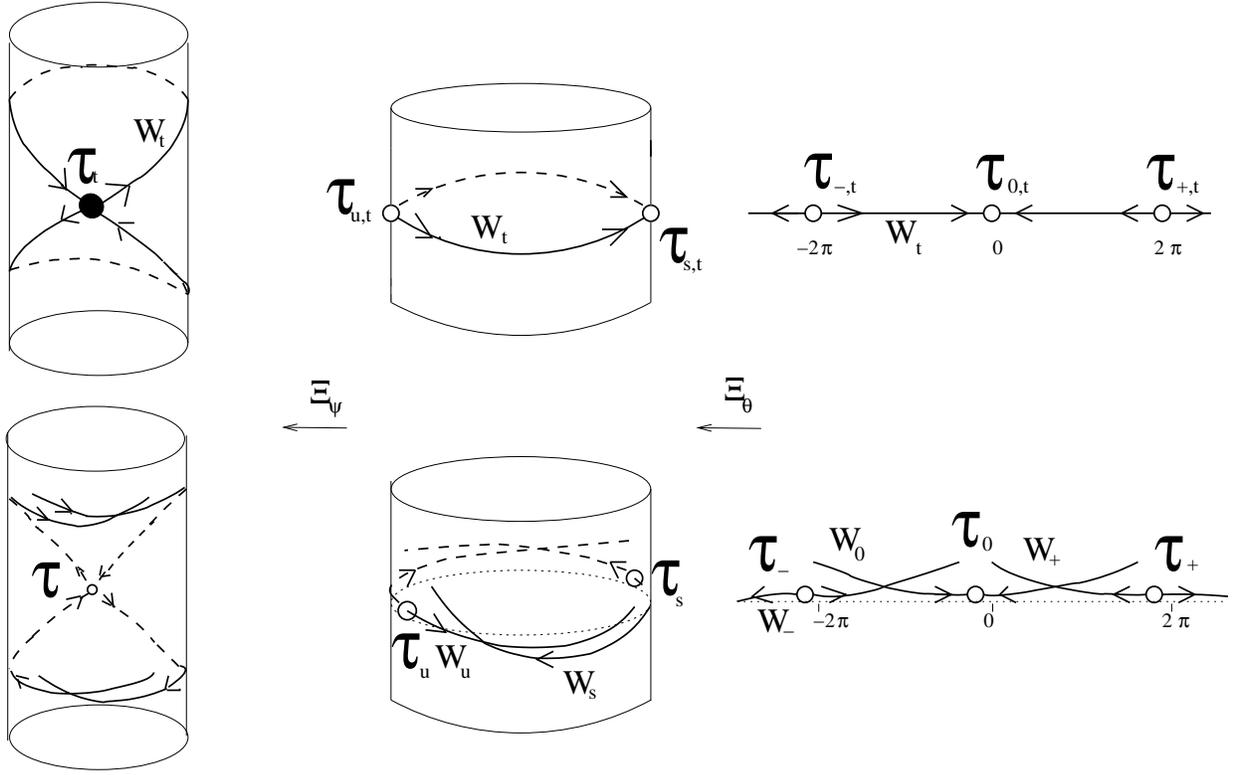}
\caption{Splitting scenario}
\label{fig1}
\end{figure}

Clearly, as it was the case with $H_\nu$, it suffices still
suffices knowing $H_\theta$ on the interval $x\in[0,2\pi)$ only.
The relation (\ref{transs}) applied to the truncated Hamiltonian
\[
H_{\theta,{\tt t}}=\lambda\psi(x)y + {y^2\over 2} + \spr{\omega_1}{I} +
{1\over2}\spr{\Theta I}{I}
\]
simply implies that if one writes
\[
\lambda\psi(x)y + {y^2\over 2}={y\over 2}(2\lambda\psi(x)+y),
\]
one sees the sputnik manifold $W'_{\tt t}$, where $y=-2\lambda\psi(x)$.
Naturally this manifold is the zero section for the Hamiltonian $H'_{\theta}\equiv H_{-\psi}\circ\Xi_\theta$.
This will
transform the latter expression to
\[
-\lambda\psi(x)y + {y^2\over 2}, \] which is tantamount to the shift $x\rightarrow x+2\pi$ in view of
$2\pi$-antiperiodicity of $\psi$. Similarly changing $y\rightarrow -2\lambda\psi(x)+y$ in the perturbation in
(\ref{wwham}) is tantamount to changing $x\rightarrow x+2\pi$ and $\varphi\rightarrow\varphi-2\pi\theta$.

It is convenient to treat $H_\theta$ as a multi-valued real-analytic function
on $T^*(\T'\times\T^n)$, whose branch is specified by fixing an even value of
$j$ in (\ref{transs}) and $j=0$ suffices for consideration. The branches differ by the shift of the angle $\varphi$
by an integer multiple of $4\pi\theta$. The splitting problem is well-posed for
a chosen branch of $H_\theta$ and the magnitude of splitting is clearly the
same on each branch. Technically, first one restricts $x$ to an interval ${\cal
I}_0={\cal I}_u$. Theorem \ref{mt} furnishes a Lagrangian manifold $W_0$ as a
graph over ${\cal I}_0\times\T^n$, containing an invariant torus ${\cal T}_0$
near $x=0$, for which it is the unstable manifold. $W_0$ is
described by a generating function $S_0$. Theorem \ref{mt} also results in the
stable sputnik manifold $W_0'$ of the torus ${\cal T}_0$, described by a generating function
$S'_0$, such that one-forms $dS_0$ and $dS'_0$ belong to the same cohomology
class $\xi=[dS_0]$.

Then $x$ is restricted to an interval ${\cal I}_{-}={\cal I}_s$
and one gets the Lagrangian manifold $W_{-}$ as a graph over
${\cal I}_{-}\times\T^n$ containing an invariant torus ${\cal
T}_{-}$ near $x=-2\pi$, described by the generating function
$S_{-}$. By (\ref{transs})  $S_{-}=S_0'\circ \boldsymbol l_{2\pi}
\circ\boldsymbol l_{-2\pi\theta}$. Moreover starting from the pair
$(S_0,S_-)$ (corresponding to $j=0,-1$) using (\ref{transs}) one
can define pairs of manifolds $(W_j,W_j')$ for all $j$, containing
invariant tori ${\cal T}_j$ ($W_j$ or $W'_j$ being respectively
unstable or stable manifolds for even or odd values of $j$
respectively) as graphs over $\boldsymbol l_{2\pi j}({\cal
I}_0)\times\T^n$, with generating functions $S_j=S_{j-2}\circ
\boldsymbol l_{-4\pi}\circ \boldsymbol l_{4\pi\theta}$, all
characterized by the same $\xi$. In particular, $S_+=S_-\circ
\boldsymbol l_{-4\pi}\circ \boldsymbol l_{4\pi\theta}$ corresponds
to $j=1$.

All the manifolds $W_j$ for even or odd $j$ are identified respectively with $W_u$ or $W_s$ via the
transformation $\Xi_\theta^{-1}$, it suffices to introduce the splitting function in the same way as it was
described above, identifying $\beta=\pm$ with $j=\pm1$ respectively. With the same notation one has
$\fS_\beta=S_0-S_{\beta}$, well defined on ${\cal I}_{\beta}\times\T^n$.

\medskip
\noindent
This completes the construction of the normal form, and calls for a structural stability theory for
Hamiltonians like $H_\theta$, which underlies the proof of Theorem \ref{main}. Another goal is to make this theory amenable to
the presence of the sputnik symmetry, in order to be able to conclude that $\xi=[dS_j],\,\forall j$. Both issues are studied at
length in the next section.

\label{ls}

\section{KAM theory on semi-infinite bi-cylinders over tori}
\label{hypkam}\setcounter{equation}{0}
As the forthcoming theory is self-contained, the notation in this section
may be occasionally different from the preceding sections. E.g. the function
$\psi(x)$ is introduced axiomatically, rather than by (\ref{tchi}). For
structural stability no symmetry properties of $\psi$ are required, which on
the other hand are
essential for the splitting problem. The details of the set-up arising in connection with the $2\pi$-antiperiodicity
of $\psi$ are not addressed until Section \ref{splitt}.

\subsection{Energy-time coordinates}\label{et}
\subsubsection*{Time-map}
Consider a  real-analytic function $\psi(x)$ for $x$ in a closed
real interval ${\cal I}$, for definity containing the points
$\pm\pi$ in the interior. Suppose $\psi(0)=0,\,\psi_x(0)=1,$
$\psi(x)\neq0$ on ${\cal I}\setminus\{0\}$ and is uniformly
bounded with its first two derivatives. Then $\psi(x)$ allows a
holomorphic extension into some closed rectangular domain ${\cal
D}\subseteq\C$, symmetric with respect to the real axis and
containing ${\cal I}$ together with a ball $\B_r$ at the origin
for some ${r}$, such that $\psi$ and its first two derivatives are
bounded inside ${\cal D}$ and apart from that $|\psi(x)|\geq{r}/2$
in ${\cal D}\setminus \B_{{r}}$. The pair $(\psi,{\cal D})$ is
regarded as fixed.

For $x\in{\cal D}\setminus\{0\}$ consider a map
$s$ from ${\cal D}$ into a Riemannian surface ${\cal S}$ of the logarithmic
type and its inverse $x$ as follows:
\begin{equation}
s:\;x\rightarrow\int_{\pi}^x{d\zeta\over
\psi(\zeta)},\;\;\;x=s^{-1},
\label{stime}
\end{equation}
as well as
\begin{equation}
\chi:{\cal S}\rightarrow {\cal D},\;\;\chi=\psi\circ x.
\label{chip}
\end{equation}
The maps $x,\,\chi$ can be represented by homonymous functions of
a complex variable $s\in\C$, which are $2{i}{\pi}$-periodic and
well defined in a family of semi-infinite strips about the rays
$\Re{s}\in(-\infty,T],\,\Im{s}=\pi j,\,j\in\Z$ for some $T$, see
Fig. 2. They are real-valued on the above rays and vanish
exponentially as $\Re s\rightarrow-\infty$. Without loss of
generality $T>1$.

Let us further consider only such values of the parameters $r,T$
that $0<{r}_{\psi}< r < 2{r}_{\psi}<1$, $1<T_{\psi}<T<2T_{\psi}$,
for some fixed pair $({r}_{\psi},T_{\psi})$ defined in terms of
$(\psi,{\cal D})$ and such that $\B_{2r_{\psi}}\subset{\cal D}$
and $T_{\psi}>-2\log r_{\psi}$. The functions $x(s),\,\chi(s)$ are
holomorphic in the half-plane $\Re{s}\leq-2T_{\psi}$. Further
estimates will ignore constants depending on the pair $(\psi,
{\cal D})$ as well as constants $n,\tau$ in the Diophantine
condition (\ref{dio}) by using the $\lesssim, \gtrsim$ and
$\asymp$ symbols in an obvious way.

For $\rho>0$ and $s\in \C$ let
\begin{equation}
\begin{array}{l}
\Lambda_{\te,\rho}= \{
\Re s \leq T,|\Im s|\leq\rho\}\cup\{
\Re s\leq -2T_{\psi}\} \cup  \{
\Re s \leq T,|\Im s -\pi|\leq\rho\},\\ \hfill \\
\Lambda^-_{\te,\rho}=-
\Lambda_{\te,\rho},\hspace{3mm}\hat\Lambda_{\te,\rho}=
\Lambda_{\te,\rho}\cap\Lambda^-_{\te,\rho}
\end{array}
\label{lbd}
\end{equation}
be further referred to as complex bi-strips. Their projections on
the union of the real axis and the line $\Im s=\pi$ will be
denoted by omission of the index $\rho$. The index $T$ may also be
omitted in qualitative argument. On the other hand,
$\Lambda_\infty$ stands for a pair of lines $\R\cup\R+i$ (the $+$
or $-$ sign henceforth having a priority over $\cup,\cap$). The
difference between the case of a finite $T$ and $T=\infty$ will be
emphasized. Also define a bounded one-strip (rectangle)
\begin{equation}
\Pi_{\te,\rho}=\{s\in\C:\,|\Re{s}|\leq T,\,|\Im{s}|\leq\rho\},
\label{pr}
\end{equation}
with the same index drop rules. Clearly $\hat\Lambda_{\te,\rho}=\Pi_{\te,\rho}\cup\,\Pi_{\te,\rho}+i\pi$.

For any real-analytic function $\tilde u$ on ${\cal D}$, the
composition $u=\tilde{u}\circ s$ returns real values for $\Im
s=0,\pi$, let's coin the term ``bi-real-analytic'' for that. With
the above notations for the domains, the functions $x(s),\chi(s)$
will be referred to as bi-real-analytic for
$s\in\Lambda_{T,\rho}$, for some $(T,\rho)$. For a function
$\psi$, given by (\ref{tchi}) one will naturally have
$\rho<{\pi\over2}$ in Section 4.

Suppose $s\in\Lambda$ and $h$ is a canonically conjugate momentum to $x$. Consider a canonical
transformation $\Xi_s$ from $T^*\Lambda$ into $T^*{\cal I}$ as follows
\begin{equation} \Xi_s:\;\left\{
\begin{array}{lll}
x&=&  x(s), \\
y&=&{h\over\chi(s)}.
\end{array} \right.
\label{time}
\end{equation}
where the maps $x,\chi$ have been defined by (\ref{stime},
\ref{chip}).

Let us extend the maps $s,x$ to maps $\boldsymbol s,\boldsymbol x$
between $\Lambda\times\T^n$ and ${\cal I}\times\T^n$ acting as the
identity on $\varphi\in\T^n$, in accordance with the general
convention of using bold symbols referring to the whole base
space. Extend accordingly the transformation $\Xi_s$ to
$\Xi_{\boldsymbol s}$, incorporating the pair $(I,\varphi)\in
T^*\T^n$. Let
\begin{equation}
{\cal C}_\te=\Lambda_\te\times\T^n,\;\;{\cal
C}_\te^-=\Lambda_\te^-\times\T^n,\;\; \;\;\hat {\cal C}_\te=\ {\cal
C}_\te\cap {\cal C}^-_\te
\label{cld}
\end{equation}
be referred to as bi-cylinders over tori, further just
``bi-cylinders''. In particular, ${\cal
C}_\infty=\Lambda_\infty\times\T^n$. In this section only
semi-infinite bi-cylinders ${\cal C}_\te$ will be dealt with.
Bi-infinite and bounded bi-cylinders ${\cal C}_\infty$ and $\hat
{\cal C}_\te$ will come into play in Section \ref{splitt}. In
qualitative argument, the index $T$, if finite (unlike $T=\infty$)
is often omitted further.

\subsubsection*{Analyticity domains}
Let us describe more precisely the analyticity domains for the map $s$ in order to further define the necessary
function spaces on them. Technical difficulties will arise from the fact that the bi-cylinders ${\cal C}_\te,
\,{\cal C}_\infty$ are not compact. E.g. a ``near-identity'' transformation $\ar$ of Alexandroff
compactification of ${\cal C}_\te$ should not necessarily preserve $\{s=-\infty\}$, i.e. the differential $d\ar$
may be unbounded. Similarly a Hamiltonian of general position on $T^*{\cal C}$ may be unbounded as
$s\rightarrow-\infty$, unless the momentum $h=0$. The reason is clearly because $ds(0)$ does not exist. In
a series of papers \cite{{CG},{Ga},{GGM}}, etc. these difficulties were overcome via improper integration
techniques.

\medskip\noindent
Analyticity domains of functions involved will be
characterized by positive parameter vectors $\fp$ as follows. Let
$\mathfrak{p}=(r,T,\rho,\sigma)\in\R^4_{++}$. Introduce partial order
$\mathfrak{p}'=(r',T',\rho',\sigma')\,\leq\,\mathfrak{p}$ if $r'\leq
r,\,T'\leq T,\rho'\leq \rho,\,\sigma'\leq\sigma$. If
$\mathfrak{p}'\leq\mathfrak{p}$ and $|\mathfrak{p}-\mathfrak{p}'|
\equiv\inf(r-r',T-T',\rho-\rho',\sigma-\sigma')>0$, write $\mathfrak{p}'<\mathfrak{p}$. Addition
of parameter vectors, as well as multiplication by positive real numbers is
defined component-wise, as well as the difference
$\mathfrak{p}-\mathfrak{p}'$ for $\mathfrak{p}'<\fp$. For
$\de\in\R_{++},\,\de<|\fp|$ the notation $\fp'=\fp-\de$ means subtracting $\de$
component-wise.
In the sequel the components and dimension of the parameter vectors $\mathfrak{p}$ may vary; $\fp$ can
incorporate $T=\infty$.

\medskip\noindent
Given $T$ such that both points $x:\,\Re s(x)=T$ (of opposite
signs) are in the interior of ${\cal I}$, one may want to be able
to describe the widest complex strip $\Lambda_{T,\rho}$ for some
$\rho$ such that the image $x(\Lambda_{T,\rho})$ be contained in
${\cal D}$. This can be done as follows, see Fig.
2.\begin{figure}[h]
\includegraphics{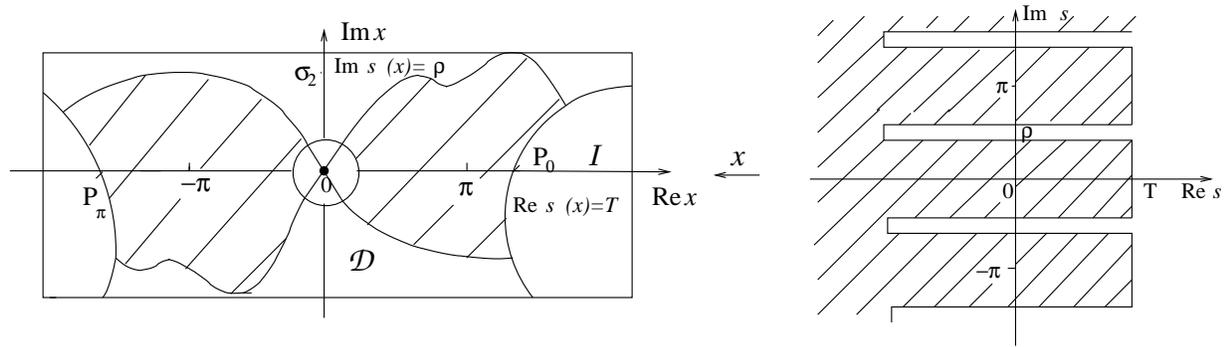}
\caption{The map $x$. The figure illustrates how the quantities
$\rho, \sigma_2$ can be determined (\ref{rho}) relative to the
pair $(\psi,{\cal D})$. The shaded region on the left, including
the circle around the origin is the domain ${\cal I}_{r,T,\rho}$
(\ref{xdomain}).}\label{fig2}
\end{figure}

\noindent
Let the level set $\Re s(x)=T$ for $x\in{\cal
D}$ intersect the real axis transversely at a pair of points
$P_0\in(\pi,\infty)$ and $P_\pi\in (-\infty,0)$. Let $\gamma_{0,\te}$ be a
connected component of the above level set containing the former point and
$\gamma_{\pi,\te}$ - containing the latter point (the two $\gamma$'s may coincide). The
points $P_0,P_\pi$ are connected to the origin by the level curves $\Im
{s}(x)=0,\pi$ respectively. Let $P_\zeta\in\gamma_{0,\te}\mbox{ or } \gamma_{\pi,\te}:\,\Im
{s}(P_\zeta)=\zeta$. For the
level set $\Im {s}(x)=\zeta,\,x\in{\cal D}$,
let $\gamma^*_\zeta$ be the connected component, whose closure contains the origin.
Let
\begin{equation}
\begin{array}{lll}
\rho_+ &=&\sup\{\zeta:\,0<\zeta\leq\pi,\, P_\zeta\in \gamma^*_\zeta, \, P_{-\zeta}\in
\gamma^*_{-\zeta}\}, \\ \hfill \\
\rho_- &=&\sup\{\zeta:\,0<\zeta\leq\pi,\, P_{\pi+\zeta}\in \gamma^*_{\pi+\zeta}, \,  P_{\pi-\zeta}\in
\gamma^*_{\pi-\zeta}\}, \\ \hfill \\
\rho&=&\inf(\rho_-,\rho_+).
\end{array}
\label{rho}
\end{equation}
In other words, the quantity $\rho$ simply shows for how long the
points $P_{0,\pi}$ can be moved  along the connected components of
the level curves $\Re s=T$, to which they belong, so that the
whole segment of the level curve of $\Im s(x)$, connecting them to
the origin remains contained in ${\cal D}$. Either
$\rho\in(0,\pi/2]$ or $\rho=\pi$, which corresponds to the case
when $\Re s=T$ is a simple closed curve contained in ${\cal D}$,
which together with its interior forms the image of the half-plane
$\Re s\leq T$ in ${\cal D}$. If $\rho<{\pi\over2}$ (the only case
of interest for the splitting problem) the equality $|\Im
s(x)|=\rho$ or $|\Im s(x)-\pi|=\rho$ can be achieved on four
different level curve segments $\gamma^*_\zeta$, where $\zeta=\pm
\rho$ or $\zeta-\pi=\pm \rho$ in Fig. 2. Marking these curve
segments simply as $\gamma^*_j,\,j=1,2,3,4$, define
\begin{equation}
\sigma_2=\inf_{j=1,2,3,4} \sup_{x\in  \gamma^*_j }|\Im x|.
\label{sone}
\end{equation}
In essence, these are the parameters $\rho,\sigma_2$ entering the main estimate (\ref{sir}) of Theorem
\ref{main}; the quantity $\sigma_2$ will not reappear until the end of Section 4. If one is willing to go into more detail, one should consider the above quantities as four-vectors
to account for each $\gamma^*_j$.

\medskip
\noindent
For $\fp=(r,T,\rho)$ denote
\begin{equation}
{\cal I}_{\fp}=\B_r\cup x(\Lambda_{\rho,T}).
\label{xdomain}
\end{equation}
For complex extensions of the bi-cylinders defined by (\ref{cld}) introduce the
notations
\begin{equation}
{\cal C}_{\te,\rho,\sigma}=\Lambda_{\te,\rho}\times\T^n_\sigma,
\label{cylinders}
\end{equation}
as well as ${\cal
C}^-_{\te,\rho,\sigma},\,{\cal
C}_{\infty,\rho,\sigma},\,
\hat{\cal C}_{\te,\rho,\sigma}$ analogous to (\ref{lbd}, \ref{cld}).

\subsubsection*{Function spaces}
Defined below are the necessary spaces of bi-real-analytic functions on
the bi-cylinders ${\cal C}$ as well as their maps. This is done simply
via the composition of real-analytic functions on or diffeomorphisms of ${\cal I}\times\T^n$
 with the bi-real-analytic  map $\boldsymbol x$. One needs
the following formalism in order to proceed toward an implicit function theorem
for structural stability of vector fields or Hamiltonians on ${\cal C}$ or
$T^*{\cal C}$ to be further used for exponentially small splitting estimates.
However, the theorem in question is interesting in its own right
as a ``non-compact'' version of KAM theory.

\medskip
\noindent Let ${\fB}^j({\cal D}),\,j\in\Z_+,$ be spaces of functions
real-analytic and uniformly bounded in ${\cal D}$, whose Taylor series at $x=0$
starts at order $j$; $j=0$ will be further omitted. With topology induced by the
supremum norm, ${\fB}^j({\cal D})$ are Banach spaces. Any $u\in
{\fB}^j({\cal D})$ can be represented as $u(x)=x^j v(x)$, where $v\in
{\fB}({\cal D})$, or alternatively as $u(x)=\psi^j(x) w(x)$, where $w\in
{\fB}({\cal D})$.  One can take the supremum of $|v|$ or $|w|$ for an
equivalent norm of $u$, the comparison constants depending on the pair $(\psi,{\cal D})$ only. For $\fp=(r,T,\rho)$ define the spaces ${\fB}^j({\cal I}_{\fp})$ in the same
way as ${\fB}^j({\cal D})$. As ${\cal I}_{\fp} \subseteq{\cal D}$,
clearly ${\fB}^j({\cal D})\subseteq{\fB}^j({\cal
I}_{\fp})\subset{\fB}^j({\cal I}_{\fp'})$ for $\fp'<\fp$. For
coherence with the forthcoming notation, let us write $\fB_{\fp}^j({\cal
I})$ instead of ${\fB}^j({\cal I}_{\fp})$.

If a function $u\in {\fB}_{\fp}^j({\cal I})$ has an extra analytic
dependence in $\varphi\in\T^n_\sigma$, $2\pi$-periodic in each
component of $\varphi$, one adds an extra component $\sigma$ in
the above parameter vector $\fp$ and writes
$u\in\bsitj,\,u=u(x,\varphi)$. Define the set $\bsj$ of all
holomorphic functions $u$ on ${\cal C}_{T,\rho,\sigma}$, such that
$u=\tilde u\circ \boldsymbol x$ for some $\tilde u\in \bsitj$.
E.g. consider a graph $y=D_x\tilde S(x)$ of a one-form in the
variables $(y,x)\in T^*{\cal I}$. If $S(s)=\tilde S[x(s)]$, then
the corresponding graph in the variables $(h,s)$ obtained via
(\ref{time}) is $h=dS(s)$, so $y=\chi^{-1}(s) dS(s)=O(1)$. I.e. a
$\varphi$-independent one-form $dS(s)$ over ${\cal C}$ vanishes
exponentially as $s\rightarrow-\infty$.

$\bsj$ is a closed subspace in the Banach space of all bounded holomorphic
functions on ${\cal C}_{T,\rho,\sigma}$ and thus a Banach space itself, with
the supremum norm $|\cdot|_{\fp}$. Note that if $u\in \bsj$, then the function
$u[s(x),\varphi]$ allows analytic continuation into the neighborhood $\B_r$ of
$x=0$, vanishing at $x=0$ to the $j$th order. This can be taken for an
independent definition of the spaces $\bsj$. Moreover, if $u(s,\varphi)\in
\bsj$, a multiplier $\chi^j(s)$ can be factored out, i.e.
\begin{equation}
u(s,\varphi)=\chi^j(s) v(s,\varphi),\;\;\;v\in \bs,
\;\;\;|v|_{\fp}\,\asymp\,|u|_{\fp}.
\label{dcp}
\end{equation}
Section 4 will deal with bi-real-analytic functions on bi-infinite and bounded bi-cylinders. To this effect, if
$u\in{\fB}_{r,T,\rho,\sigma}$ allows a uniformly bounded analytic continuation as $T\rightarrow\infty$, then
write $u\in \bsinf$ or $u\in\fB_{r,\rho,\sigma}({\cal C}_\infty)$.

Besides ${\fB}_{\te,\rho,\sigma}(\hat{\cal C})$ stands for the
space of bi-real-analytic functions $u(s,\varphi)$, which are
$2\pi$-periodic in each component of $\varphi$ and uniformly
bounded in the bounded bi-cylinder $\hat{\cal
C}_{\te,\rho,\sigma}$. Also let ${\fB}_{\te,\rho,\sigma}(\Pi\times
\T^n)$ be the space of real-analytic functions $u(s,\varphi)$,
which are $2\pi$-periodic in each component of $\varphi$ and
uniformly bounded in the bounded one-cylinder
$\Pi_{\te,\rho}\times\T^n_\sigma$, defined by (\ref{pr}). Clearly
${\fB}_{\te,\rho,\sigma}(\hat{\cal C})=
[{\fB}_{\te,\rho,\sigma}(\Pi\times \T^n)]^2$. With the
supremum-norm $|\cdot|_{\fp},$ each of the above spaces is a
Banach space. Component-wise supremum norm $|\cdot|_{\fp}$ or the
equivalent Euclidean norm $\|\cdot\|_{\fp}$ will be used for
vector functions.

\medskip
\noindent
For any  $u\in \bs$, there exists a
unique decomposition
\begin{equation}
u(s,\varphi)\,=\,u_0(\varphi)+u_{1}(s,\varphi), \hspace{3mm}\mbox{where}\hspace{3mm}u_0\in{\fB}_\sigma(\T^n),\;\;u_{1}\in {\fB}^1_{\fp}({\cal C}).
\label{dcmp}
\end{equation}
Using it, define the average $\langle u \rangle$ ``at infinity'' as
\begin{equation}
\langle u \rangle\,\=\,\int_{\T^n} u_0(\varphi)d\varphi. \label{aver}
\end{equation}
For $u\in\bs,$ its component $u_{1}$ satisfies an obvious
exponential estimate\footnote{Clearly not any real-analytic
function of $s$ vanishing at infinity at an exponential rate will
be a member of one of the above spaces. E.g. for $u(s)=se^s$, the
function $\tilde{u}=u[s(x)]$ is not analytic at $x=0$. } in ${\cal
C}_{\mathfrak{p}}$:
\[
|u_{1}(s,\varphi)|\,\lesssim \, e^{s} |u_{1}|_{\mathfrak{p}}.
\]

\medskip
\noindent Let us further describe the maps of the bi-cylinder
${\cal C}$ induced by real-analytic diffeomorphisms of ${\cal
I}\times\T^n$ after a change $x=x(s)$. Given
 $\fp=(r,T,\rho,\sigma)$, a sufficiently small $\de\in\R_{++}$ and $\fp'=\fp-\de$, let
\begin{equation}
\tilde{\ar}={\tt id}+\tilde{\br}:\;\q\rightarrow
\q+\tilde{\br}(\q),\;\;\q=(x,\varphi)\in{\cal D}\times
\T^n_\sigma,\;\;\tilde{\br}=(\tilde b,\tilde B)\in
[\bsit]^{n+1},\;\;|\tilde{\br}|_{\fp}\lesssim\de \label{primmap}
\end{equation}
be a smooth map of ${\cal I}_{r',T',\rho'}\times\T^n_{\sigma'}$
into ${\cal I}_{r,T,\rho}\times\T^n_{\sigma}$, well defined for a
small enough constant in the above estimate for
$|\tilde{\br}|_{\fp}$. It will always be assumed that
$\de<\de_{\psi}$, where the latter is ``small enough'' in terms of
the pair $(\psi,{\cal D})$. The natural norm for $\tilde{\br}$ is
the $C^1$-norm in $\bsitprime$, which is easy to estimate knowing
the $C^0$-norm on some intermediate space $\fB_{\fp''}(\It)$ with
$\fp'<\fp''\leq\fp$. Details regarding intermediate parameter
values will be mostly bypassed.

Let $\defpit$ be the set of all such diffeomorphisms and define
$\defp$ as the set of all maps
\begin{equation}
\ar:\;{\cal C}_{\fp'}\rightarrow{\cal C}_{\fp},\;\exists
\tilde{\ar}\in\defpit:\, \ar=\boldsymbol
s\circ\tilde{\ar}\circ \boldsymbol x. \label{cylmap}
\end{equation}
The analyticity indices can be dropped in the qualitative
argument. For  $\tilde{\ar}$ one can come up with a unique
representation $\tilde{\ar}=\tilde {\ar}_1\circ\tilde{\ar}_0$,
where $\tilde{\ar}_0:\,x\rightarrow x+ b_0(\varphi)$ acts on
$\varphi$ as the identity, while $\tilde{\ar}_1={\tt
id}+\tilde{\br}$ preserves $x=0$, i.e. $\tilde{\br}=(\tilde
b,\tilde B)\in{\fB}^1_{\fp}(\It)\times [{\fB}_{\fp}(\It)]^{n}$.
Then $\ar=\ar_1\circ\ar_0$, where the transformation $\ar_1$
preserves $\{s=-\infty\}$. Naturally one can write $\ar_1={\tt
id}+\br$, where $\br=(b,B)\in[\bs]^{n+1}$. Indeed for the change
$\varphi\rightarrow\varphi+\tilde B(x,\varphi)$ all one has to do
is to define $B(s,\varphi)\equiv\tilde B[x(s),\varphi]$. As far as
the variable $s$ is concerned, the change of the $x$-variable
corresponding to $\tilde{\ar}_1$ can be written as $x\rightarrow
x+\tilde b(x,\varphi)$ for $\tilde b\in{\fB}^1_{\fp}({\cal
I}\times\T^n)$, i.e. one can write $\tilde b=\psi(x)v(x,\varphi),$
with $v\in{\fB}_{\fp}({\cal I}\times\T^n)$. Thus given $x=x(s)$
one gets
\[
s\rightarrow s+\int_{x}^{x+\psi(x) v(x,\varphi)}{d\zeta\over\psi(\zeta)} =
s+ \sum_{j=0}^\infty [\psi(x)v(x,\varphi)]^{j+1}  {D^j_x\over (j+1)!}{1\over \psi(x)}
\equiv
s+b(s,\varphi),
\]
for some $b\in\bs$, with the norm $|b|_{\fp}\asymp
|v|_{\fp}\lesssim |\tilde{\br}|_{\fp}$. In particular the change
of $s$ under $\ar_1$ is asymptotically an identity as
$s\rightarrow-\infty$. The above expression can be viewed as a
$C^1$ functional from $\fB^1_{\fp}(\It)$ into $\bs$, mapping zero
into zero and whose differential is bounded away from zero in some
neighborhood of zero. Then by the inverse function theorem any
function $b\in\bs$ with $|b|_{\fp}\lesssim\de_{\psi}$ generates a
diffeomorphism $x\rightarrow x[s+ b(s,\varphi)]\equiv x+\tilde
b(x,\varphi)$, with $\tilde b\in {\mathfrak B}^1_{\fp}({\cal
I}\times\T^n)$.

For the transformation ${\ar}_0$ it's easy to see that writing
\[
 s\,\rightarrow \,a_0(s,\varphi)=a_0[s,b_0(\varphi)]=\int_{\pi}^{x(s)+b_0(\varphi)} {d\zeta\over \psi(\zeta)}
\]
is as far as one can get, as the series expansion analogous to the preceding
formula will not converge uniformly in $s$ for $\Re{s}\leq -2T_{\psi}$, i.e
near $x=0$.
\begin{remark} Unless $b_0\equiv0,$ the quantity $a_0[s,b_0(\varphi)]$ is neither
real-valued, nor continuous for real $s\leq -2 T_{\psi}$.
Continuity can be achieved by extending it to $\{s=-\infty\}$,
then the defining component $a_0$ or $\ar_0$ maps $\{\Re s\leq -2
T_{\psi},\,\Im{s}=0,\pi\}\cup \{s=-\infty\}$ into $\{\Re s\leq -
T_{\psi},\,\Im{s}=0,\pi\}\cup \{s=-\infty\}$ and the differential
$d\ar_0=d\boldsymbol s\circ d\tilde{\ar}_0 \circ d\boldsymbol x$
is unbounded. Further calculations will use the expressions
\begin{equation}
d\ar_0^{-1}\,=\,\left[\begin{array}{cc}
  {\psi[s(x)+b_0(\varphi)]\over\chi(s)} & -{d b_0(\varphi)\over\chi(s)} \\
  0 & {\tt id}_n
\end{array}\right],\hspace{3mm}
\psi[s(x)+b_0(\varphi)]=\chi(s) +[1+\chi(s)\eta_1(s)]b_0(\varphi)+\eta_2[s,b_0(\varphi)]b_0^2(\varphi),
\label{calc}
\end{equation}
where the quantities $\eta_1,\eta_2$ viewed as functions of
$(s,\varphi)$ are in $\bs$ by the assumptions on $\psi$.
\end{remark}As
one is interested in the coordinate changes
$\ar\in\mathfrak{D}_{\fp,\de}({\cal C})$ only  as far as their
action on functions from $\bs$ is concerned, they are naturally
represented by an element $\hat{\br}=(b_0,\br)$ of
$\bst\times[\bs]^{n+1}\equiv \bsz$, with the product topology and
vector supremum norm $|\cdot|_{\fp}$, the origin corresponding to
the identity transformation.

Then if $u\in\bs$ and $|\hat{\br}|_{\fp}\lesssim\de$, $u'=u\circ
\ar(\hat{\br})\,\in\,\bsprime$ with $\fp'=\fp-\de$. Moreover with
$\fp''=\fp-{1\over2}\de$ one can write
\begin{equation}
|u-u\circ \ar|_{\fp'}\lesssim |du|_{\fp''}
|\hat{\br}|_{\fp'}\lesssim {1\over \de}
|u|_{\fp}\,|\hat{\br}|_{\fp'}, \label{estm}
\end{equation}
by the Cauchy inequality.

Apart from $\ar=\ar_1\circ\ar_0$, the general form for the transformation $\ar$ can be also taken as
\begin{equation}
\ar(\hat{\br})=\ar_0(b_0)+\br:\;\left\{ \begin{array}{llclll}
  s & \rightarrow & a_0[s,b_0(\varphi)]&+&b(s,\varphi),  \\
  \varphi & \rightarrow & \varphi &+& B(s,\varphi).
\end{array}  \right.
\label{ar}
\end{equation}

\medskip
\noindent In order to deal with functions of $s$, which are
unbounded at infinity, let us introduce a function space
$\bsm\cong \bst\times\bs$ as follows, see (\ref{calc}):
\[ u(s,\varphi)\in {\fB}^{-}_{\fp}({\cal C})\;\mbox{ iff }\;
u(s,\varphi)={v(s,\varphi)\over\chi(s)},\;\;
v(s,\varphi)\in \bs.
\]
The norm on $\bsm$ is simply $|v|_{\fp}$. Since
one can write in the spirit of (\ref{dcmp})
$v(s,\varphi)=v_0(\varphi)+\chi(s)v_1(s,\varphi)$, with $v_1\in\bs$, then
\[
u(s,\varphi)={v_0(\varphi)\over\chi(s)}+v_1(s,\varphi),
\]
and $\sup(|v_0|_\sigma,\,|v_1|_{\fp})$ can be taken for the norm
$|u|_{\fp}$ as well. Also let $\bsg\equiv\bsm\times[\bs]^n$. An
element of this space describes a bi-real-analytic vector field on
the bi-cylinder $\cyl$. Clearly $\bsg\cong\bsz$, the elements of
the latter space representing the maps of $\cyl$. If $\g\in\bsg$
is a vector field and $\ar\in\defp$ then $\g\circ\ar$ is not in
$\bsgprime$, however. Indeed, as the result of the transformation
$\ar_0$ the quantity ${v_0(\varphi)\over\chi(s)}$ in the first
component changes to ${\displaystyle
{v_0(\varphi)\over\psi[x(s)+b_0(\varphi)]}}$ which can blow up for
a finite $s$. However $d\ar^{-1}\g\circ\ar$ corresponding to the
``new'' vector field does belong to $\bsgprime$, see (\ref{calc}).
Also, a simple calculation shows that in order to estimate the
norm for partial derivatives of a function $u\in
{\fB}^{-}_{\fp}({\cal C})$ one can still use the Cauchy formula
${\displaystyle |d u|_{\fp'}\,\lesssim\,\de^{-1}|u|_{\fp}.}$

\medskip
\noindent As far as Hamiltonian functions on $T^*\cyl$ are concerned, consider
the Banach space $\bskit$ (with the sup-norm) of bounded real-analytic
Hamiltonian functions on $T^*(\It)$, given by Taylor series with coefficients
in $\bsit$, uniformly convergent for the momenta
$\tilde{\p}=(y,I)$ inside $\B^{n+1}_\kappa,\,\kappa>1$.  Define the space $\bsk$ of
Hamiltonians on $T^*\cyl$ as the subset of holomorphic functions on $T^*\cyl$,
such that
\[
H\in\bsk\mbox{ iff }\,\exists\,\tilde H\in\bskit:\,H=\tilde
H\circ\Xi_{\s}.\] Thus the members of $\bsk$ are given by power
series in $({h\over\chi(s)},I)$, with coefficients in $\bs$.

\medskip
\noindent The final remarks on the notation are that sometimes, if
it is clear to which of the above spaces a function $u$ belongs,
the norm of $u$ may be referred to simply as $|u|$ rather than
$|u|_{\fp}$. If $u\in\bs$ and has bounded partial derivatives, the
notation $|u|_{1,\fp}$ will stand for the $C^1$-norm. The notation
$|u|_\infty$ will stand for the supremum norm of $u$, restricted
to the real values of all its variables, except $s$ which is
either real or $\Im s=\pi$.

\subsubsection*{Conjugacy problem}
In the formal framework developed above one can set up a conjugacy problem for a class of bi-real-analytic
perturbations of a constant vector field \begin{equation}
\x_0=\lambda{\partial\over\partial s}+\spr{\omega}{{\partial\over\partial\varphi}}
\label{vf}
\end{equation}
on the semi-infinite bi-cylinder ${\cal C}$, with $\lambda\in\R_{++}$ and a
Diophantine $\omega\in\di$. In the same way as (\ref{two}), the unperturbed
vector field can be taken slightly more general, i.e.
\[
\x_\theta=\lambda{\partial\over\partial
s}+\spr{\omega+\theta(s)}{{\partial\over\partial\varphi}},\] with an
angle-independent $n$-vector function $\theta(s)$, whose each component is a
member of the space $\bsone$. This case is
reducible to (\ref{vf}) after a change \[
 s=s',\;\;\varphi=\varphi'+\lambda^{-1}\int_{-\infty}^0
\theta(s+t)dt.\]

The question of structural stability of the vector field $\x_0$
under the group ${\mathfrak D}({\cal C})$ of bi-real-analytic maps
of ${\cal C}$ is roughly as follows: given a vector field
$\x=\x_0+\g$, where $\g\in\bsm$ (with the parameter vector
$\mathfrak{p}=(r,T,\rho,\sigma)$) and $|\g|_{\fp}$ is small
enough, does there exist a coordinate change $\ar\in  \defp$, i.e.
$\ar(\hat{\br})=\ar_0(b_0)+\br$ with $\ar_0=(a_0,{\tt id})$ and
$\hat{\br}=(b_0,\br) \in\bszprime$, ${\fp}'=\fp-\de$, such that
\begin{equation}
d\ar^{-1}\x\circ\ar\,=\,\x_0?
\label{question}
\end{equation}
The general answer to this question is {\em no}, as it is for the
torus, for one can take $\x=\x_0+\xii$, with a constant
$\xii\in\R^{n+1}$ ($\x_0$ being further identified with a constant
vector $(\lambda,\omega)\in\R^{n+1}$). Note that within the map
class $\defp$, the answer is {\em no} even if only the
``longitudinal'' component of the vector $\xii$ is nonzero, as the
scalings of the variable $s$ are outside this class. Hence,
conjugacy should be sought modulo $\xii\in\R^{n+1}$, asking for a
pair $(\ar(\hat{\br}),\xii)$, such that
\begin{equation}
d\ar^{-1}(\x_0+\g+\xii)\circ\ar\,=\,\x_0. \label{fldcj}
\end{equation}
The problem can be relatively easily shown to satisfy the input of
an implicit function theorem of Nash-Moser type, following the papers of Zehnder
\cite{Z1}, \cite{Z2}, who made further generalizations in the abstract set-up
to embrace the KAM theory with its small divisors. One essential
modification is that here one should deal with the differential operator
\begin{equation}
\Dlw\,\=\,\Diff
\,=\, \lambda D_s + D_\omega, \label{dhw}
\end{equation}
rather than just $D_\omega$, defined by (\ref{dom}). Auxiliary
results apropos of solvability of linear PDEs involving the
operator $\Dlw$ in the set-up of various function spaces
introduced earlier are presented in Appendix A.

However, estimates resulting from an application of the abstract
theorem are unsatisfactory for the analysis of the normal form
near a simple resonance (\ref{nfh}) with its hierarchy of orders
of magnitude and parameter dependencies. With extra scruple one
can benefit by quasi-linearity of underlying equations,
intermittent use of $C^1$ and $C^0$ estimates, similarly to the
classical KAM case \cite{Ru}. This is done in Section \ref{kam},
resulting in particular in Corollary \ref{cjt}.

\subsubsection*{Application to $H_\theta$}
Let the function $\psi$ be given by (\ref{tchi}). Then the domain
${\cal D}$ can be taken as a closed rectangle of some semi-width,
bounded from above by $\sigma_1$, symmetric with respect to the
interval ${\cal I}_0=[-2\pi+r,2\pi-r]$ of the real axis for some
$r<\inf(\sigma_1,1)$. The pair $(\psi,{\cal D})$ as well as the
parameter bounds $r_{\psi},\,T_{\psi},\de_{\psi}$ are well
defined, in particular one can ensure $T_\psi>1$, because the mean
value of $\psi$ over $\T'$ is zero. Then given $T\in
(T_{\psi},2T_{\psi})$ one can use (\ref{rho}-\ref{xdomain}) see
Fig. 2, to determine the constants $\rho$ and $\sigma_2$ as well
as the domain ${\cal I}_{\fp}$.

The application of
the transformation $\Xi_{\s}$ to the Hamiltonian
(\ref{wwham}), whose $x$-variable is restricted on ${\cal I}_{\fp}$
with the notation $(\p;\q)=(h,I;s,\varphi)$ results in the ``new'' Hamiltonian
$H_{\s}={H}_\theta\circ\Xi_{\s}\in\bsk$ as follows:
\begin{equation}
\begin{array}{llcccc}
H_{\s}(\p,\q)&=&\lambda h + \spr{\omega_1}{I} +
\spr{Q_{2}\tilde{\p}}{\tilde{\p}}&+& f_\theta[\x(\q)]+
\spr{\g_\theta(\tilde{\p},\x(\q))}{\tilde{\p}}\\ \hfill \\
&\equiv&H_{\lambda,\tt t}({\p},\q)&+&V_{\tt t}({\p},\q),
\end{array}
\label{shamt}
\end{equation}
where $\tilde{\p}=({h\over\chi(s)},I)=(y,I),$ and the function
$\chi(s)$ is defined by (\ref{tchi}, \ref{chip}). Or, including
the momentum-dependent part of $\g_\theta$ into the
``unperturbed'' Hamiltonian:
\begin{equation}
\begin{array}{llcccc}
H_{\s}(\p,\q)&=&\lambda h + \spr{\omega_1}{I}
+O_2(\tilde{\p};\q)&+& [f(\q)+
\spr{\g(\q)}{{\p}}]\\ \hfill \\
&\equiv&H_{\lambda}({\p},\q)&+&V({\p},\q),
\end{array}
\label{sham}
\end{equation}
where $f(s,\varphi)=f_\theta(x(s),\varphi),\; \g(s,\varphi)={\rm
diag}(\chi^{-1}(s),{\tt id}_n)\,{\g}_\theta(0;x(s),\varphi)$.
Namely $\g\in\bsg$ and the difference between $V$ and $V_{\tt t}$
is that the latter may contain super-linear terms in $\p$.
\medskip
\noindent For the sputnik Hamiltonian $H'_{\s}\equiv
H'_\theta\circ\Xi_{\s}\in\bsk$ one gets
\begin{equation}
H'_{\s}=H_{\s}\circ L^{-2}_{\chi}, \mbox{ with } L_{\chi}\equiv
L_\psi\circ\Xi_{\s}, \;\;\;{L}_\chi:\,(h,I,s,\varphi)\rightarrow(
h+\lambda\chi^2(s),I,s,\varphi). \label{lab}
\end{equation}
Then similarly to the two previous formulas
\begin{equation}
\begin{array}{llllll}
H'_{\s}(\p,\q)&=&H_{-\lambda,\tt t}({\p},\q)+V'_{\tt t}({\p},\q) \\ \hfill \\
&=&-\lambda h + \spr{\omega_1}{I}
+O'_2(\tilde{\p};\q)+[f'(\q)+
\spr{\g'(\q)}{{\p}}],
\end{array}
\label{shams}
\end{equation}
i.e $H_{-\lambda,\tt t}$ is the same as $H_{\lambda,\tt t}$, but
for the sign of the first term and $V'_{\tt t}= V_{\tt t}\circ
L_\chi^{-2}$, whereupon the terms $f'\in\bs$ and $\g'\in\bsg$
correspond to the zero and first order terms of the Taylor
expansion in $\tilde{\p}$; the rest of the expansion is absorbed
by the term $O'_2(\tilde{\p};\q)$. Just the same, $V'$ will denote
the expression in the square brackets in (\ref{shams}) above.

By Lemma \ref{nfl} one can take the radius $\kappa$ of convergence
of both Taylor series in $\tilde{\p}$ as large as necessary, for
instance $\kappa>2$ and for all $\tilde{\p}\in\B_\kappa$ assume
that uniformly over the (complexified) base space, the absolute
value of each eigenvalue of $D^2_{\tilde{\p}\tilde{\p}} H_{\s}$
and $D^2_{\tilde{\p}\tilde{\p}} H'_{\s}$ is uniformly bounded from
below by $1/2$ and from above by a constant times $M_0R_0^{-1}$.

\subsection{Hyperbolic KAM theorem}\label{kam}
This section develops a KAM-type approach to Hamiltonian systems in $T^*({\cal
C})$, such as the transformed simple resonance Hamiltonian (\ref{sham}). The
section contains the statement of Theorem \ref{mt}, the principal part of its
proof and a number of corollaries,  one of each is the solution to the
conjugacy problem (\ref{fldcj}). Yet the prototype of Theorem \ref{mt} can be
found in \cite{RW2} the theory exposed below is quite different in flavor, as
it avoids compactification of the base space and elucidates the connection with
other methods of study of manifolds asymptotic to invariant tori, developed for
instance in \cite{CG}, \cite{El}, \cite{Sa}, \cite{DG}. The theorem contains
significant technical improvements, allowing for a supposedly optimal parameter
dependence, necessary in order to make the theory applicable to a generic
simple resonance normal form and for the purpose of studying the lower bounds
for exponentially small splitting \cite{DGJS}, \cite{RW1}.

Theorem \ref{mt} per se represents an alternative to the traditional approach to whiskered tori, largely due to
Graff \cite{Gr}.  The present approach appears to be more natural for describing the
whiskers as semi-infinite cylinders over tori
globally\footnote{Graff's theorems apply to partially hyperbolic tori of all dimensions $\leq\,n$, but are in
essence local near a hyperbolic equilibrium, where the whiskers, being very ``short'' cylinders over tori, can
be naturally described by naive generating functions. The present theory takes advantage of the fact that in the
simple ``one-hyperbolic'' case discussed, the truncated normal form Hamiltonian (\ref{tnfh}) is integrable and
the generating function $S_{\tt t}$ (\ref{sep0}) is defined globally over a ``long'' cylinder; this necessitates a considerably different
analysis. The present approach seems to be extendable to the case of lower-dimensional tori, despite in the
latter case one certainly cannot hope to have a simplistic global description for the whiskers. Nevertheless,
one may try to consider only narrow strips thereof in tubular neighborhoods of transverse homoclinic orbits,
generically existing in the phase space of the non-integrable ``hyperbolic'' sub-system, corresponding to the
truncated normal form if one lets dim$\,x>1$ in (\ref{tnfh}).\label{hm}}, as it yields their representation via
generating functions, solving the Hamilton-Jacobi equation in the setting of the Banach space families,
introduced earlier. This representation with an extra sputnik property is ideally suited for describing the
splitting problem in Section 4. In addition, it can serve as an initial set-up for the variational construction
of orbits, shadowing the whiskers' intersections \cite{Be}.

Theorem \ref{mt} applied to  the Hamiltonian (\ref{sham}) implies
structural stability of its principal part, when the
momentum-subquadratic part of the perturbation $(f,\g)$ is small.
In case it is zero, the zero section of the bundle $T^*{\cal C}$
is an invariant Lagrangian manifold. Under the perturbation, this
manifold persists as a graph of a closed one-form $dS$ over ${\cal
C}$. Throughout the proof, a rapidly vanishing sequence  of
generating functions $\{S_j\}_{j\in\Z_+}$ is constructed, starting
from $S_0\equiv0$, such that the graph of the one-form
$d\sum_{j=0}^NS_j$ approximates the Lagrangian manifold in
question with increasing accuracy as $N\rightarrow\infty$. The
manifold itself is then described by the generating function
$S=\sum_{j=0}^\infty S_j$. In order to obtain the sequence
$\{S_j\}$, a KAM type iterative procedure is developed. A single
step of it is described by the Iterative lemma. Namely, for
$j\in\Z_{++}$ one constructs an ``affine'' canonical
transformation on $T^*({\cal C})$, i.e.

\begin{equation}
{\Psi}_j\,=\,{\Psi}_j(\ar_j,S_j):\ \left\{
\begin{array}{llllllll}
\q_{j-1}&=&\ar_j(\q_j),\\ \p_{j-1}&=&\ti{(d\ar_j)}\p_j+dS_j,
\end{array} \right.
\label{iter}
\end{equation}
with $(\p_0,\q_0)\equiv(\p,\q)$.

A closed one-form $dS_j$ is given by the generating function
$S_j(\q)=S_j(s,\varphi)=\spr{\xi_j}{\varphi}+\hat{S}_j(s,\varphi)$,
where $\xi_j\in\R^n$ is the $H^1({\cal C},\R)$ cohomology class
representative for $dS_j$, and the function $\hat{S}_j(s,\varphi)$
is a zero-form on ${\cal C}$, i.e. is $2\pi$-periodic in each
component  of $\varphi$. A single pair $(\ar_j,\hat{S}_j)\in
{\fD}_{\fp_{j-1},\de_j}({\cal C})
\times{\fB}_{\mathfrak{p}_j}({\cal C})$, where
$\mathfrak{p}'<\mathfrak{p}_{j'}<\mathfrak{p}_{j}<\fp_0\equiv\mathfrak{p}=(r,T,\rho,\sigma),$
the initial set of analyticity parameters,
$\de_j=|\fp_{j-1}-\fp_j|$ for $1\leq j<j',$ and
$\lim_{j\rightarrow\infty}\mathfrak{p}_j=\mathfrak{p}'$,
component-wise. The function $\hat{S}_j$ is chosen in order to
eliminate to the leading order a momentum-independent term $f$ in
the perturbation. The transformation $\ar_j$ is an approximate
(first order) solution the conjugacy problem (\ref{fldcj}). Its
existence is ensured by the specific choice of the quantity
$\xi_j$. Upon application of the transformation $\Psi_j$, the
principal part of the Hamiltonian (modulo a constant) picks up a
term equal to a small constant times $h$, thus slightly changing
the value of $\lambda$.

As $j\rightarrow\infty$, the
Hamiltonian vector field in question, restricted to the manifold $\p=dS$ thus becomes
conjugate to (\ref{vf}) (with a slightly changed value of $\lambda$) via a
canonical transformation $\Psi  =\Xi_1\circ\Xi_2\circ\ldots$.
Note that the classical Kolmogorov's theorem \cite{Ko} allows  a similar
geometric interpretation \cite{Z2} with the base space, of course being a
torus, rather than a bi-cylinder.

\subsubsection*{Set-up and statements of Theorem \ref{mt} and Iterative lemma}
The unperturbed Hamiltonian in (\ref{sham}) belongs to a certain class. Fix a
Diophantine  $\omega\in\di$ and define a class $\fN_{\omega}$ as follows.

\begin{definition}[Unperturbed Hamiltonian] A function $H_\lambda$ on $T^*{\cal C}$
belongs to the class $\fN_{\omega}$ if  modulo a constant, it can be represented as follows:
\begin{equation}
H_\lambda(h,I,s,\varphi)\,=\,\lambda h+\spr{\omega}{I}+O_2(\tilde{\p};\q), \label{templ}
\end{equation}
where $\lambda\neq0$ (further assumed positive) and
\begin{enumerate}

\item $O_2(\tilde{\p};\q)=O_2({h\over\chi(s)},I;\q)\,\in\bsk$,
{\rm for some parameter vector} $(\kappa,\mathfrak{p})$;

\item  $\exists\,R,M\in\R_{++},$
{\rm such that}  $\forall\,(\tilde{\p},\q)\in\B^{n+1}_\kappa\times{\cal
C}_{\mathfrak{p}},$ $\|\langle D^2_{\ie\ie} O_2(\tilde{\p};\q)\rangle^{-1}\|\leq
R^{-1}$ {\rm and} $\|D^2_{\tilde{\p}\tilde{\p}} O_2(\tilde{\p};\q)\|\leq M$.
\end{enumerate}
\label{nfdef}
\end{definition}

\noindent Given $H_\lambda\in\fN_{\omega}$, consider its small perturbation

\begin{equation}
H= H_\lambda+V,\;\;\;V(\p,\q)=f(\q)+\spr{\g(\q)}{\p},
\label{pt}
\end{equation}
where

\begin{equation}
f\in\bs,\;\g\in\bsg, \hspace{3mm}|f|_{\fp}\,\leq \,\mu,\;|\g|_{\fp}\,\leq
\,\mu\nu^{-1},\hspace{3mm}\mbox{for some} \;\;0\leq\mu<\nu\leq1.
\label{aspt}
\end{equation}
The parameter $\nu$ is further used to obtain the desired smallness condition
(\ref{munot}) generally indicating that the above described iterative
procedure allows larger upper bounds for the norm of $\g$ than the norm of $f$
in the perturbation $V$.

\begin{theorem}[Hyperbolic KAM theorem]
\label{mt} If $\mu$ is small enough,  there
exists a canonical transformation
\begin{equation}
\Psi  \,=\,\Psi   (\ar ,S ): \left\{
\begin{array}{llllllll}
\q&=&\ar(\q'),\\ \p&=&\ti{(d\ar )}\p'+dS ,
\end{array}
\right. \label{zed1}
\end{equation}
such that for any $\kappa'<\kappa,\,\mathfrak{p}'<\mathfrak{p}$
and some new parameter values $\lambda',R',M',$ different  from
$\lambda,R,M$ respectively by $O(\mu)$, one has $H \circ\Psi   \in
\fN_{\omega},$ with the new parameter set
$\{\kappa',\fp',\lambda',R',M'\}$ and $\de=|\fp-\fp'|$:

\begin{enumerate}
\item The transformation $\ar=\ar(\hat{\br})\in\defp$, with  $|\hat{\br}|_{1,\fp'}\,=\,O(\mu)$.
\item The one-form $dS$ is defined by the generating  function
$S (\q)=\spr{\xi }{\varphi}+\hat{S} (s,\varphi)$, with
$\xi \in\R^n,$ $\hat{S}\in\bsprime$, and $|\xi |,|\hat{S} |_{1,\fp'}\,=\,O(\mu)$.
\end{enumerate}
\end{theorem}
In the above non-technical formulation of Theorem \ref{mt},
``$\mu$ small enough'' means that it satisfies the following
smallness condition (\ref{munot}). The symbols $O(\mu)$ depend on
the parameter values from both the old and the new parameter sets.
Further without loss of generality, one can assume that the
quantities
$\delta=\sigma-\sigma',\de=|\fp-\fp'|,\lambda,R,M^{-1},|\omega|^{-1}\leq1$.
The exact estimates are summarized below. In applications, one or
more of them can turn out to be functions of a small parameter
$\eps$, and the magnitudes of the analyticity loss in the
variables $(s,\varphi)$ can differ considerably. Then the
following estimates can be adjusted if necessary, see e.g.
footnote \ref{da} below.

\medskip
\noindent {\bf Parameter statement of Theorem \ref{mt}}

\medskip
\noindent {\em Let
\begin{equation}
\begin{array}{lllllllllll}
\kappa'=\kappa-1,&\varsigma&=& \inf(\gamma\delta^{\tau},\lambda),&
&\eta&=&R\inf(M^{-1}\varsigma\de,\nu).
\end{array}
\label{eeta}
\end{equation}
There exists a constant $C$, depending only on $n,\tau,\psi,{\cal D}$, such that if
\begin{equation}
\mu\;\;\;\leq \;\;\; C^{-2}\eta^2\;\;\lesssim\;\;
(R/M)^2\de^2[\inf(\varsigma,\nu)]^2, \label{munot}
\end{equation}
the following estimates hold:
\begin{equation}
\begin{array}{ccclllllllll}
|\hat{S}|_{\fp'}&\leq &C\mu \varsigma^{-1}, &
|d\hat{S}|_{\fp'}&\leq& C\mu(\varsigma\de)^{-1}, &
|d{S}|_{\fp'}&\leq& C\mu\eta^{-1}, \\ \hfill \\
|\hat{\br} |_{\fp'}&\leq & C\mu (\eta\varsigma)^{-1},&
|d\hat{\br} |_{\fp'}&\leq& C\mu (\eta\varsigma\de)^{-1},\\ \hfill \\
\lambda^{-1}|\lambda'-\lambda|&\leq&C\mu (\eta\lambda)^{-1}, \\ \hfill \\
M^{-1}|M'-M|, \,R^{-1}|R'-R|&\leq&C\mu (\eta\varsigma\de)^{-1}.
\end{array}
\label{transf}
\end{equation}
}

\medskip
\noindent
Theorem \ref{mt}
can be cast into the abstract generalized Newton method framework \cite{Z1},
\cite{Z2}. However, in order to obtain the desired parameter dependencies, a direct proof is given.
The main tool is furnished by
the following lemma, fulfilling a single Newton's iteration.

\begin{lemma}[Iterative lemma]
For a Hamiltonian (\ref{templ}-\ref{aspt}) with a parameter set
$\{\kappa,\fp,\lambda,R,M,\mu,\nu\}$,
if $\mu$ is small enough (condition (\ref{small}) below), there exists
a canonical transformation
\begin{equation}
\Psi \,=\,\Psi  (\ar,S): \left\{
\begin{array}{llllllll}
\q&=&\ar(\q'),\\ \p&=&\ti{(d\ar)}\p'+dS,
\end{array}
\right. \label{zed2}
\end{equation}
such that for any $\kappa'<\kappa,\,\mathfrak{p}'<\mathfrak{p}$,
and some new parameter values $\lambda',R',M',$ different  from
$\lambda,R,M$ respectively by $O(\mu)$, one has for
$\de=|\fp-\fp'|$:
\begin{enumerate}
\item The transformation $\ar=\ar(\hat{\br})\in\defp$ is such that $|\hat{\br}|_{1,\fp'}\,=\,O(\mu)$.
\item The one-form $dS$ is defined by the generating function $S(\q)=\spr{\xi}{\varphi}+\hat{S}(s,\varphi)$, with
$\xi\in\R^n,$ $\hat{S}\in\bsprime,$ and $|\xi|,|\hat{S}|_{1,\fp'}\,=\,O(\mu)$.
\item The Hamiltonian ${\displaystyle
H\circ\Psi  =  H_{\lambda'}+V',}$ with $H_{\lambda'}\in \fN_{\omega}$
satisfies Definition \ref{nfdef} and (\ref{aspt}) with
a parameter set $\{\kappa',\fp',\lambda',R',M',\mu',\nu'\}$,
where $\mu'=O(\mu^2)$ and $\nu'$ is independent of $\mu$.
\end{enumerate}
\label{il}
\end{lemma}
Clearly, the quantities $\Psi,\ar,S,$ etc. in the Iterative lemma are not the
same as their homonyms in Theorem \ref{mt}; it should not cause confusion. In
the non-technical formulation above, the quantities $\mu',\nu'$ and the symbols
$O(\mu)$ depend on the parameter values from both the old and the new parameter
sets. They are further specified as follows.

\medskip
\noindent {\bf Parameter statement of Iterative lemma}

\medskip
\noindent {\em Let $\kappa-\kappa'<1$ and $\varsigma,\eta$ be computed via formulae (\ref{eeta}).
There exists a constant $C$, depending only on $n,\tau,\psi,{\cal D}$, such that if
\begin{equation}
\mu\;\;\;\leq \;\;\; C^{-2}\eta\varsigma\de(\kappa-\kappa'),
\label{small}
\end{equation}
the following relations hold:
\begin{equation}
\begin{array}{rclcccrclllll}
|\hat{S}|_{\fp'}&\leq &C\mu \varsigma^{-1}, &
|d\hat{S}|_{\fp'}&\leq& C\mu(\varsigma\de)^{-1}, &
|\xi|&\leq& C\mu\eta^{-1}, \\ \hfill \\
|\hat{\br} |_{\fp'}&\leq & C\mu (\eta\varsigma)^{-1},&
|d\hat{\br} |_{\fp'}&\leq& C\mu (\eta\varsigma\de)^{-1},\\ \hfill \\
\lambda^{-1}|\lambda'-\lambda|&\leq&C\mu (\eta\lambda)^{-1},\\ \hfill \\
M^{-1}|M'-M|, \,R^{-1}|R'-R|&\leq&C\mu (\eta\varsigma\de)^{-1},\\ \hfill \\
\nu'&=&\inf(M^{-1}\varsigma\de,\nu), & \mu'&\leq & C^2\mu^2 \eta^{-2}.
\end{array}
\label{newmu}
\end{equation}
}

\medskip
\noindent The proof of the Iterative lemma is given in the next
section. Detail of the proof is  important in order to justify the
following remarks, concerning Theorem \ref{mt}. Once one accepts
the Iterative lemma, the rest of the proof of Theorem \ref{mt}
becomes a routine iteration scheme, given in Appendix B.

\subsubsection*{Corollaries and remarks}

Theorem \ref{mt} essentially states that the zero
section $\p'=0$ is an invariant Lagrangian manifold for the
Hamiltonian $H\circ \Psi\,\in\,\fN_{\omega}$. Due to the
``affineness'' of the transformation $\Psi  (\ar,S)$, in the ``old'' coordinates
$(\p,\q)$ this manifold is a section of $T^*{\cal C}$, given by
\begin{equation}
W(s,\varphi)\,=\,\overline{\{(h,I,s,\varphi)\in \R^{n+1}\times {\cal C}\,:\;h=D_s \hat S_1(s,\varphi),\,
I=\xi+D_\varphi\hat S(s,\varphi)\}},
\label{lagrange}
\end{equation}
where $\hat{S}\in\bsprime$, i.e. $\hat S(s,\varphi)=\hat
S_0(\varphi)+\hat S_1(s,\varphi)$ in the sense of (\ref{dcmp}).
Recall that ${\cal C}=\Lambda\cup\T^n$ is a bi-cylinder, i.e. $\Im
s=0$ or $\pi$.

The manifold $W(s,\varphi)$ has been compactified by incorporating an invariant torus ${\cal T}$ corresponding to
$\{s'=-\infty\}$.  The flow on ${\cal T}$ is conjugate to rotation with the
frequency $\omega$. With $S'=S\circ\ar$ and $\q'=(s',\varphi')$, another parameterization for $W$ is
\begin{equation}
W(\q')\,=\,\overline{\{(\p,\q):\;\q=\ar(\q'),\,\p=\ti{[d\ar(\q')]}\,dS'(\q'),\,\q'\in{\cal
C}\}}, \label{lagrangep}
\end{equation}
and get ${\cal T}$ as
\begin{equation}
{\cal T}=\bigcap_{T\leq 0} W[(s',\varphi')\in
\Lambda_T\times\T^n].
\label{KAMtorus}
\end{equation}
\begin{remark} In the sequel ${\cal T}$ will be often described as located ``near'' $s=-\infty$: this verbiage refers
precisely to the representations (\ref{lagrange}-\ref{KAMtorus}) and the underlying
estimates.\end{remark}The representations (\ref{lagrange}) constitute the basis for the splitting
analysis in Section \ref{splitt}. For easier cross-reference let us recap the
above as a corollary.
\begin{corollary}[Hamilton-Jacobi equation]
Let $H_0(0,\q)=0$. The function $S(\q)=\spr{\xi}{\varphi}+\hat{S}(\q)$ satisfies the Hamilton-Jacobi equation
on $\bsprime:$
\begin{equation}
H(\partial_{\q}S(\q),\q)\,=\,c_0,
\label{hjea}
\end{equation}
with a real $c_0$ bounded by (\ref{c0}).
\label{hjc}
\end{corollary}
Let us now look back at the conjugacy problem (\ref{fldcj}). It
corresponds precisely to the case $f\equiv0$ in the perturbation
(\ref{pt}). Then the parameters $M,R,\nu$ are redundant (can be
all set to $1$ in the estimates) and the unperturbed Hamiltonian
$H_\lambda$ can be thought momentum-linear. Then given the
perturbation $\g\in\bsg$, the aim of the conjugacy problem is to
find a transformation $\Psi(\ar)=\Psi(\ar,0)$ to conform with
(\ref{zed1}) that is with $S\equiv0$, as well as a constant
$\xii\in\R^{n+1}$ (which is not unrelated to the vector
$(\lambda-\lambda',\xi)\in\R^{n+1}$) such that $(H_\lambda+
\spr{\g+\xii}{\p})\circ \Psi=H_{\lambda'}$. Then the proof of
Lemma \ref{il} and Theorem \ref{mt} can be straightforwardly
adjusted to yield the following corollary.

\begin{corollary}
There exists a constant $C(n,\tau,\psi,{\cal D})$, such that
if $\g\in\bsm$, with
\[
|\g|_{\mathfrak{p}}\;\leq\; C^{-2}\varsigma\de,
\]
there exists  a pair $(\ar,\xii)\in({\mathfrak
D}_{\fp,\de}\times\R^{n+1})$, such that the map $\ar=\ar(\hat{\br})$ effects (\ref{fldcj}) and
\begin{equation}
\begin{array}{rcl}
|\xii|,|\hat{\br}|_{1,\mathfrak{p}'}& \leq & C|\g|_{\mathfrak{p}}(\lambda\de)^{-1}.
\end{array}
\label{mucon}
\end{equation}
\label{cjt}
\end{corollary}

\noindent Besides, if $H=\tilde H\circ \Xi_{\s}$, then the
transformation $\Xi_{\s}\circ\Psi\circ \Xi_{\s}^{-1}$ effects the
structural stability of the normal form
\[
\tilde H_\lambda(y,I,x,\varphi)\;=\;\lambda\psi(x)y + \spr{\omega}{I} +
O_2(\tilde{\p};x,\varphi)
\]
under small perturbations, proving the existence of an invariant manifold which is a graph over
the variables $(x,\varphi)\in[-2\pi+r,2\pi-r]\times\T^n$, containing an
invariant torus near $x=0$ \cite{RW2}.

The next corollary is quite obvious with respect to $\tilde H$. It claims that
if the latter is perturbed by a pair $(\tilde f,\tilde{\g})$ such that $\tilde
f, d\tilde f, \tilde{\g}$ all vanish at $x=0$, then the invariant torus at
$x=0$ satisfies the Hamilton equations and does not move: in particular the
constants $\xi,c_0$ in Theorem \ref{mt} and Corollary \ref{hjc} should be zero.
For instance, this is the case in Arnold's example \cite{A1} (the
Hamilton-Jacobi formalism for such a degenerate perturbation was developed in
\cite{Sa}). This fact can also be established by going through the proof of Lemma \ref{il}
and will be used further to claim Corollary \ref{sput}, essential for the
splitting problem.

\begin{corollary}[Degenerate perturbation]\label{dp}
Suppose $f\in\bstwo$ and $\g=(g,G)\in\bs\times[\bsone]^n$. Then apropos of the
transformation (\ref{zed1}) one has $\xi=0$, $S\in\bstwoprime$, $b_0\equiv0$,
$\br=(b,B)\in\bs\times[\bsoneprime]^n$ and $c_0=0$ in Corollary \ref{hjc}.
\end{corollary}

\noindent
The following remarks address yet more technical issues.
\begin{remark}
{\tt (Estimates)} The smallness condition (\ref{munot}) for $\mu$
appears to be optimal as far as the parameter dependencies are
concerned: if $\lambda,\nu=1,\,\de=\delta$ it reproduces the
standard KAM theorem optimal smallness condition, see e.g.
\cite{Po1}. The use of an extra parameter $\nu$ has been essential
here to express that the order of magnitude of $\g$ is inherently
somewhat greater than that of $f$, as far as the perturbation is
concerned. This fact can  obstruct accessibility of the condition
(\ref{munot}) if one pursues Kolmogorov's approach to KAM theory,
see e.g. \cite{DG}. Similar (standard KAM) estimates resulting
from the general abstract implicit function theorem machinery are
also worse \cite{Z1}, \cite{Z2}. Under the assumptions of
Corollary \ref{dp} the estimates of the Parameter statement of
Theorem \ref{mp} shall be modified as follows. Apart from $\xi=0$,
one should use $\varsigma=\lambda$ and formally set $R=1$ in all
the estimates. Finally, if $\lambda<0$ in Definition \ref{nfdef}
then obviously $|\lambda|$ should substitute $\lambda$ in the
estimates.
\end{remark}

\vspace{-5mm}
\begin{remark}
{\tt (Local\footnote{This is certainly not true for ``large'' perturbations:
each hyperbolic manifold (unstable or stable) should have a counterpart (stable or unstable) or {\em
sputnik}, see the coming Corollary \ref{sput}.} uniqueness and parameter
dependence)} Given $H_\lambda$ and a small perturbation $(f,\g)$ obeying the
smallness condition (\ref{munot}),  the pair $(\ar,dS)$ is unique (as it is in
standard analytic KAM theorem). Indeed, the unique solution $u(v)$ of the PDE
in Proposition \ref{mp}, provides the right inverse of the operator $\Dlw$,
which is also its left inverse, guaranteeing uniqueness, see \cite{Z1}. In
other words, local uniqueness follows from the uniqueness of PDE solutions (modulo a constant) in
Appendix A. Similarly, in the case of a continuous (e.g. on the pair $(f,\g)$)
or real-analytic dependence of $H$ in an extra parameter (e.g. $\mu$), the pair
$(\ar,S)$ retains the same type of dependence in the parameter. Local
uniqueness is indispensable for the splitting problem to be well-posed.
\end{remark}

\vspace{-5mm}
\begin{remark}
{\tt (Other settings)} It is known in KAM theory that the
non-degeneracy assumption in Definition \ref{nfdef} allows many
variations \cite{Ru1}. Theorem \ref{mt} can be adapted to these
settings in the same way as the standard KAM theorem. For instance
$H_\lambda$ in (\ref{pt}) can be only linear in the actions $I$,
provided that the perturbation $V$ does not depend on $I$ either,
the so-called ``isochronous'' \cite{CG} case. In this case the
transformation $\ar$ acts on the $\varphi$-variables as the
identity, and in the smallness condition (\ref{munot}) one can
certainly set $R=1$. One can pass smoothly (so-called
``twistless'' case \cite{Ga}) to the isochronous case from the
set-up of Definition \ref{nfdef} by introducing an extra
parameter, which would multiply all the terms containing $I$ in
the Hamiltonian, except $\spr{\omega}{I}$. This can be verified by
examining the proof of Lemma \ref{il} in the same fashion as it is
shown for the standard KAM theorem \cite{Lo1}.

As this paper has been motivated by the need for a general theory
for exponentially small splitting near resonances, and the main
estimate (\ref{expsmall}) of the ensuing Theorem \ref{estimate}
will not persist in the $C^r$ category (even to the first order of
perturbation theory) the latter setting has not been considered.
However there seems to be no obstruction to generalizing Theorem
\ref{mt} for $C^r$ Hamiltonians (with $r$ large enough) using
standard smoothing techniques \cite{{Z1},{Po1}}.
\end{remark}

\vspace{-3mm} \noindent Suppose the Hamiltonian $H_\lambda$ has
extra structure admitting an unperturbed sputnik, similarly to
(\ref{shamt}-\ref{shams}).
\begin{assumption}[Sputnik on ${\cal C}_T$]
There exists a canonical transformation
\begin{equation} L_*=L_*(\ar^*,S^*):\;\left\{
\begin{array}{llllllll}
\q&=&\ar^*(\q_*),\\ \p&=&\ti{(d\ar^* )}\p_*+dS^*(\q),
\end{array}
\right. \label{sptloc}
\end{equation}
where $\ar^*$ is a map of ${\cal C}_T$, $\ar^*={\tt id}+\br^*$,
with $\br^*\in[\bs]^{n+1}$, while $dS^*$ is exact, with
$S^*\in\bstwo$, such that both Hamiltonians $H$ and $H'\equiv
H\circ L_*$ satisfy the input of Theorem \ref{mt} with the same
analyticity parameters $\kappa,\fp$ and
equivalent\footnote{Equivalent here is meant in the same sense as
the equivalence for orders of magnitude or norms, that is up to a
constant factor. Clearly $\lambda,R,M$ in Theorem \ref{mt} can be
just bounds, rather than the actual values of non-degeneracy
parameters. I.e. a change $\lambda\rightarrow-\lambda$ or
$M\rightarrow 2M$ is inconsequential. Besides the transformation
$L_*$ may in principle entail some extra analyticity loss of the
order $\de$. Necessary amendments are easy to make for a concrete
example.\label{feq}} non-degeneracy parameters $\lambda,R,M$ and
smallness parameters $\mu,\nu$. Moreover, see
(\ref{templ}-\ref{aspt}), suppose $H_\lambda\circ
L_*=H_{\lambda'}+V',$ where $H_\lambda,H_{\lambda'}\in\fN_\omega$,
$|\lambda'|\asymp|\lambda|$, while the perturbative term $V'$
satisfies the assumptions of Corollary \ref{dp}. \label{sploc}
\end{assumption}

\begin{remark} Regarding the normal form Hamiltonian (\ref{nfh}) alias
(\ref{shamt}-\ref{shams}), the sputnik transformation
$L_*=L^{-2}_\chi$ arises from the fact that the lower separatrix
is a graph of over the upper one for the truncated normal form
Hamiltonian (\ref{tnfh}). To this effect, the action of $L_*$ is
$H_{\lambda,{\tt t}}\circ L_*= H_{-\lambda,{\tt t}}$, natural for
symplectic flows. So in this case $\ar^*={\tt id}$. However
further in Section 4, the diffeomorphism $\ar^*$ will incorporate
a shift of the variable $\varphi$ by
$2\beta\pi\theta,\,\beta\in\{+,-\}$ (and of $s$ by $-\beta i\pi$),
see (\ref{transs}) and (\ref{transss}) to come. Assumption
\ref{sploc} is a generalization, which may be useful for instance
in case of higher multiplicity resonances when there is more than
one characteristic exponent $\lambda$, see also footnote
\ref{hm}.\end{remark}Note that as $\ar^*$ acts as the identity on
the set $\{s=-\infty\}$ (possibly causing reparameterization of
the angles $\varphi$ in terms of $\varphi_*$) it is not only $S^*$
but also $S^*\circ \ar^*\in\bstwo$, and the unperturbed invariant
tori for $H$ and $H'$ can be identified. As $H'$ satisfies the
input of Theorem \ref{mt}, there exist a canonical transformation
$\Psi'(\ar^1,S^1)$ and an invariant manifold $W'$ containing an
invariant torus ${\cal T}'$, versus the transformation
$\Psi(\ar,S)$ and the pair $(W, {\cal T})$ for $H$, described by
(\ref{lagrange}, \ref{KAMtorus}). Note that the claim
$\Psi'=L^{-1}_*\circ\Psi$ would be false, as the latter
transformation would still make the manifold $W$ the zero section,
rather than $W'$.

\begin{corollary}[Sputnik]\label{sput} Under
Assumption \ref{sploc}, let $\Psi(\ar,S):\,H\circ \Psi\in
\fN_\omega$, $\Psi'(\ar^1,S^1):\,H'\circ \Psi'\in \fN_\omega$
according to (\ref{zed1}) Theorem \ref{mt}. Then the closed
one-forms $dS$ and $dS^1$ belong to the same cohomology class
$\xi\in H^1({\cal C},\R)$. Both $dS$ and $dS'\equiv dS^* +
d[S^1\circ (\ar^*)^{-1}]$ belong to the same cohomology class and
satisfy the Hamilton-Jacobi equation (\ref{hjea}) of Corollary
\ref{hjc} for $H$, on the same energy level $c_0$. The closure
$W'$ of the graph of the form $dS'$ intersects the manifold $W$
defined by (\ref{lagrange}) for the form $dS$ at the torus ${\cal
T}$ defined by (\ref{KAMtorus}), i.e ${\cal T}'={\cal T}$.
\end{corollary}
{\tt Proof:} Let $H_2=H\circ\Psi\circ L_*$. Since $H\circ\Psi\in
\fN_\omega$, by the explicit form (\ref{sptloc}) of $L_*$, under
the assumptions on the pair $(\ar^*,S^*)$, the Hamiltonian
$H_2=H'\circ (L_*^{-1}\circ\Psi\circ L_*)$ satisfies the
conditions of Corollary \ref{dp}, for $\mu$ small enough.

Indeed, $(H\circ\Psi)\circ L_* =
(H_\lambda+[c_1h+O_2(\tilde{\p};\q)])\circ L_*,$ where $c_1$ is a
constant and the term in square brackets is $O(\mu)$. Then the
assumption on $L_*$ (in particular its preservation of
$\{s=-\infty\})$ implies the previously made statement, namely
$$H_2(\p_*,\q_*)= \lambda_2 h_*+\spr{\omega}{I_*}+O_2(\tilde{\p}_*;\q_*)+V_2(\p_*,\q_*),$$ where
$\lambda_2=\lambda'+O(\mu)$ and
$\tilde{\p}_*=({h_*\over\chi(s_*)},I_*)$.

In fact, the pair $(f_2,\g_2)$, comprising the perturbation $V_2$
above, which has arisen as the result of substitution of
(\ref{sptloc}) into the Hamiltonian $H\circ\Psi$ can be bounded in
terms of $|d\hat{\br}|_{\fp'}\de^{-1}$ times a constant,
determined by on $\Psi_*$, see (\ref{transf}). As $S^*(\q_*)$ is
proportional to $\chi^2(s_*)$, the quantities
$f_2(\q_*),\g_2(\q_*)$ behave as $s_*\rightarrow -\infty$ amenably
to Corollary \ref{dp}.

Then for $\mu$ small enough\footnote{Note from the remark on
estimates that the smallness condition for the norm of
$(f_2,\g_2)$ to guarantee the transformation $\Psi_2$ is in fact
more relaxed than the right hand side of (\ref{munot}). However,
one should not worry about the precise smallness condition here:
Theorem \ref{mt} warrants the existence of the transformations
$\Psi,\Psi'$ as long as (\ref{munot}) is satisfied, and this is
all one needs, plus local uniqueness.} there exists a canonical
transformation $\Psi_2=\Psi_2(\ar_2,S_2)$ in the form
(\ref{zed1}), such that $H_2\circ\Psi_2\in\fN_\omega$. The flow of
the Hamiltonian $H_2(\p_*,\q_*)= \lambda_2
h_*+\spr{\omega}{I_*}+\ldots$ contains an invariant Lagrangian
manifold ${\displaystyle
W_*=\overline{\{\p_*=dS_2(\q_*),\,\q_*\in{\cal C}\}}}$ near the
zero section $\{\p_*=0,\,\q_*\in{\cal C}\}$, with the generating
function $S_2\in\bstwoprime$, containing an invariant torus ${\cal
T}$, corresponding to the limit as $s_*\rightarrow-\infty$ and
$\p_*=0$. The same flow also contains the invariant manifold
$\displaystyle W=\overline{
\{\p_*=-dS^*[\ar^*(\q_*)],\,\q_*\in{\cal C} \} }$, i.e. the
closure of the zero section for $H\circ\Psi$ (set $\p=0$ in
(\ref{sptloc})), containing the same torus. Note that both
$[dS^*]=[dS_2]=0$.

The claim now is that $W_*=W'$, by local uniqueness. Indeed,
$H\circ\Psi\circ L_*\circ\Psi_2 = H\circ L_* \circ(
L^{-1}_*\circ\Psi\circ L_*\circ \Psi_2)\in\fN_\omega$ and the
transformation in the parentheses is near identity. On the other
hand, $(H\circ L_*)\circ\Psi'\in \fN_\omega$, so one must have
$\Psi'=L^{-1}_*\circ\Psi\circ L_*\circ \Psi_2$, for $\mu$ small
enough. Both the left and the right hand side of the latter
identity must have the (unique) form (\ref{zed1}) in terms of the
pair $(\ar^1,S^1)$, which by Theorem \ref{mt} is well defined as
long as $\mu$ satisfies (\ref{munot}). Then $W$ and $W'$ clearly
belong to the same energy level.

This essentially completes the proof. One may notice that there is
a natural semidirect product structure that the canonical
transformation composition induces on pairs $(\ar,S)$. Namely in
the ``old'' coordinates $(\p,\q)$ the manifold $W_*$ is
represented in terms of the generating function
$S'=S^*+S^1\circ(\ar^*)^{-1}$, on the other hand equal
$S+S^*\circ\ar^{-1}\circ S_2\circ (\ar\circ\ar^*)^{-1}$, which (as
$[dS^*]=[dS_2]=0$) clearly implies $[dS]=[dS']$, i.e
$S(s,\varphi)=\spr{\xi}{\varphi}+\hat{S}(s,\varphi)$ and
$S'(s,\varphi)=\spr{\xi}{\varphi}+\hat{S}'(s,\varphi)$, with the
same $\xi$. However, despite $\hat S$ and $\hat S'$ are both
elements of $\bsprime$, one cannot claim that necessarily
$S_0(\varphi)=S'_0(\varphi)$ in the sense of the decomposition
(\ref{dcmp}). $\Box$
\begin{remark} {(\tt Notation)}
The final remark in this section is that the theory developed above can
obviously be applied to the restriction of the Hamiltonian $H(\cdot,s)$ to
``one strips'' $s\in\{\Re s\leq T, |\Im s|\leq\rho\}\cup\{\Re
s\leq-2T_{\psi}\}$ or $s\in\{\Re s\leq T, |\Im s-\pi|\leq\rho\}\cup\{\Re
s\leq-2T_{\psi}\}$. Let us reserve the notation $\Lambda_{\beta,\te,\rho}$ with
$\beta=+$ or $-$ respectively for the lower or the upper one of the strips above. The
indices $(T,\rho)$ can be omitted in accordance with the notational convention
introduced in Section \ref{et}; also $T=\infty$ can be used. In
the same fashion the subscript $\beta$ may be added to the notations ${\cal C},
H, \Psi, S,\ar,$ etc. Clearly, as far as the application of Theorem \ref{mt} to
the restrictions $H_\beta$ over ${\cal C}_\beta$ is concerned, it results in
the same constants $c_0,\xi$ as well as the $s$-independent components $\hat
S_{0}(\varphi)$ of the generating functions $\hat S_\beta$, which then
represent analytic continuations of one another for $\beta\in\{+,-\}$. These
$\beta$-notations will be used in Section 4.
\end{remark}

\subsection{Proof of  Iterative lemma}
The proof of Lemma \ref{il} consists of several steps inherent in
KAM-type theorems \cite{{Ko},{Z2}}. In the classical case, Proposition
\ref{rus} plays the key role. Here these are Propositions \ref{mp}-\ref{mp1}.
The notations $\Psi,\ar,S$, etc. in this section pertain to the formulation of
the Iterative lemma and its parameter statement, rather then Theorem \ref{mt}.

\medskip
\noindent First, notice that any ``affine'' transformation
$\Psi[\ar(\hat{\br}),S]$ described by (\ref{zed2}) (where
$\ar\in\defp$, with the norms $|\hat{\br}|_{1,\fp'}$ and
$|S|_{1,\fp'}$ small enough, say $O(\mu)$) will act in such a way
that $H\circ\Psi\in\bskprime$ for $H\in\bsk$, despite the
unboundedness as $s\rightarrow-\infty$. This follows from the
definitions of the spaces $\bs,\defp,\bsk$ in Section \ref{et} and
can be verified directly. In particular, if $S\in\bsprime$ and has
bounded partial derivatives, the quantity ${D_s
S(s,\varphi)\over\chi(s)}$ is bounded. Thus a substitution
$(h,I)\rightarrow (h,I)+d S$ into (\ref{pt}) will not affect its
structure as a Taylor series in $\tilde{\p}=({h\over\chi(s)},I)$.

Also note that ${h\over \chi(s)}\circ\Psi[\ar(\hat{\br}),0] =
{h[1+v(s,\varphi)]\over \chi(s)},$ where $v\in\bsprime$ and is $O(\mu)$. This
is equivalent to an earlier made statement that
$d\ar^{-1}\,\g\circ\ar\in\bsgprime$ for $\g\in\bsg$.  Combining it with the
fact that by (\ref{estm}) a single Taylor coefficient, member of $\bs$ in the
Hamiltonian will only change by $O(\mu)$ as a result of the action of the
transformation $\ar$, one can see that $H\circ\Psi\in\bskprime$ if $\mu$ is small enough.

\subsubsection*{Homological equation}
The quantities $(\ar,S)$ are to eliminate the perturbation $V$ to the leading
order. In order to do so they should solve approximately the ``homological
equation'':
\begin{equation}
H_\lambda\circ\Psi(\ar,S)\,-\,H_\lambda\,+\,V\,-\,c_0\,-\,c_1 h\,=\,O_2(\tilde{\p})\,+\,O(\mu^2),
\label{hml}
\end{equation}
which is only possible for some specific values of the constants
$c_0,c_1$ to be found. Writing (\ref{hml}) out in essence requires
only a direct substitution of (\ref{zed2}) into (\ref{pt})
followed by an estimate for the ``remainder'' $O(\mu^2)$. In order
to get the equation for $\hat{S}$, it is enough to plug
$\p\rightarrow\p+dS=\p+(0,\xi)+d\hat{S}$ into the Hamiltonian
(\ref{pt}) assuming $S=O(\mu)$.

Furthermore if $\ar=\ar_0(b_0)+\br$ as in (\ref{ar}), with $\hat{\br}=(b_0,\br)=O(\mu)$, a calculation using
(\ref{calc}) and $d\ar^{-1}=({\tt id}+d\ar_0^{-1}d\br)^{-1}(d\ar_0)^{-1}$  yields
\begin{equation}
d\ar^{-1}\left[\begin{array}{c}
  \lambda \\ \omega
\end{array}\right]=\left[\begin{array}{c}
  \lambda \\ \omega
\end{array}\right]-{1\over\chi}\left[\begin{array}{c}
  (-\lambda+D_\omega)b_0  \\ 0
\end{array}\right] - \left[\begin{array}{c}
  \Dlw b - \eta_1 b_0 \\ \Dlw B
\end{array}\right] + O(|\hat{\br}|_{\fp'}^2)
\label{calc1}
\end{equation}
where the function $\eta_1(s)\in\bs$ is determined solely by
$\psi$. As the constant vector $(\lambda,\omega)$ formally belongs to $\bsg$, the result of having
multiplied it from the left by $d\ar^{-1}$ is a vector-function from
$\bsgprime$, whose leading order is given above.

It is convenient to think of $\hat{\br}\in\bszprime\cong\bsgprime$ as an element of
the latter space, formally writing
\begin{equation}
\begin{array}{rll}
\hat{\br}&=& \left[\begin{array}{c}{b_0\over\chi}+b\\B\end{array}\right], \\ \hfill
\\
\Dlw\hat{\br}&=&{1\over\chi}\left[\begin{array}{c}
  (-\lambda+D_\omega)b_0  \\ 0
\end{array}\right]+\left[\begin{array}{c}
  \Dlw b - \eta_1 b_0 \\ \Dlw B
\end{array}\right],
\end{array}
\label{calc2}
\end{equation}
as the latter expression appears in (\ref{calc1}). In other words,
in order to do the estimates throughout the rest of the proof, one
can set $d\ar^{-1}= {\tt id}-d\hat{\br}$, although $s\rightarrow
s+{b_0(\varphi)\over\chi(s)}$ would not be a legitimate
transformation for the $s$-variable, defined by (\ref{cylmap}).

Then one ends up having a pair of first order linear PDEs: one for
$\hat{S}$ and one for $\hat{\br}$. The first PDE is amenable to
Proposition \ref{mp}, the second one to Propositions \ref{mp} and
\ref{mp1} together, alias Proposition \ref{mpp}. It requires the
appropriate choice of constants $c_0$ and $c_1$ respectively, as
well as the one-form $dS$ cohomology class representative $\xi\in
\R^n$ in (\ref{zed2}). The latter quantity $\xi$ is chosen to
ensure that the right hand side in the equation for $B$ has a zero
average in the sense of (\ref{aver}) enabling the choice of the
constant $\boldsymbol c\in\R^{n+1}$ in Proposition \ref{mpp}
simply as $\boldsymbol c=(c_1,0)$ (owing to the non-degeneracy
assumption on the quadratic part of $H_\lambda$, see Definition
\ref{nfdef}). Thus (\ref{hml}) is equivalent to the following
system of equations

\begin{equation}
\begin{array}{lll}
c_0 &=& \langle f\rangle + \spr{\xi}{\omega},\\
\Dlw \hat{S} &=& \langle f\rangle-f,\\
 \langle\tilde G\rangle&=&0,\\ \Dlw \hat{\br}  &=& \g_1 + D^2_{\p\p}H(\p,\q)\vert_{\p=0}
\boldsymbol \xi-\boldsymbol c,
\end{array}
\label{invert}
\end{equation}
where
\begin{equation}
\boldsymbol \xi=(0,\xi),\hspace{3mm}\g_1 = \g + D^2_{\p\p}H(\p,\q)\vert_{\p=0} d\hat{S};
\hspace{3mm}\g_1=(g_1,G_1). \label{rhs}
\end{equation}

Following Kolmogorov \cite{Ko} one starts solving (\ref{invert}) with the  equation  for $\hat{S}$, then
finds $\xi$ to satisfy the penultimate one,
then solves the last equation, the constants $c_{0,1}$ being determined along
the way.

The equations for the quantities $\hat{S},\,\hat{\br} $ are
clearly amenable to Propositions \ref{mp}, \ref{mpp}. The norm of
the solution, according to these propositions, is simply estimated
by $\varsigma^{-1}$ times the norm of the right-hand side.
Applying the propositions results in analyticity loss. As a matter
of fact, one encounters the analyticity loss six times along the
way: solving the equation for $\hat{S}$, evaluating the
derivatives, solving the equation for $\hat{\br}$, evaluating the
derivatives, making sure that not only $\ar\in\defp,$ but the
$C^1$ estimate for the norm of $\hat{\br}$ is valid throughout the
maximum range of $\ar$, and finally inverting it. Hence, strictly
speaking one should introduce five intermediate spaces between say
$\bs$ and $\bsprime$, and the parameters $\delta,\de$ should be
scaled by factor 6. These standard steps are bypassed, and all the
estimates for $\hat{S}$ and $\hat{\br},$ as well as their
derivatives, no matter that they may also be valid in some
intermediate (smaller) spaces, are all written in the target
spaces $\bsprime$ and $\bsgprime\cong\bszprime$ right away. The
scaling of the analyticity loss parameters is absorbed into the
constant $C$ in (\ref{small}, \ref{newmu}). The estimates follow
from Propositions \ref{mp}-\ref{mpp}, Definition \ref{nfdef}, and
the bounds (\ref{aspt}). Here are the details.

First, by Proposition $\ref{mp}$ and the Cauchy inequality, $\hat{S}$ and its
partial derivatives are in $\bsprime$, with the estimates\footnote{ One can go slightly
more subtle estimating the derivative $D_s \hat{S}$, as it does not depend on
$\hat S_0$, where $\hat{S}=\hat{S}_0+\hat{S}_1$ in the sense of  the
decomposition (\ref{dcmp}). The norm of $\hat S_1$ however,  is estimated by
Proposition \ref{mp} without any small divisors as $\mu\lambda^{-1}$. Then one
can take a minimum of the following two estimates. One is to apply the Cauchy
inequality, acquiring a factor $\de^{-1}$. The other is to deduce from the
equation itself that $|D_s \hat{S}|\leq
\lambda^{-1}(|\omega||D_\varphi\hat{S}|+\mu).$ The same thing can be done
further estimating the $C^1$-norm of $\hat{\br}$. So a $C^1$-estimate for the
solution of the equation $\Dlw u=v$ on ${\cal C}$ can be obtained by dividing
the norm of $v$ by
$\inf[\varsigma\delta,\,\lambda\sup(\de,\lambda|\omega|^{-1}\delta)]$ rather than
$\de$, causing a straightforward modification of (\ref{newmu}) and
(\ref{transf}).\label{da}}
\begin{equation}
|\hat{S}|_{\fp'}\,\lesssim \,\mu \varsigma^{-1},\hspace{3mm}
|d\hat S|_{\fp'}\,\lesssim\,\mu(\varsigma\de)^{-1}. \label{dshat}
\end{equation}
Then the
quantity $\g_1$ in (\ref{rhs}) belongs to the space $\bsgprime$ (in
fact, any intermediate space between
$\bsgprime$ and $\bsg$) with the estimate
\[
|\g_1|_{\fp'}\lesssim \mu[M(\varsigma\de)^{-1}+\nu^{-1}].
\]
So the following expression is well
defined:
\[
\xi \;= \; -\,\langle\,D^2_{\ie\ie}H(\p,\q)\vert_{\p=0}\rangle^{-1} \langle
G_1\rangle.
\]
In order to estimate it, note that $\xi\in\R^n$ depends only on the right-hand
side as a real function of $\q$. Thus,
\begin{equation}
|\xi|\;\lesssim \;\mu
R^{-1}\left(M[\inf(\gamma,\lambda)]^{-1}+\nu^{-1}\right),\hspace{3mm}|dS|_{\fp'}\;\lesssim\;\mu\eta^{-1},
\label{ds}
\end{equation}
where in the first estimate the bounding constant depends on $(\sigma,\rho)$;
this is not the case in the second, rougher estimate, see (\ref{eeta}) for the
formula for the quantity $\eta$. This gives the first three of the estimates
(\ref{newmu}), as well as
\begin{equation}
|c_0|\,\lesssim\,    \mu
|\omega|R^{-1}\left(M[\inf(\gamma,\lambda)]^{-1}
+\nu^{-1}\right)\leq \mu|\omega|\eta^{-1}, \label{c0}
\end{equation}
the bounding constant in the first estimate depending on
$(\sigma,\rho)$, but not in the second estimate. Being finite, the
constant $c_0$ is further dropped.

The upper bound on the value of $c_1$ is then obtained from the last equation in (\ref{invert}) after using Proposition
\ref{mpp}, with the norm of the right hand side in $\bsgprime$ bounded by $\mu \eta^{-1}$:
\begin{equation}
|c_1|\,\lesssim\,   \mu \eta^{-1}. \label{c1}
\end{equation}
Then the $C^0$ norm of $\hat{\br}$ (as an element of $\bszprime$
or of its representation (\ref{calc2}) as an element of
$\bsgprime$) is bounded in terms of $\mu (\eta\varsigma)^{-1}$ and
the $C^1$ norm - in terms of $\mu (\eta\varsigma\de)^{-1}$, which
also ensures $|\hat{\br}|_{\fp'}\lesssim\de$, i.e. $\ar\in\defp$,
as well as its inverse if one takes a big enough constant in
(\ref{newmu}).

Finally, the presence of the additional multiplier
$(\kappa-\kappa')$, in the smallness condition (\ref{small}) is to
ensure that the image of $\B^{n+1}_{\kappa'}$ in (non-canonical)
momenta $\tilde{\p}=({h\over\chi(s)},I)$ under the map
$\Psi(\ar,S)$ is contained in $\B^{n+1}_{\kappa'}$, which is
tantamount to requiring
$|\hat{\br}|_{1,\fp'}\lesssim\kappa-\kappa'$. The smallness
condition (\ref{small}) by itself guarantees that $|c_1|<\lambda$.

\subsubsection*{Analysis of transformed Hamiltonian}
\label{ran} It remains to estimate the term $O(\mu^2)$ in (\ref{hml}). Let us
assume the smallness condition (\ref{small}) for $\mu$ and use the above
obtained bounds for $S$ and $\hat{\br}$. Let $C$ be large enough, say $100(n+1)$ times the
bounding constant for all the inequalities in the preceding section.

\begin{enumerate}
\item From (\ref{calc1}) it's easy to see that ${\tt
id}+d\hat{\br}$ is essentially a ``stretch factor'' for
non-canonical momenta $\tilde{\p}$. Thus for the growth and
non-degeneracy parameters $M',R'$ of the Taylor series
$H\circ\Psi\,\in\bskprime$ one has ${\displaystyle M^{-1}|M'-M|,
\;R^{-1}|R'-R|\,\lesssim\,|\hat{\br} |_{1,\fp'}}$. Together with
(\ref{c1}) it gives the corresponding estimates of (\ref{newmu}).
The contribution from the ``shift'' by $dS$ in (\ref{zed2}) is
negligible, as the non-degeneracy assumptions in Definition
\ref{templ} are global in $\tilde{\p}\in\B^{n+1}_\kappa$.

\item
The  ``new'' perturbative momentum-zero-order  term $f'$ is formed by several
contributions. The first one comes from the momentum-super-linear part of $H$ and is bounded
by $M|S|_{1,\fp'}^2$. The contribution from the linear terms is bounded by
$|c_1||\hat S|_{1,\fp'}+|\g|_{\fp}|S|_{1,\fp'}$. Note that
as $D_s \hat S\in\bsoneprime$, the ``new''
momentum-independent term $f'$ is in $\bsprime$ indeed. The final contribution is $f\circ\ar-f$,
with a bound $\mu|\hat{\br}|_{1,\fp'},$ by (\ref{estm}). Combining it with (\ref{ds}) and (\ref{c1}) yields
\[
|f'|_{\fp'}\,\leq\, C^2\mu^2\eta^{-2}.
\]

\item
The ''new'' perturbative  term ${\g}'$ in the first order in the momentum has a component, coming from the
momentum-super-linear part of $H$, bounded by
$M|\hat{\br}|_{1,\fp'}|S|_{1,\fp'}$.
Another contribution comes from the acquired term $c_1 h$; its norm can be
bounded by $|c_1||\hat{\br}|_{1,\fp'}$. Finally, the remainder $d\ar^{-1}(\g\circ\ar -\g)$
has to be taken into account, with the bound
$\mu\nu^{-1}|\hat{\br}|_{1,\fp'},$ by (\ref{estm}). The first contribution clearly dominates the second
one, and one can write
\[
|{\g}'|_{\fp'}\;\leq\;{1\over2}C^2\mu^2 [ M\eta^{-2}(\varsigma\de)^{-1} +
\nu^{-1}\eta^{-1}(\varsigma\de)^{-1}]\leq C^2\mu^2 \eta^{-2}
\sup[M(\varsigma\de)^{-1},\nu^{-1}].
\]
The last pair formulas complete the set of estimates (\ref{newmu})
and  the proof. The final remark to make here is that formally
setting $d\ar^{-1}={\tt id}-d\hat{\br}$ with
$\hat{\br}\in\bsgprime\cong\bszprime$ does not bring in extra
error. E.g. it can be taken precisely for the sought quantity, to
which after $\hat{\br}$ has been determined one can unambiguously
match a transformation $\ar\in\defp$ as long as
$|\hat{\br}|_{\fp'}$ is small enough. $\Box$
\end{enumerate}

\section{Splitting problem}\setcounter{equation}{0}\label{splitt}
This section contains the principal part of the proof of Theorem \ref{main}.
The theorem follows from the analytic splitting theory in $T^*({\cal
C}_\infty)$, developed further on the basis of the main results of the previous section.

\subsection{Preliminaries}
\subsubsection*{Energy-time coordinates}
Let us start out with the necessary additions to the set-up in the
beginning of Section \ref{et}.  One still uses formulae
(\ref{stime}, \ref{time}) for the energy-time coordinates $(h,s)$.
However there is extra structure underlying the splitting problem.
Namely, suppose the function $\psi(x)$ introduced in Section
\ref{et} and determining the transformation $s$ is
$2\pi$-antiperiodic (as it is in (\ref{tchi}) due to reversibility
of the truncated normal form (\ref{tnfh})). Then $\psi(x)$ is
$4\pi$-periodic, and the domain ${\cal D}$ for the $x$-variable is
a complex extension of $\T'=\R/4\pi\Z$. Technically, assume that
it contains a pair of balls of radius $r\in(r_{\psi},2r_{\psi})$
centered as $x=0$ and $x=2\pi$ and that outside these balls in
${\cal D}$ one has $|\psi(x)|\geq r/2$. In addition, without loss
of generality one can assume that
\begin{equation}
\mbox{P.V.}\int_{-\pi}^\pi {d\zeta\over \psi(\zeta)}=0,
\label{pv}
\end{equation}
where P.V. indicates that the integral is taken in the principal
value sense\footnote{Otherwise the lower limit of integration
$\pi$ in the defining formula (\ref{stime}) should be substituted
by some $a\in(0,2\pi)$ and the integral in (\ref{pv}) shall be
taken from $a-2\pi$ to $a$. Such an $a$ always exists by
continuity, $2\pi$-antiperiodicity and positivity on $(0,2\pi)$ of
$\psi(x)$.}.

Rewrite the $2\pi$-antiperiodicity property (\ref{ant}) as
$l_{2\pi}\circ\psi=\iota\circ\psi$, where $\iota:\,s\rightarrow
-s$ is sign inversion. Extend the latter to a diffeomorphism
$\io:\, (s,\varphi)\rightarrow (-s,\varphi)$.  As the image of the
map $s$ defined by (\ref{stime}) acting on ${\cal
D}\setminus\{0,2\pi\},$ one can simply consider an imaginary
circle $\C/2i\pi\Z$, as all the functions of $s$ we are dealing
with here are $2i\pi$-periodic. Otherwise a branch of $s$ can be
fixed by drawing a branch cut in ${\cal D}$ as a ``semicircle'' in
$\T'$, connecting the points $x=0$ and $2\pi$, but not containing
$x=\pi$. The inverse map $x$ is represented by a homonymous
function of a complex variable $s$, which is $2{i \pi}$-periodic
and analytic in the bi-infinite bi-strip
$\Lambda_{\te,\rho}\cup\Lambda^-_{\te,\rho}$, see (\ref{lbd}). For
now, $4\pi$-periodicity of the function $\psi$ does not allow one
to distinguish the values of $x(s)$ modulo $4\pi$. The half-width
$\rho$ and the parameter $\sigma_2$ can be still defined by
(\ref{rho}, \ref{sone}) with the choice of, say $T=T_{\psi}\geq
-2\log r_{\psi}$, i.e. the level set $\Re s=T$ will be contained
inside the ball of radius $2r_{\psi}$ centered at $x=2\pi$. In
other words,
\[
\rho=\sup\{\zeta>0: \mbox{ both curves $\Im s(x)=\pm\zeta$ lie in ${\cal
D}$, connecting the points $x=0$ and $x=2\pi$}\},
\]
while $\sigma_2=\inf \sup |\Im x|$ over the above pair of curves. In this case one should definitely have
\begin{equation}
\rho<{\pi\over2},
\label{best}
\end{equation}
for a bounded ${\cal D}$. One reason, for instance is that the
level curves $\Re s=\mp T$ for $T\geq T_{\psi}$ are contained
inside balls of radius $2r_{\psi}$ centered at $x=0$ and $x=2\pi$
respectively and are not homotopic in ${\cal
D}\setminus\{0,2\pi\}$. Furthermore, in terms of the maps $s,\,x$
(\ref{stime}) the $2\pi$-antiperiodicity of the function $\psi(x)$
combined with (\ref{pv}) result in:
\begin{equation}
s\circ l_{2\pi}=  l_{{ i \pi}}\circ\iota\circ s,\;\;x\circ\iota=
l_{2\pi}\circ x\circ l_{{i \pi}},
\label{2pi}
\end{equation}
where $l_{i\pi}:\,s\rightarrow s+i\pi$ (also defining a
diffeomorphism $\boldsymbol l_{i\pi}:\,(s,\varphi)\rightarrow
(s+i\pi,\varphi)$ and a canonical transformation $L_{i\pi}$ acting
on the momenta as the identity) and $l_{\pm i\pi}$ are identified
on $\C/2i\pi\Z$ as well as $l_{\pm2\pi}$ on $\T'$.

\begin{figure}[ht]
\includegraphics{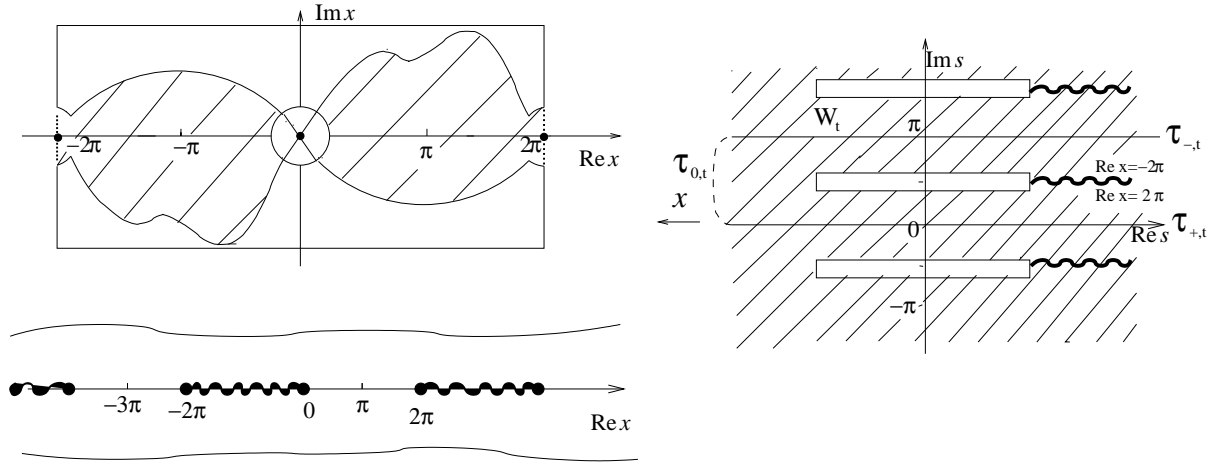}
\caption{The map $x$. Branch cuts on the right image are
unnecessary if one does not distinguish the values of $x$ modulo
$4\pi$. The cuts along the horizontals $\Im s=\pm{\pi\over2}$ are
the (coinciding) images of the vertical dotted diameters of the
small circles centered at $x=-2\pi,2\pi$. The cut on the right
image, where $\Im s={\pi\over2}$ [$\Im s=-{\pi\over2}$]
corresponds to the dotted radius positioned above [below] the real
axis on the upper left image. The lower left image illustrates how
various branches of $x(s)$ as a multi-valued function can be
constructed.}\label{fig3}
\end{figure}

Fig. 3 provides an illustration, and the antiperiodicity property
is expressed there by the fact that the area where
$\Re{x}\in[-2\pi,0)$ is congruent to the area where
$\Re{x}\in[0,2\pi)$ flipped about the real axis and translated
left by $2\pi$. Then as far as the definition (\ref{rho}) of the
quantity $\rho$, relevant to the pair $(\psi,{\cal D})$ is
concerned, by continuity there must exist a level curve
$\gamma^*_\zeta$ of $\Im s(x)$, emanating from $x=0$ with
$|\zeta|\leq\pi/2$ which will exit a bounded domain ${\cal D}$
before (ever) arriving to the point $x=2\pi$. This necessitates
the existence of  singular points of the functions $x(s)$ off the
real axis for $|\Re{s}|\leq T_{\psi},\,\rho\leq|\Im
{s}|\leq{\pi\over2}$, inside the unshaded rectangular regions in
Fig. 3, as it is the case with the classical pendulum, where the
singular points are $s=\pm i{\pi\over2}$. Dealing with the
classical pendulum where $\psi(x)=2\sin(x/2)$, the semi-width of
${\cal D}$ can be taken arbitrarily large which means
$\rho=\rho(\psi,{\cal D})$ approaches ${\pi\over2}$ from
below\footnote{As $\rho$ approaches ${\pi\over2}$ the functions
from the space $\bs$, etc. naturally grow unbounded. For optimal
splitting estimates one wants to have $\rho$ as large as possible,
which results in various technical nuances in the literature, see
\cite{Ga}, \cite{DGJS} etc.}. For the pendulum $\psi(x)$ is also
odd, so one should have $D_x\psi(\pi)=D_s\chi(0)=0$ and
$\chi(s)=2{\,\rm sech}\,{s},$ an even function of $s$.

However, further analysis of Hamiltonian (\ref{wwham},
\ref{sham}), see also (\ref{transs}), does require the ability to
distinguish the values of $x$ modulo $4\pi$ (unless
$\theta\in\Z^n$). In order to do so one can treat the function
$x(s)$ whose domain is shown in the upper left image in Fig. 3 as
a multi-valued, bi-real-analytic, $2{i \pi}$-periodic function
with values in $\C$. Its single branch labelled by $j\in \Z$ will
be defined by fixing $x(0)=4\pi j$. This is equivalent to taking
${\cal D}$ as a bi-infinite strip about the real axis, drawing
branch cuts at all the translations of $[-2\pi,0]$ by $4\pi j$ and
taking the lower limit of integration in (\ref{stime}) equal to
$\pi+4\pi j$, see the lower left image in Fig. 3. As a
multi-valued map, $x(s)$ has branch points inside the unshaded
rectangular regions in Fig. 3. There are branch cuts emanating
from a corner of each shaded region, whose exact appearance
depends on how exactly the covering of ${\cal D}$ (identification
of $x\in\C$ modulo $4\pi$) is defined. See the caption to Fig. 3.

\subsubsection*{Application to $H_\theta$}
In this case the pair $(\psi,{\cal D})$ is well defined, see the
end of Section 3.1. Recall that the Hamiltonian $H_\theta$
(\ref{wwham}) is viewed as a multi-valued function on
$T^*(\T'\times \T^n)$, where in view (\ref{transs}) a branch is
identified by a (real) value of $H_\theta$ at
$(y,I,x,\varphi)=(0,0,0,0)$. The application of the transformation
$\Xi_{\s}$ to a chosen branch of $H_\theta$ results in the ``new''
Hamiltonian $H_{\s}\circ \Xi_{\s} =H_{\s}(h,I,s,\varphi)\in\bsk$,
with $\fp=(r,\infty,\rho,\sigma)$ given by (\ref{sham}).
Application to any other branch of $H_\theta$ is tantamount to the
shift of the angles $\varphi$ by a multiple of $4\pi\theta$ and
does not require extra consideration. As the values of the
variable $x$ cannot be identified modulo $4\pi$ (\ref{transs})
unless $\theta\in\Z^n$, as far as the domain for the variable $s$
of $H_{\s}$, is concerned it should contain branch cuts, e.g. as
shown in the right image in Fig. 3. The splitting problem for
$H_{\s}$ can be briefly described as follows.

\medskip
\noindent In the absence of perturbation, the bi-infinite
bi-cylinder ${\cal C}_\infty$ is an invariant manifold $W_{\tt t}$
which is the unstable manifold to an invariant torus ${\cal
T}_{0,\tt t}$ at $s=-\infty$, see Figs 1, 3. As
$s\rightarrow\infty$ along the two different lines $\Im s=0$ and
$\Im s=\pi$, one arrives into a pair of different tori ${\cal
T}_{\beta,\tt t}$ ($\beta\in\{+,-\},$ the $+$ sign corresponding
to the former line; compare with Remark 3.7) corresponding to
$x=\pm 2\pi$ respectively in Fig. 1. $W_{\tt t}$ is a part of the
stable manifold for these tori, which can be analytically
continued further by choosing one of them and then flipping the
branch cuts in Fig. 3 with respect to the imaginary axis,
whereupon the chosen branch (above or below the branch cut along
the line $\Im s={\pi\over2}$ in Fig. 3) can be analytically
continued into the strip, whereof it was separated by a branch
cut. And so on.

The perturbed situation, with respect to Figs 1, 3 is
qualitatively  as follows. Let $H_0$ be the restriction of
$H_{\s}$ on $T^*{\cal C}_T$, for $T<\infty$. Theorem \ref{mt}
stipulates the existence of the perturbed manifold $W_0$, defined
as a graph over the semi-infinite bi-cylinder ${\cal C}_T$ in
terms of a generating function $S_0$. $W_0$ contains an invariant
torus ${\cal T}_0$ near (in the sense of Remark 3.2)
$\{s=-\infty\}$, for which it is an unstable manifold.
Furthermore, in order to apply Theorem \ref{mt} (twice) on the
bi-cylinder ${\cal C}_T^-$ going into $\{s=+\infty\},$ one should
use the analytical continuation of $H_{\s}(\cdot, s)$ in $s$ from
either the strip $|\Im s|\leq\rho$ or the strip $|\Im
s-\pi|\leq\rho$. Such an analytic continuation is roughly
tantamount to flipping the branch cuts in Fig. 3 with respect to
the imaginary axis. Denote these analytic continuations as
$H_{\beta}$ for $\beta=+,-$ respectively.  By (\ref{transs}) these
analytic continuations should in particular differ from one
another by the $4\pi\theta$-shift of the $\varphi$-variables and
can be both described using the sputnik Hamiltonian $H_0'=H_0\circ
L^{-2}_\chi$ to $H_0$, the set-up of Corollary \ref{sput} being
guaranteed by (\ref{transs}). An application of Theorem \ref{mt}
to each Hamiltonian of the pair $H_\beta$ results in a manifold
$W_{\beta}$ defined as a graph over the semi-infinite bi-cylinder
${\cal C}_T$ via a generating function $S_{\beta}$ and containing
an invariant torus ${\cal T}_{\beta}$ near $\{s=+\infty\}$ (i.e.
near $x=-2\pi,2\pi$ for $H_\theta$, with $\beta={\rm sign}\,x$)
where $W_\beta$ is the stable manifold.

The Hamiltonians $H_{0}$ and $H_{\beta}$ coincide for $s$ in the bounded strip around
$\R$ and $\R+i\pi$, respectively for $\beta=+,-$. This enables one to define
the splitting function on a finite bi-cylinder $\hat{\cal C}$ as $S_0-S_\beta$ for
$s$ in the corresponding strip. An
important issue that the cohomology classes of the one-forms
$dS_{0,\beta}$ are all equal to one another follows easily from
(\ref{transs}) and Corollary \ref{sput}.

\medskip
\noindent Technically, let us start out by calling
$H_\beta(\cdot,s)$ the $2i\pi$-periodic restriction of
$H_{\s}(\cdot,s)$ into the complex one-strips
$\Lambda_{\beta,\infty,\rho}$ described in Remark 3.7. Instead of
dealing with $H_{\beta}$ near $\{s=+\infty\}$, consider the
Hamiltonians $H_{\beta}\circ \fI_\beta$, where
\[
\fI_\beta:\,(h,I,s,\varphi)\rightarrow\left\{\begin{array}{ll}(-h,I,-s,\varphi), &\beta=+,\\
(-h,I,-s+2i\pi,\varphi), &\beta=-.\end{array}\right.
\]
Addition of $2i\pi$ to $-s$ in the second line is optional inside
the functional dependencies, as all the functions of $s$ involved
are $2i\pi$-periodic. It has been done above simply to make sure
that the restriction $\io_\beta$ of $\fI_\beta$ to the base space
maps the strip $|\Im s-\pi|<\rho$ into itself; it uses
$2i\pi$-periodicity of $H_{\s}(\cdot, s)$ in $s$. Then combining
(\ref{transs}) with (\ref{time}) one gets

\begin{equation}
H_\beta \circ \mathfrak{I}_\beta =
H_{-\beta}\circ{L}^{-2}_\chi\circ{L}^\beta_{2\pi\theta}\circ
L^{-\beta}_{i\pi}\equiv H_{-\beta}\circ{L}_{*,\beta},
\label{transss}
\end{equation}
the transformation $L_\chi$ having been defined earlier by
(\ref{lab}). Namely (the correction below $[+2i\pi]$ not appearing
for $\beta=+$)
\begin{equation} H_\beta(-h,I,-s[+2i\pi],\varphi)=
H_{-\beta}(h-2\lambda\chi^2(s-\beta i\pi),I,s-\beta i\pi,\varphi+\beta 2\pi\theta).
\label{ss}
\end{equation}
Indeed, the last two formulas follow from (\ref{2pi}) regarding
the presence of branch cuts for the map $x(s)$ by simply matching
$\beta$ with ${\rm sign}\,x$, for real $x$. Namely, in the second
formula in (\ref{2pi}) the shifts $l_{\pm2\pi}$ and $l_{\pm i\pi}$
have been identified. To make a choice of the sign for the formula
(\ref{transss}) all one has to do is to act by $l_{-2\pi}$ on
$x\in(0,2\pi)$ and  by $l_{2\pi}$ on $x\in(-2\pi,0)$; in the same
fashion $l_{-i\pi}$ acts on $s:\,\Im s=i,$ corresponding to
$x\in(-2\pi,0)$, and  $l_{i\pi}$ acts on real $s$ corresponding to
$x\in(0,2\pi)$.

The two functions $H_\beta(\cdot,s)\circ \fI_\beta$ allow analytic
continuation in $s$ into the one-strip $\Lambda_{\beta,\te,\rho}$,
see Remark 3.7, where they are characterized by the same array of
non-degeneracy and smallness parameters (in particular $\nu=1$).
The pair of transformations ${L}_{*,\beta}$ defined in
(\ref{transss}) plays the role of sputnik transformations, as they
obviously satisfy Assumption \ref{sploc}. On the other hand, the
Hamiltonians $H_\beta(\cdot,s)$ represent the same analytic
function $H_0(\cdot,s)$. The rest of the development is clear:
Theorem \ref{mt} and Corollary \ref{sput} are satisfied by $H_0$
and a pair of its sputniks $H_\beta\circ\fI_\beta$. We proceed
with some generalization.

\subsection{Splitting theory on ${\cal C}_\infty$}\label{st}
Let $\beta\in\Z_2\equiv\{+,-\}$, $\kappa\gg1$,
$\fp=(r,\infty,\rho,\sigma)$. Consider a Hamiltonian $H\equiv
H_0\in\bsk$ in the cotangent bundle $T^*{\cal C}_\infty$ of the
bi-infinite bi-cylinder ${\cal C}_\infty$, in the form (\ref{pt}).
Let $H_\beta(\cdot,s)$ be the restrictions of the function
$H_0(\cdot,s)$ into bi-infinite strips $|\Im s|\leq \rho$ and
$|\Im s-\pi|\leq\rho$ for $\beta=+,-$ respectively. Suppose, the
Hamiltonians $H_\beta\circ\fI_\beta=-\lambda
h+\spr{\omega}{I}+\ldots$ allow analytic continuation into
one-strips $s\in\Lambda_{\beta,\infty,\rho}$ introduced by Remark
3.7, namely into the region $\{\Re s\leq -2T_{\psi}\}$,
$2i\pi$-periodically. For the above analytic continuations let us
still use the notations $H_\beta\circ\fI_\beta$, the latter
quantities being well defined in the cotangent bundle over the
one-cylinders ${\cal C}_{\beta,\fp}$. Suppose, the restrictions of
the quantities $H_0(\cdot,s)$ and $H_\beta\circ\fI_\beta(\cdot,s)$
over $\Re s\leq T<\infty$ satisfy the conditions of Theorem
\ref{mt}, with the same analyticity parameters
$\kappa,r,T,\rho,\sigma$ and equivalent (footnote \ref{feq})
non-degeneracy and smallness parameters $\lambda,R,M,\mu,\nu$.

In particular this means that for $\mu=0$ the by-infinite
bi-cylinder  ${\cal C}_\infty$ would be an invariant Lagrangian
manifold for $H_0$, asymptotic to a torus ${\cal T}_0$ at
$s=-\infty$ and a pair of tori ${\cal T}_\beta$ at $s=+\infty$,
see Fig. 3. In addition, assume the following.

\begin{assumption}[Sputniks on ${\cal C}_\infty$] There exists a pair of canonical transformations
\begin{equation} L_{*,\beta}=L_{*,\beta}(\ar^*_\beta,S^*_\beta):\;\left\{
\begin{array}{llllllll}
\q&=&\ar^*_\beta(\q_*),\\ \p&=&\ti{(d\ar^*_\beta
)}\p_*+dS^*_\beta(\q),
\end{array}
\right. \label{sptn}
\end{equation}
where $\ar^*_\beta$ are diffeomorphisms of ${\cal C}_\infty$, such that
 $\ar^*_\beta-{\tt id}\in[\bs]^{n+1}$, as well as
$S^*_\beta\in\bstwo$, such that
\[
H_\beta\circ\fI_\beta=H_{-\beta}\circ L_{*,\beta}.
\]
\label{SP}
\end{assumption}

\noindent Then one has the following lemma.
\begin{lemma}
Let $\mathfrak{i}\in\{0,\beta\}$. There exist invariant Lagrangian manifolds
\begin{equation}
\begin{array}{lll}
W_0&=&\overline{\{(h,I;s,\varphi)\in\R^{n+1}\times {\cal
C}:\;\,\;h=D_s S_0(s,\varphi),\,I=D_\varphi S_0(s,\varphi)\}},
\\
W_\beta&=&\overline{\{(h,I;s,\varphi)\in\R^{n+1}\times {\cal
C}^-_\beta: \;\,h=D_s S_\beta(s,\varphi),\,I=D_\varphi
S_\beta(s,\varphi)\}},
\end{array}
\label{lm}
\end{equation} contained in the level set $H^{-1}(c_0)$ for some $c_0$ satisfying (\ref{c0}).  One has $S_0\in\bsprime$, $S_\beta\circ\io\in\bsprimeb$, they satisfy the corresponding bounds of (\ref{transf}) and $[dS_{\mathfrak
i}]=\xi\in\R^n,\,\forall\mathfrak i$. The manifolds $W_{\mathfrak i}$ contain invariant tori ${\cal
T}_{\mathfrak{i}}$ (whereupon the flow is conjugate to a rotation with the frequency $\omega$) near $s=-\infty$
for $\mathfrak{i}=0$ and $s=+\infty$ for $\mathfrak{i}=\beta$. \label{things}
\end{lemma}
{\tt Proof:} The lemma is an immediate consequence of Theorem
\ref{mt} and Corollary \ref{sput}; see also Remark 3.7 following
the latter. Indeed, Assumption \ref{SP} implies that Assumption
\ref{sploc} is satisfied in the sense that $H_\beta\circ\fI_\beta$
is a sputnik of $H_0$ (restricted as a function of $s$ from the
bi-strip $\Lambda_{\te,\rho}$ to a one-strip
$\Lambda_{-\beta,\te,\rho}$ to yield $H_{-\beta}$) under the
sputnik transformation $L_{*,\beta}.\,\Box$

\medskip
\noindent The above lemma is central for the splitting problem,
for now one can introduce the {\em splitting distance} as an exact
one-form on the bounded bi-cylinder ${\displaystyle \hat{\cal
C}_{T'}=\bigcup_{\beta\in\{+,-\}} \hat{\cal C}_{\beta,T'}}$ by
defining it separately on each of the above components
(corresponding to $\beta={\rm sign}\,x$ as far as the original
simple resonance splitting problem is concerned):
\begin{equation}
d{\fS}(s,\varphi) = d[S_0(s,\varphi)-S_{\beta}(s,\varphi)],\;\;(s,\varphi)\in
\hat\Lambda_{\beta,\te',\rho'}\times\T^n_{\sigma'}. \label{spd}
\end{equation}
Let us call the function ${\fS}\in\bshatprime$ the {\em splitting
potential}. Note that the (complexified) domains $\hat{\cal
C}_{\beta,\te',\rho',\sigma'}$ are disjoint for different $\beta$.
In  terms of the notations introduced in (\ref{lbd}, \ref{pr}) and
Remark 3.7, one has
$\hat\Lambda_{+,\te,\rho}=\Pi_{\te,\rho},\,\hat\Lambda_{-,\te,\rho}=\Pi_{\te,\rho}+i$,
where $\Pi_{\te,\rho}$ is simply a symmetric rectangle in $\C$
with the half-length $T$ and half-width $\rho<{\pi\over2}$.

To make things easier, let us simply view ${\fS}(s,\varphi)$ as a
double-valued function, bounded and real-analytic for
$(s,\varphi)\in \Pi_{\te',\rho'}\times \T^n_{\sigma'}$ by changing
$s\rightarrow s+i\pi$ for $s\in\hat\Lambda_{-,T',\rho'}$ (the
strip about the line $\Im s=\pi$). Hence the rest of the
statements and estimates will be valid for either one of the two
branches of $\fS$ over the domain $\Pi\times\T^n$, the index
$\beta$ being mostly omitted. The functions one is dealing with
are members of the space
$\fB_{\te',\rho',\sigma'}(\Pi\times\T^n)$, whose element $u$ can
be represented as a uniformly convergent Fourier series in
$\varphi$ with coefficients $u_k(s),\,k\in\Z^n$, bounded and
holomorphic for $s\in\Pi_{\te',\rho'}$, with
$u_{-k}(s)=u^*_{k}(s)=u_{k}(s^*)$ (${}^*$ marking the complex
conjugate). Thus in the sequel parameter vectors $\fp$ will have
three components $(T,\rho,\sigma)$, and as usual $\fp'<\fp$ and
$\de=|\fp-\fp'|$. Besides the use of a certain finite number of
intermediate values $\fp''$ of the analyticity parameters such
that $\fp'<\fp''<\fp$ is implied by default along the way.

\begin{lemma}
The splitting potential ${\fS}$ satisfies a homogeneous quasi-linear PDE in $\Pi_{\te',\rho'}\times
\T^n_{\sigma'}$:
\begin{equation}
([\lambda+\lambda_\mu(s,\varphi)]D_s+
\spr{\omega+\omega_\mu(s,\varphi)}{D_\varphi}){\fS}=0,
\label{qlpde}
\end{equation}
where the pair $\g=(\lambda_\mu,\omega_\mu)\in [\fB_{\te',\rho',\sigma'}(\Pi\times\T^n)]^{n+1}$
 satisfies the bounds
\begin{equation}
|(\lambda_\mu,\omega_\mu)|_{\fp'}\leq C \mu\eta^{-1},
\label{moreb}
\end{equation}
where $\eta$ is defined by (\ref{eeta}) and $C$ is of the same order as in Theorem \ref{mt}.
\label{hom}
\end{lemma}
{\tt Proof:} Follows by Corollary \ref{hjc}: each single $S$ is a solution of the Hamilton-Jacobi equation
for $H$ on the energy level $c_0$, thus (\ref{qlpde}) is obtained by subtracting the equation for the former
function from the same equation for the latter one. The estimate (\ref{moreb}) follows from Definition
\ref{nfdef} and the bound for $\g$ in (\ref{aspt}) as well as the bound (\ref{transf}) for the norm of $dS$
(which also turns out to be the estimate for $M|dS|$, see (\ref{dshat}-\ref{c1}) for detail. As the result one
may have to multiply the constant $C$ in Theorem \ref{mp} by a factor, depending on $n$ and $\tau$ only. $\Box$

\medskip
\noindent A prototype of the following lemma is due to Eliasson
\cite{El}.
\begin{lemma}
A branch of the function $\fS$ for
$(s,\varphi)\in\Pi\times\T^n$ has at least $n+1$ critical points $\varphi_c=\varphi_c(s)$,
given $s$, i.e. where $D_\varphi \fS(s,\varphi_c)=0$.
\label{orb}
\end{lemma}
{\tt Proof:} This statement is a consequence of the fact that
given $s\in\Pi$, the function $\fS(s,\varphi)$ is $2\pi$-periodic
in each component of $\varphi$, by Lemma \ref{things}, essentially
stating that the one-form $d\fS$ is exact. The number $n+1$ of
critical points is the Ljusternik-Schnirelmann characteristic of
the torus $\T^n$. $\Box$

\medskip
\noindent The statement of the next lemma is similar to Corollary
\ref{cjt}, claiming the structural stability of the constant
vector field $\x_0=\lambda{\partial\over\partial
s}+\spr{\omega}{{\partial\over\partial \varphi}}$ on the bounded
one-cylinder $\Pi\times\T^n$ under small perturbations. It is
crucial for the exponentially small estimate (\ref{sir}). The
prototype of this result was proved by Sauzin \cite{Sa} regarding
the so-called characteristic vector field.
\begin{lemma}
There exists a constant $C=C(n,\fp)$ but independent of
$\omega$, such that for $\x=\x_0+\g$, with $\g\in \fB_{\fp}(\Pi\times\T^n)$ such
that
\begin{equation}
|\g|_{\mathfrak{p}}\;\leq\; C^{-2}\lambda\de,
\label{gsmall}
\end{equation}
there exists a diffeomorphism $\ar={\tt id}+{\br}$ with
 $\br\in[\bspiprime]^{n+1}$,
effecting the conjugation
${\displaystyle d\ar^{-1}\circ\x\circ\ar=\x_0,}$
with
\begin{equation}
\begin{array}{rcl}
|\br|_{1,\fp'}& \leq &C|\g|_{\mathfrak{p}}(\lambda\de)^{-1}.
\end{array}
\label{mucg}
\end{equation}
\label{cj}
\end{lemma}
{\tt Proof:} This lemma is yet another implicit function theorem regarding the
operator $\Dlw$ introduced in (\ref{dhw}). One can rewrite the conjugacy
problem in question as
\[
\g\circ({\tt id}+\br)-\Dlw \br=0.
\]
Or in Hamiltonian terms, one seeks a canonical
transformation $\hat\Psi:\,\q=\ar(\q'),\,\p=\ti{(d\ar)}\p',$ such that a linear
Hamiltonian $\hat{H}(\p,\q)=\lambda h +\spr{\omega}{I}+\spr{\g(\q)}{\p}$ is
conjugate to $\lambda h +\spr{\omega}{I}$.

The vector field conjugacy problem is certainly amenable to an abstract
implicit function theorem \cite{Z1}, \cite{RW3} however the latter would not
provide the optimal condition (\ref{gsmall}), as well as for the conjugacy problem (\ref{fldcj}).
In order to get (\ref{gsmall}) one should follow the standard iterative scheme mimicking the proof of
Corollary \ref{cjt} (which in turn is a particular case of the proof of Theorem
\ref{mt}) basing it however on  Proposition
\ref{mppp} rather than Proposition \ref{mp}.

The latter proposition analyzes the possibility of finding a
solution $u$ to a PDE $\Dlw u=v$ on $\Pi_{\te,\rho}\times
\T^n_\sigma$, such that $|u|_{\fp'}$ can be bounded irrespective
of $\omega$. The kernel of the operator $\Dlw$ on $\bspi$ consists
of all functions, which are represented by Fourier series in the
variable $\phi=\varphi-{\omega\over\lambda} s$, and the norm of
such a function in $\Pi_{\te,\rho}\times \T^n_\sigma$ clearly does
depend on $\omega$. To avoid it one is naturally led to solving a
Cauchy problem set up by conditions (\ref{prt}, \ref{cf}) in
Appendix A, as a way to determine the required inverse of the
operator $\Dlw$.

A single application of Proposition \ref{mppp} furnishes the approximate (first
order) solution $\br_1$, satisfying \ref{mucg}, whereupon the
perturbation $\g$ changes to $\g_1$, the norm of $\g_1$ bounded by
$|\g|_{\fp}|\br_1|_{1,\fp'}$. This fact constitutes the analogue of the Iterative
lemma \ref{il}, whereupon the standard dyadic iterative procedure (see e.g. Appendix B) is
run. $\Box$

\medskip
\noindent Results of the type of Lemmas \ref{hom}-\ref{cj}, as far
as exponentially small splitting is concerned have been a target
of a number of works of Lazutkin starting from \cite{La} and
followers, see e.g. \cite{DGJS}. Indeed from Lemmas
\ref{hom},\ref{cj} one can further easily deduce the following
upper bound for the infinity norm (the supremum over the real
values of the variables only) for each branch of the splitting
function $\fS$.
\begin{theorem}
Suppose, the assumptions of Section \ref{st} are satisfied and the smallness condition
(\ref{munot}) holds, with a large enough $C=C(n,\tau,\psi,\fp)$.
Then
\begin{equation}
|\fS|_\infty\;\;\leq\;\; C \mu\eta^{-1}\sum_{k\in\Z^n\setminus\{0\}}
\exp\left(-|\spr{k}{{\omega\over\lambda}}|\rho' - |k|\sigma'\right),
\label{expsmall}
\end{equation}
where $\eta$ is defined by (\ref{eeta}), $\rho'=\rho-\de<\pi/2$,
$\sigma'=\sigma-\delta$.
\label{estimate}
\end{theorem}
{\tt Proof:} The estimate is clearly the same for each branch of
the double-valued function $\fS$ on $\Pi\times\T^n$. Consider one
branch. Let $\ar$ be the conjugating diffeomorphism of Lemma
\ref{cj} (the estimates of the lemma are uniform in $\beta$). Then
the function $\fS'={\fS}\circ\ar$ is constant along the flow lines
of the constant vector field $\x_0$ on the bounded one-cylinder
$\Pi\times\T^n$. Then it is a real-analytic function on $\T^n$:
one can formally write ${\fS}'={\fS}'(\varphi-{\omega\over\lambda}
s)={\fS}'(\phi)$, where $\phi=\varphi-{\omega\over\lambda} s$ and
expand it into the Fourier series
\begin{equation}
{\fS}'(s,\varphi)=\sum_{k\in\Z^n\setminus\{0\}} {\fS}'_{k} e^{-i
\spr{k}{{\omega\over\lambda}}s}e^{i\spr{k}{\varphi}}.
\label{pp}
\end{equation}
Then ${\fS}'(s,\varphi)$ is a quasi-periodic function of $s\in \Pi_{\te',\rho'}$. Besides the standard
complex analysis technique for estimating the Fourier coefficients yields (once
again scaling the analyticity loss parameters by, say factor 4)
\begin{equation}
\left| {\fS}'_{k}e^{-i
\spr{k}{{\omega\over\lambda}}s} \right|\,\leq\, e^{-|k|\sigma'}|{\fS}|_{\fp'},\;\;\forall
\,s\in\Pi_{\te',\rho'}.
\label{ppone}
\end{equation}
This implies, as $|\Im s|\leq\rho'$ that
\[
| {\fS}'_{k} | \;\;\leq\;\;|{\fS}|_{\fp'}
\exp\left(-|\spr{k}{{\omega\over\lambda}}|\rho' - |k|\sigma'\right),
\]
and consequently (\ref{expsmall}) as one can use (\ref{transf}) for
$|{\fS}|_{\fp'}$, while
${\displaystyle
|{\fS}|_\infty \asymp|{\fS}'|_\infty\leq\sum_{k\in\Z^n\setminus\{0\}}|{\fS}'_{k}|}.$
$\Box$

\subsubsection*{Conclusion of the proof of Theorem \ref{main}}
The splitting problem for the Hamiltonian $H_{\s}$ given by
(\ref{sham}) satisfies the conditions of Theorems \ref{mt} and
\ref{estimate}, with $\omega=\omega_1={\omega_0\over \sqrt{\eps
R_0}},$ $R=1$ and $M=M_0/R_0$.

Assuming (\ref{sc}) and letting $\mu=\seps\lambda^2
(M_0/\sqrt{R_0})\inf(\varsigma_0\delta_0,R_0)$, $\nu=1$, it is
easy to check that for the application of Theorem \ref{mt} (with
$\kappa>1$, well-defined quantities $r_\psi,T_\psi$, as well as
$\rho$ defined by (\ref{rho}), $\sigma=\sigma_1+{1\over2}\delta_0$
and $\de\lesssim\delta_0$) one has $\eta^{-1}\lesssim
\sup(M_0R_0^{-1} \sqrt{\eps R_0}
(\lambda\varsigma_0\delta_0)^{-1},1)\lesssim\lambda^{-1}.$ Thus
(\ref{sc}) with the proper choice of the bounding constant,
depending on $(n,\tau,\psi,{\cal D},\lambda,\sigma)$ ensures the
applicability of Theorem \ref{mt}.

Besides, the sputnik property, in order to satisfy Assumption
\ref{SP} has come a long way: (\ref{trans}, \ref{transs},
\ref{transss}). To establish the fact that all the three tori
${\cal T}_{0,\beta}$ claimed by Lemma \ref{things} correspond to
the same torus ${\cal T}$ (see Fig. 1) for the Hamiltonian
(\ref{nfh}) one should chase back through relations
(\ref{transss}, \ref{transs}, \ref{trans}) and notice that the
manifolds $W_\beta$ arise from the sputnik manifold $W'_0$ to
$W_0$, described by Corollary \ref{sput}, with the sputnik
transformation $L_*=L^{-2}_\chi$, alias $L^{-2}_\psi$, solely via
a translation $x\rightarrow x-\beta 2\pi$ in terms of the
variables $(x,\varphi)$ of Hamiltonian (\ref{wham}).

Finally, we show how the transformation $\Xi_\theta$ (\ref{two}) with the underlying base space transformation
$\ar_\theta:\,(x,\varphi)\rightarrow(x,\varphi+\theta x)$ affects the estimate (\ref{expsmall}). Suppose
$\fS_{\s}$ is the splitting potential, defined according to (\ref{spd}) for the Hamiltonian $H=H_{\s}$. Theorem
\ref{mt} followed by Lemmas \ref{hom}-\ref{cj} imply that the Hamilton-Jacobi equation for the quantity
$\fS_\psi=\fS_{\s}\circ \s^{-1}\circ\ar_\theta^{-1}$ in the variables $(x,\varphi)$ of Hamiltonian (\ref{wham})
is conjugate to
\[
(\lambda\psi(x) D_x + \spr{\omega_1+\lambda\theta\psi(x)}{D_\varphi}){\fS}'_\psi(x,\varphi)=0,
\]
by a near-identity change $\ar$ of variables $(x,\varphi)$ (where ${\fS}'_\psi=\fS_\psi\circ \ar$)
with the total analyticity loss of the order of $\delta_0$.
Therefore instead of (\ref{pp}) one has
\[
{\fS}' (s,\varphi)=\sum_{k\in\Z^n\setminus\{0\}} {\fS}'_{k} e^{-i
\spr{k}{{\omega_1\over\lambda}s+\theta x(s)}}e^{i\spr{k}{\varphi}}.
\]
where ${\fS}'={\fS}'_\psi\circ \s$. Then with $\rho'=\rho-\de$ and $\sigma_2'=\sigma_2-\de$, by mimicking
(\ref{ppone}) one gets the estimate
\[
| {\fS}'_{k} | \;\leq\;|{\fS}'_\psi|_{\fp'}
\exp\left(-|\spr{k}{{\rho'\over\lambda}\omega_1+\sigma'_2\theta}| - |k|\sigma_1\right),
\]
which completes the proof of Theorem \ref{main}, upon removing
primes for the pair $(\rho,\sigma_2)$ in the formulation of the
theorem. Note that the lower bound $2n+2$ for the number of
homiclines comes from the application of Lemma \ref{orb} to the
two (upper, lower) separatrix branches, see Fig. 1.  $\Box$

\section{Appendices}
\renewcommand{\thesection}{}
\renewcommand{\thesubsection}{\Alph{subsection}}
\renewcommand{\theequation}{\Alph{subsection}.\arabic{equation}}

\subsection{First order linear PDEs on bi-cylinders}
\setcounter{equation}{0}
The following set of proposition addresses the issue of the existence of the
right inverse for the operator $\Dlw$ defined by  (\ref{dhw}), with $\omega\in\di$
and $\lambda>0$.
\begin{proposition}
Let $\mathfrak{p}=(r,T,\rho,\sigma)$ and $v\in\bs$. There
exists a real $c,\,|c|\leq |v|_\fp$, such that the solution of
the equation ${\displaystyle \Dlw u=v-c}$ exists in $\bsprime$ for
$\mathfrak{p}'=(r,T,\rho,\sigma')$ with
$0<\sigma'<\sigma$. Let $\sigma-\sigma'=\delta$ and  $\varsigma = \inf(
\gamma\delta^{\tau},\,\lambda)$. Then
\[
|u|_{\fp'}\,\lesssim \, \varsigma^{-1}|v|_{\fp}.
\]
\label{mp}
\end{proposition}

\noindent {\tt Proof:} Following (\ref{dcmp}), let $v=v_0+v_1$, where
$v_0\in{\fB}_\sigma(\T^n)$ and $v_1\in\bsone$. Seek
$u(s,\varphi)=u_0(\varphi)+u_1(s,\varphi)$, such that $D_\omega u_0=v_0-c$ and
$\Dlw u_1=v_1$. The first equation obeys Proposition \ref{rus}, if one
chooses $c=\langle v_0\rangle$. The second
one is solved in ${\cal C}_{\mathfrak{p}}$ by the method of characteristics,
regarding the fact that one can write $v_1(s,\varphi)=\chi(s) w(s,\varphi)$, for some $w\in \bs$:
\begin{equation}
\begin{array}{lll}
u_1(s,\varphi)\,=\,\int_{-\infty}^0 v_1(s+\lambda t,\varphi+\omega t)dt&=&
\lambda^{-1}\int_{-\infty}^0 \chi(s+t)w(s+t,\varphi+\omega\lambda^{-1}
t)dt \\ \hfill \\
&=& \lambda^{-1}\int_{-\infty+{i}\Im{s}}^s \chi(\zeta)w[\zeta,\varphi-\omega\lambda^{-1}
(s-\zeta)]d\zeta \\ \hfill \\
&=&\lambda^{-1}\int_0^{x(s)}
w[ s(x),\varphi-\omega\lambda^{-1}(s- s(x))]dx
\end{array}
\label{char}
\end{equation}
The latter integral was obtained via the substitution $\zeta= s(x)$, see
(\ref{stime}) the integration fulfilled along the level curve $\Im
 s(x)=\Im s$, see Fig. 2. The latter integral is
 bounded by a $(\psi,{\cal D})$-depending constant times the norm of $v_1$, see (\ref{dcp}). Thus $u_1\in
\bsone$, with the norm
$|u_1|_{\mathfrak{p}}\,\lesssim\,\lambda^{-1} |v|_{\mathfrak{p}}$.
Bi-real analyticity of $u$ follows by construction. $\Box$

\begin{proposition}
Given  $v\in{\fB}_\sigma(\T^n)$  the
solution of the equation ${\displaystyle (-\lambda+D_\omega) u=v}$ exists in
${\fB}_\sigma(\T^n)$ and
\[
|u|_{\sigma}\lesssim \lambda^{-1}|v|_{\sigma}.
\]
\label{mp1}
\end{proposition}

\noindent {\tt Proof:} Clearly $u(\varphi)$ can be found explicitly as a
Fourier series in $\varphi$, whose coefficients $u_k$ are expressed via the
Fourier coefficients $v_k$ of $v(\varphi)$ as follows: $\;{\displaystyle
u_k={v_k\over-\lambda+ i\langle k,\omega\rangle},\;k\in\Z^n.}$ If $v_{-k}={v}^*_k$ (complex conjugate),
this property i.e. real analyticity is clearly retained by $u$.
 $\Box$

\begin{proposition}
Let $\mathfrak{p}=(r,T,\rho,\sigma)$ and  $v\in\bsm$. There exists
a real constant $c,\,|c|\lesssim |v|_{\fp}$, such that the
solution of the equation ${\displaystyle \Dlw u=v-c}$ exists in
$\bsmprime$ for $\mathfrak{p}'=(r,T,\rho,\sigma')$ with
$0<\sigma'<\sigma$. With the same $\delta$ and  $\varsigma$ as in
Proposition \ref{mp}, one has
\[
|u|_{\fp'}\,\lesssim \,\varsigma^{-1}|v|_{\fp}.
\]
\label{mp2}
\end{proposition}
\noindent {\tt Proof:} By Definition of the space $\bsm$, $v$ admits a
unique decomposition $v(s,\varphi)= {v_0(\varphi)\over
\chi(s)}+v_1(s,\varphi)$, where $v_0\in\bst$ and
$v_1\in \bs$. Seek the solution ${\displaystyle
u(s,\varphi)={u_0(\varphi)\over \chi(s)}+u_1(s,\varphi)}$. Then
\[
{1\over\chi} \left( -\lambda{d\chi\over\chi}+D_\omega\right)u_0+\Dlw u_1 =
{v_0\over\chi}+v_1 -c.
\]
As ${\displaystyle {d\chi(s)\over\chi(s)}= 1 + \chi(s)\eta_1(s)}$, see
(\ref{calc}, \ref{calc1}) where the function $\eta_1(s)\in\bs$, $u_0$ can
be taken as  the solution of the equation ${\displaystyle
(-\lambda+D_\omega)u_0=v_0}$, which exists by Proposition \ref{mp1}, while
$u_1$ should satisfy
\[
\Dlw u_1 = v_1+\lambda \eta_1 u_0-c.
\]
The first two terms in the right hand side are members of $\bs$, so
Proposition \ref{mp} does the job, with the constant $
c=\langle v_1\rangle +
\langle v_0\rangle\psi_{xx}(0)$. $\Box$

\medskip
\noindent Combining Propositions \ref{mp} and \ref{mp2}, one gets
\begin{proposition}
Let  $\boldsymbol v\in\bsg$. There exists  a constant $\boldsymbol
c\in\R^{n+1},\,|\boldsymbol c|\lesssim |\boldsymbol v|_{\fp}$,
such that the solution of the $(n+1)$-vector equation
${\displaystyle \Dlw \boldsymbol u=\boldsymbol v-\boldsymbol c}$
exists in $\bsgprime$ for $\mathfrak{p}'=(r,T,\rho,\sigma')$ with
$0<\sigma'<\sigma$. With the same $\delta$ and  $\varsigma$ as in
Proposition \ref{mp}, one has
\[
|\boldsymbol u|_{\fp'}\,\lesssim\,\varsigma^{-1}|\boldsymbol
v|_{\fp}.
\]
Let $\boldsymbol v=(v,V)$. If for $j=1,\ldots,n,\,$ $\langle
V_j\rangle=0$, then $\boldsymbol c=(c,0)$, where $c$ is the same as in
Proposition \ref{mp2}. \label{mpp}
\end{proposition}

The proposition herein pertains to a conjugacy  problem on the
bounded one-cylinder $\Pi\times\T^n$. Let $\fp=(T,\rho,\sigma)$. A
function $v\in\bspi$ is given as a Fourier series\[ v(s,\varphi)=
\sum_{k\in\Z^n} v_k(s) e^{i\spr{k}{\varphi}},
\]
where $v_{-k}(s)=v^*_k(s)=v_k(s^*)$; ${}^*$ marks the complex
conjugate. Given $\omega\in\R^n$ (not necessarily Diophantine) let
\begin{equation}
\Z^n_\omega\equiv\{k\in\Z^n:\spr{k}{\omega}>0\}.
\label{prt}
\end{equation}
For $k\in \Z^n_\omega$ and a fixed $\rho'<\rho$ (let also
$\fp'=(\T',\rho',\sigma')<\fp, \de=|\fp-\fp'|$) denote
\begin{equation}
u_k[v]=-{i\over\lambda}\int_{0}^{\rho'}  v_k(i\zeta)
e^{-{\spr{k}{\omega}\over\lambda}\zeta} d\zeta. \label{cf}
\end{equation}
For $k\in -\Z^n_\omega,$ let $u_k[v]= u^*_{-k}[v]$, for $k$
such that $\spr{k}{\omega}=0,$ let $u_k[v]=0$.

\begin{proposition}
Let $\lambda>0$. A Cauchy problem
\[
\Dlw u =
v,\;\;\;u(0,\varphi)=  \sum_{k\in\Z^n}u_k[v]
e^{i\spr{k}{\varphi}}
\]
has a unique solution $u\in \bspiprime$, with $|u|_{\fp'}\leq C\lambda^{-1}
|v|_{\fp}$,  $|du|_{\fp'}\leq C(\lambda\de)^{-1}
|v|_{\fp}$ where $C$ may depend on $n$ and $\fp$ but is independent of $\omega$.
\label{mppp}
\end{proposition}
{\tt Proof:} Seek ${\displaystyle u(s,\varphi)=\sum_{k\in\Z^n}
u_k(s) e^{i\spr{k}{\varphi}}}$. Then $u_k(s)$ satisfy
\[
\lambda u'_k+i\spr{k}{\omega}u_k=v_k,\;\;u_k(0)=u_k[v].
\]
For
$k$ such that $\spr{k}{\omega}=0$ let
\[
u_k(s)={1\over\lambda}\int_0^s v_k(t)dt.
\]
For $k\in\Z^n_\omega$ let
\[
u_k(s)={1\over\lambda}\int_{i\rho'}^s v_k(t)\,e^{i{\spr{k}{\omega}\over\lambda}(t-s)}dt,
\]
where the integral can be taken along the part of the imaginary
axis until $\Im t=\Im s$ and then along the horizontal line. For
$k\in-\Z^n_\omega$ take the lower limit of integration as
$-{i\rho'}$. Then the integrand is always bounded by
$\sup_{s\in\Pi_{T,\rho}}|v_k(s)|$ in the absolute value. The
initial conditions are satisfied: by real analyticity of $v$, one
has $v_{-k}(-i\zeta)=v_k(i\zeta)$ for $\zeta\in\R$. In particular,
the solution $u$ is also real analytic. Note that for the elements
$u(s,\varphi)$ of $\bspi$, represented by Fourier series in
$\varphi$ with analytic coefficients $u_k(s)$ the supremum norm in
$\Pi_{\te,\rho}\times\T^n_\sigma$ is equivalent to the norm
defined as ${\displaystyle
\sum_{k\in\Z^n}\sup_{\Pi_{\te,\rho}}|u_k(s)|}$, the comparison
constants depending in particular on $\sigma$. $\Box$

\subsection{Conclusion of the proof of Theorem \ref{mt}}
\setcounter{equation}{0} First note that the smallness condition
(\ref{munot}) has simply combined the smallness condition
(\ref{small}) on $\mu$ for the Iterative lemma to be valid with
the lemma's remainder estimate (\ref{newmu}) on $\mu'$, simply to
ensure that $\mu'<\mu$. The rest to ensure that actually $\mu'$ is
many enough times smaller than $\mu$, so that the Iterative lemma
can be applied again and again, with a smaller and smaller
analyticity loss. This is achieved simply by choosing the constant
$C$ in (\ref{munot}) small enough. Without loss of generality (for
this part of the proof) assume that
$\kappa-\kappa'=r-r'=\rho-\rho'=T-T'=\sigma-\sigma'=\de=\delta$,
as well as $\delta<\lambda$ and  $\nu>M^{-1}\varsigma\delta$.

Take a geometric sequence $\{\delta_j=2^{-j}\delta\}_{j\geq1}$. Let
$\sigma_j=\sigma-\sum_{l=1}^j\delta_l$. Define sequences $\{\kappa,r,T,\rho\}_{j\geq1}$ in the  same way. Denote
$\mathfrak{p}_j=(r_j,T_j,\rho_j,\sigma_j)$. Identify the parameters
$\kappa,r,T,\rho,\sigma,\lambda,R,M,\nu,\mu$ with themselves, endowed with
zero indices. Let $C_0>1$ be the constant, whose existence is stated by Lemma \ref{il}.

If $\mu$ satisfies the smallness condition (\ref{munot}) with some $C\geq
2^{2\tau+3}C_0$, the assumption (\ref{small}) of Lemma \ref{il} is satisfied
for a single application of the lemma, with an analyticity loss
$\delta_1=2^{-1}\delta$ and a parameter $\varsigma_1=2^{-\tau}\varsigma$
(playing the role of $\varsigma$ in (\ref{eeta}) implying that $\eta_1\geq
2^{-(\tau+1)}\eta$. This results in a coordinate change
$\Xi_1=\Xi_1(\ar_1,S_1)$. At the output, according to (\ref{newmu}) one will
have the perturbation  parameters ${\displaystyle \nu_1\geq
2^{-(\tau+1)}M^{-1}\delta^{\tau+1}}$ and ${\displaystyle \mu_1\leq
2^{2\tau+2}C_0^2C^{-2}\mu\leq 2^{-2\tau-4}\mu,}$ as well as the new quantities
$\lambda_1,\,R_1,\,M_1$, such that
\[
m_1\,=\,\sup\left(\lambda^{-1}_0|\lambda_0-\lambda_1|,\,R_0^{-1}|R_1-R_0|,\,M_0^{-1}|M_1-M_0|\right)\leq
2^{2\tau+2}C_0C^{-2}\leq{1\over8}.
\]
Now let $C\geq2^{2\tau+5}C_0$, and assume that the Iterative
lemma can be applied repeatedly for $j\geq2$, with an input
parameter set
$\{\kappa_{j-1},r_{j-1},\rho_{j-1},T_{j-1},\sigma_{j-1},\lambda_{j-1},R_{j-1},M_{j-1},\mu_{j-1},\nu_{j-1}
\}$ and an analyticity loss $\delta_j$,
resulting in a transformation $\Psi_j=\Psi_j(\ar_j,S_j)$ and an output parameter
set $\{\kappa_j,r_j,T_j,\rho_j,\sigma_j,\lambda_j,R_j,M_j,\mu_j,\nu_j \}$, such
that
\[
\begin{array}{llllll}
\nu_j&\geq & (2M)^{-1} 2^{-j(\tau+1)}\delta^{\tau+1},\\
\mu_j&\leq&2^{-j(2\tau+4)}\mu , \\
m_j&=&\sup\left(\lambda^{-1}_{j-1}|\lambda_{j}-\lambda_{j-1}|,\,
R_{j-1}^{-1}|R_j-R_{j-1}|,\,M_{j-1}^{-1}|M_j-M_{j-1}|\right)
&\leq&{1\over8}.
\end{array}
\]
This would imply that the sequence $\{\mu_j\nu^{-1}_j\}_{j>1}$ vanishes
geometrically.

Suppose, the above assumption is true for $l=1,...,j-1$ (the case
$j=1$ has been checked). Then on the $j$th application of the
lemma, one can let $\varsigma_j=2^{-j\tau}\varsigma$, hence
$\eta_j\geq 2^{-j(\tau+1)-2}\eta$, because the inductive
assumption implies that $\lambda_{j-1}>{\lambda\over2}$,
$R_{j-1}>{R\over2}$, $M_{j-1}<2M$. The condition (\ref{small})
where each parameter involved has been endowed with an index $j$
is satisfied; in fact, the right hand side of it majorates a
vanishing geometric sequence with a ratio $2^{2\tau+3}$, whereas
the left hand side is majorated by a vanishing geometric sequence
with a ratio $2^{2\tau+4}$ (by the induction assumption). Another
application of the Iterative lemma yields
\[
\mu_{j+1}\leq 2^{(2\tau+2)+5} C_0^2 C^{-2}\mu_j\leq
2^{-(2\tau+4)}\mu_j.
\]
The fact that $m_{j}<{1\over8}$ is easy to verify, similar to the case $j=1$.
This justifies having $2M$ in the above assumption about $\nu_j$ and completes the
prof of
the induction assumption.

Now the statement of Theorem \ref{mt} and its Parameter statement follow from
chasing through the iterative scheme (\ref{iter}) and the estimates of the
Parameter statement of  Lemma \ref{il}. The existence of the limits
$\ar=\ar_1\circ \ar_2\circ\ldots$ and $S=\sum_{j=1}^\infty S_j$ follows from
the fast convergence of the estimates for their norms and completeness of the
spaces $\bsprime,\,\bszprime$. Finally, the estimates (\ref{transf})
pretty much reproduce the corresponding estimates of the Iterative lemma.
Indeed, due to the geometric convergence of the series $\sum_{j\geq1}S_j$ and
the composition
$\ar_1(\hat{\br}_1)\circ\ar_2(\hat{\br}_2)\circ\ldots$, it suffices to estimate the norms of
$S_1$ and $\hat{\br}_1$ only. $\Box$

\newpage
\bibliographystyle{unsrt}

\end{document}